\documentclass[11pt,english]{article}

\usepackage[T1]{fontenc}
\usepackage[latin1]{inputenc}
\usepackage{babel}

\usepackage{setspace}
\usepackage{indentfirst}

\makeatletter

\usepackage{verbatim}

\makeatother

\usepackage{enumerate}
\usepackage{amsmath}
\usepackage{amssymb}
\usepackage{stmaryrd}
\usepackage{wasysym}
\usepackage{pifont} 
\usepackage{graphicx}
\usepackage[all]{xy}

\usepackage{amsthm}
\theoremstyle{plain}
\newtheorem{theo}{Theorem}[section]
\newtheorem{lem}[theo]{Lemma}

\newtheorem{prop}[theo]{Proposition}

\newtheorem{defn}[theo]{Definition}
\theoremstyle{definition}
\newtheorem{example}[theo]{Example}
\theoremstyle{remark}
\newtheorem{rmk}[theo]{Remark}

\newcommand{\re}{{\mathfrak R}{\mathfrak e}}
\newcommand{\im}{{\mathfrak I}{\mathfrak m}}
\renewcommand{\epsilon}{\varepsilon}

\newcommand{\schouten}[1]{ [\! [ #1  ] \! ]  }
\def\dar[#1]{\ar@<2pt>[#1]\ar@<-2pt>[#1]}

\newcommand{\toto}{\rightrightarrows}
\newcommand{\gm}{\Gamma}

\newcommand{\R}{{\mathbb R}}

\newcommand{\diff}{{\rm d}}

\newcommand{\be }{\begin{eqnarray*}}
\newcommand{\ee }{\end{eqnarray*}}

\title{From Lie groupoids to resolutions of singularities. Applications
  to symplectic and Poisson resolutions.}

\author{Camille Laurent-Gengoux  \\
 \\
D\'epartement de math\'ematiques \\
Universit\'e
de Poitiers \\
86962 Futuroscope-Chasseneuil, France \\
\texttt{laurent@math.univ-poitiers.fr} }

\begin{document}
\maketitle

\begin{abstract}
We use the techniques of integration of Poisson manifolds into symplectic Lie groupoids to build
symplectic resolutions (=desingularizations) of the closure of a symplectic leaf and characterise the resolutions obtained by this procedure.
 More generally, we show how Lie groupoids can be used to lift singularities, in particular when one imposes
a compatibility condition with an additional structure given by a multi-vector field.
\end{abstract}

\tableofcontents

\section{Introduction}

\subsection{Presentation}

In the literature, the idea of replacing a Poisson variety $X$
and its oddities by a symplectic manifold $Z$ that projects on $X $ through a Poisson map  appears in two a priori
different contexts, depending on what is singular:  the variety or the Poisson structure.
 In algebraic geometry, see \cite{Be2,Fu,GK},  given a variety $W$
  endowed with a Poisson structure  of maximal rank (= symplectic) at regular points,  
  Beauville has introduced
 the notion of {\em symplectic  resolution}, which consists  of a resolution $\Sigma \stackrel{\phi}{\to }W $ in the sense of Hironaka's Big Theorem
 with $\Sigma$ symplectic and $\phi$ a Poisson map.
 In differential geometry, see \cite{CDW,DeSilvaWeinstein,DZ}, a Poisson manifold $M$ is
 being replaced by a manifold $\Gamma $,  called  {\em symplectic  groupoid}, which has a (maybe
local)
 Lie groupoid structure over $M$, together with  some compatible
symplectic structure, the Poisson map being then simply either
the source map or the target map.


This article discusses the possibility to go from symplectic
groupoids to  symplectic resolutions, and conversely. To make a long story short,
 {\em our aim is to show how the
symplectic groupoid of a Poisson manifold can be used to desingularize
  the closure of a symplectic leaf}.
A secondary aim is to use Poisson groupoids (more generally Lie
groupoids endowed with multiplicative $k$-vector fields)
to find Poisson resolutions (more generally resolutions compatible with a $k$-vector
field) of the closure of an algebroid leaf.

We point out the most crucial differences that exist between  the theories of
symplectic resolutions and symplectic groupoids, 
and explain briefly how we avoid or unify them.

\begin{enumerate}
\item First of all, symplectic resolutions belong to the world of
algebraic geometry while the theory of  symplectic groupoids has been
  developed on real smooth manifolds mainly (but most of its results
extend to the holomorphic setting, see \cite{SX}). To avoid this
  difficulty, one possibility could be to rewrite
 the theory of symplectic groupoids in the language of algebraic
  geometry. We make the opposite choice, and decide work  inside
   the world of differential geometry. 


\item More precisely, we mimic, within differential geometry, the definition
 of symplectic resolutions introduced by Beauville.
The object that we are going to try to desingularize is the closure $\overline{\mathcal S} $ of a locally
closed symplectic leaf ${\mathcal S} $ of a Poisson manifold $(M,\pi)$. This closure
$\overline{\mathcal S} $ behaves precisely as in the algebraic case since
 \begin{enumerate}
\item regular points (= points in  ${\mathcal S} $) form a dense open
  subset of $\overline{\mathcal S}$, (as the regular part of an algebraic variety does) and
\item since  ${\mathcal S} $ is a symplectic leaf, the
 restriction to $\overline{\mathcal S} $ of the Poisson structure is symplectic at regular points (= points in  ${\mathcal S} $).
\end{enumerate}
\item Symplectic resolutions and symplectic groupoids behave differently with respect to dimensions.
 On the one hand,
 a symplectic resolution of a
  Poisson variety of dimension $k$ has the same
  dimension $k$. On the other hand, the symplectic groupoid procedure
  doubles the dimension, id. est, it
  starts from a Poisson manifold $(M,\pi)$ of dimension $n$ and builds a symplectic Lie groupoid of dimension $2n$.
  In particular, the symplectic groupoid itself can by no way be itself the symplectic resolution. 
  But the main idea of this paper is that the symplectic groupoid $\Gamma \toto M $ can give a symplectic resolution  of
  the  closure $\overline{\mathcal S}$ of a symplectic leaf
 $ {\mathcal S}$  by going through the following steps:
 \begin{enumerate}
 \item we choose carefully some
  submanifold $L$ of $M $, included into $\overline{{\mathcal S}}$, and
  whose intersection with ${\mathcal S} $ is Lagrangian in ${\mathcal
  S} $, and
\item we apply a procedure called symplectic reduction to the
   submanifold $ \Gamma_L = s^{-1}(L)$ (where $s: \Gamma \to
  M$ is the source map) which is coisotropic in the symplectic
manifold $\Gamma $. This symplectic reduction  reduces
  the dimension, and
  involves in this case a strong Lie groupoid machinery, and
\item we obtain (under some conditions)
 a symplectic manifold that
 gives a symplectic resolution of $\overline{\mathcal S}$.
The Poisson map onto $\overline{\mathcal S} $
is induced by the target map of the groupoid $\Gamma \toto M $.  \end{enumerate}
\end{enumerate}

Let us point out two advantages of this method. First, it unifies
two theories of desingularization. More precisely,
$\overline{\mathcal S}$ may be ``singular'' in two ways:  it can be
a singular variety (whatever it means in the smooth/holomorphic
context), in this case, what is ``desingularized'' is a singularity
of the variety. But $\overline{\mathcal S}$ can very well be a
smooth manifold of dimension $2l$, but the Poisson structure is
singular (= of rank $< 2l $) at singular points (= points in $
\overline{\mathcal S} \backslash {\mathcal S}$). In this case, what
is ``desingularized'' is a singularity of the Poisson structure.

Second, it  generalises to other types of desingularizations
compatible with some additional geometrical structure given by a multi-vector field:
 Poisson resolutions in particular, but also  contact resolutions
 and twisted symplectic resolutions (to be studied in an other work), to mention just the main ones.
 In order to reach that level of generality, we shall work, as much as we can, in the very
general case of a Lie groupoid endowed with a  multiplicative $k$-vector field. 


The organisation of the paper is as follows. In Section \ref{sec:R1}, we ask the following {\em  question}:  ``Given
a Lie groupoid, how can one desingularize the closure $\overline{\mathcal S} $ of an  algebroid leaf ${\mathcal S}$~?''
and our (partial)  {\em answer} is: ``With the help of an algebroid
crossing (see Definition \ref{def:algcro}), under some integrability assumptions''.

In Section \ref{sec:R2}, we recall from \cite{ILX} the following {\em fact}: ``a multiplicative $k$-vector field on a
 groupoid induces a $k$-vector field $\pi_M $ on the algebroid leaves''.
Then we raise the {\em question}: ``Given an algebroid and a multiplicative $k$-vector field
 on the corresponding groupoid (with $k  \geq 2$), how can one desingularize $\overline{\mathcal S} $ in a
way that is compatible with this  $k$-vector field ?'' and we suggest the following {\em answer}: ``With the help of a 
 algebroid crossing coisotropic with respect to $\pi_M $,  under some integrability condition''.

In Section \ref{sec:R3},  we raise the {\em question}: ``How can we construct a symplectic resolution
of the closure $\overline{\mathcal S} $ of a symplectic leaf  ${\mathcal S} $ of a Poisson manifold ?''
and we propose the following {\em answer}: `` Our previous results, applied
to the special case of a Lagrangian crossing (see Definition \ref{def:Lag_cross}) give automatically a symplectic resolution''.

In Section \ref{sec:char}, we give a characterisation of symplectic resolutions of the previous forms in terms of compatibility
with respect to a Lagrangian crossing.

We then present in Section \ref{sec:Ex} examples of such resolutions. 
Section  \ref{sec:Ex_2D} presents a trivial example that  illustrates in a very clear way
our results: we lift, in the world of real geometry, the singularity  at the origin of the real 
Poisson bracket on ${\mathbb R}^2 $ given by 
    $$ \{x,y\} = x^2 + y^2 .$$
In Section  \ref{sec:Ex_Springer}, the celebrated Springer resolution is rediscovered as a particular case of the previous constructions.
Note that the second of these examples lifts a singularity of the variety, while the first one lifts a singularity of the Poisson structure.
We then present an example of Poisson resolution, namely the Grothendieck resolution endowed with its Evens-Lu Poisson structure.
 We finish with a discussion of a famous symplectic resolution:
the minimal resolution of ${\mathbb C}^2 / G$ with $G = \frac{{\mathbb Z}}{l {\mathbb Z}} $.

\subsection{Acknowledgements}

First, I would like to thank Jiang-Hua Lu for her invitation to
the Hong-Kong University in June 2006. Many ideas of this paper have
been discussed there, and have very close relations with her joined work
with S. Evens (see \cite{LuEvens}). 

Second, I am very grateful with
the Erwin Schr\"odinger Institute (ESI) in Vienna
where the first draft of the paper was written.

I gratefully thank Ariane Le Blanc, Rui Fernandes, Baohua Fu, Yvette
Kosmann-Schwarzbach, Thierry Lambre, Mathieu Stienon, Patrice
Tauvel, Pol Vanhaecke, Friedrich Wagemann, Ping Xu, Rupert Yu and Marco Zambon
 for useful discussions.  I have  the hope (unlikely to be satisfied) that I do
not forget any one. The advises and comments of Baohua Fu
 were especially useful.
All the groupoid machinery used here were taught
to me by Ping Xu, and I have to thank him again for the huge amount
of knowledge that I learnt from him in the past four years.

Also, I have to thank  Claire Bousseau for her bright "M\'emoire de
Master 2" entitled ``D�singularisation symplectique et quotient de
${\mathbb C}^2 $ par un groupe fini de $SL_2({\mathbb C}) $'' that
clarifies (and made it possible for me to understand) many a basic
fact about symplectic resolutions. Last, Eva Miranda corrected several typos.

\subsection{Notations and basic facts about Lie groupoids}
\label{Sec:Notations}

For self-containess of this paper, we recall
several facts about Poisson manifolds, Lie algebroids, and Lie groupoids.
We recall in particular how to adapt to the holomorphic setting
the theories of symplectic groupoids and holomorphic Poisson
structures, an adaptation  systematically studied in \cite{SX}.

\paragraph{ Complex and real manifolds.}
We try to state and prove results which are valid
in both complex and real differential geometries.
By a {\em manifold}, we mean {\em a complex manifold or a real
  (smooth) manifold}. Then words like functions, vector fields,
vector bundles, sections
 should be understood as being holomorphic or
smooth, depending  on the context.
Moreover, in  complex geometry, we have
to work with local functions, local sections or local vector fields
rather than their global counterparts.
 Given a manifold $M$ and a vector bundle $A \to M $,
we denote simply by ${\mathcal F}(M)$ and $\Gamma(A \to M) $ the sheaf of local
functions and local sections respectively.
Sections over an open subset $V$ are denoted by $\Gamma(A_{|_V} \to V)$.
When we write identities of the form $f=g$ for local functions
$f$ and $g$, it should be understood that the identity takes place on
the intersection of their domain of definition. 
When we say that a local function $f$ vanishes on a submanifold $N$,
we mean of course that $f$ vanishes on the intersection of $N$
with the domain of definition of $f$. 
Similarly, terms of the form $X[f_1\cdots,f_k] $, where $X$ is
a a local $k$-vector field, and $f_1,\cdots,f_k $ are local functions
are to be considered only where it is defined. 
All these slight abuses of terminology shall be systematically omitted in the text.

All our manifolds are supposed to be Hausdorff. Since this last point can be
ambiguous while speaking about Lie groupoid,
we sometime say Hausdorff Lie groupoid to emphasise on this assumption.

\paragraph{ Fibered product, definition and convention.}
Given a triple  $M_1,M_2, N $ of manifolds and maps $\phi_1: M_1 \to N$
and $\phi_2: M_2 \to N$, we denote by $ M_1 \times_{\phi_1, N,\phi_2}
M_2 $ the fibered product:
$$ M_1 \times_{\phi_1, N,\phi_2} M_2
:=\{(m_1,m_2)\in M_1 \times M_2 |  \phi_1(m_1) = \phi_2(m_2)  \}.$$
This fibered product is itself a manifold as  soon as one of the
maps $\phi_1,\phi_2 $ is a submersion.


\paragraph{  Multi-vector fields, definition and convention.} A $k$-vector
field $\pi $ is a section of $\wedge^k TM \to M $ (the wedge product
being over $ {\mathbb R}$ or ${\mathbb C}$ depending on the
context). Throughout this paper, we will almost never have to consider
two different $k$-vector fields defined on the same space. It is therefore very
convenient for us to always denote by $\pi_M $ the  $k$-vector field
over $M$ that we consider.
Recall that a $k$-vector field $\pi_M $ defines a skew-symmetric $k
$-derivation of the algebra $ {\mathcal F}(U)$ of functions on an
open subset $U \subset M$  with the help of the skew-symmetric
$k$-derivation
  $$  F_1,\cdots,F_k  \to( x \to (\pi_M)_{|_x}  [ \diff_x F_1 ,\cdots, \diff_x F_k ] )$$
We denote by $ \{ \cdot,\cdots, \cdot  \}_M$ the previous skew-symmetric $k$-derivation.
We denote by $\schouten{\cdot,\cdot}_{TM} $ the Schouten bracket of multivector fields.

\paragraph{ Holomorphic/real Lie algebroids, definition.}
With some care, we define smooth and holomorphic algebroids.

\begin{defn}\label{def:algebroids}
 A {\em real Lie algebroid} $A$ over a smooth manifold $M$ is a smooth
vector bundle $A\to M$ together with a smooth collection, for each
$x \in M $, of linear maps   $\rho:A_x\to T_xM$ called the anchor,
and an $\R$-linear skew-symmetric bracket $[\cdot,\cdot]:\Gamma
A\times\Gamma A\to\Gamma A$ on the space of local smooth
 sections which satisfies the Jacobi identity, and
such that

\begin{itemize}
\item[(1)] $[X,fY]\,=\,f[X,Y]+\rho(X)(f)Y$
\item[(2)] $\rho([X,Y])\,=\,[\rho(X),\rho(Y)]$
\end{itemize}
We denote by $(A \to M, \rho,[\cdot,\cdot]) $ a Lie algebroid.

\smallskip

A {\em holomorphic Lie algebroid} is a smooth Lie algebroid $(A \to
M,\rho,[\cdot,\cdot ] )$ over a complex manifold $M$, where $A \to M
$ is a complex vector bundle, such that the anchor map is  holomorphic,
 and the bracket $[\cdot,\cdot] $ restricts  to a ${\mathbb C}$-linear Lie algebra structure on 
 holomorphic  sections over any open subset $U \subset M$.
\end{defn}

We warn the reader that holomorphic Lie algebroid should not be 
confused  with Complex Lie algebroid (CLA) in the sense of
\cite{Weinstein1}.


\paragraph{ The foliation of a Lie algebroid.} The distribution
$\coprod_{x \in M}  \rho (A_x)$ is integrable, though its rank is
not constant. Moreover, the leaves, called {\em algebroid leaves}, are smooth immersed
submanifolds, which are immersed complex submanifolds if the Lie
algebroid is a holomorphic Lie algebroid, see \cite{SX}.

%
%

\paragraph{ $k$-differential.} Recall from \cite{X_BV} that $\Gamma(\wedge^\bullet A \to M) $ is a (sheaf of) Gerstenhaber
 algebra  when equipped with the wedge product and the unique natural
 natural extension $\schouten{\cdot,\cdot}_A $ 
of the bracket of sections of $A \to M  $ (extension which is again called Schouten bracket).
 A {\em $k$-differential}  (see
 Section 2 in \cite{ILX}) is a linear operator $\delta :\gm (\wedge
^\bullet A)\to \gm (\wedge ^{\bullet +k-1}A)$ satisfying
\begin{equation}\label{k-dif-prop}
\left\{ \begin{array}{l}
\delta (P\wedge Q)=(\delta P)\wedge Q+(-1)^{p(k+1)}P\wedge \delta
Q, \\ \delta \left( \schouten{ P, Q}_A  \right) =\schouten{ \delta P, Q}_A +(-1)^{(p+1)(k+1)}
\schouten{ P, \delta Q }_A .
\end{array}\right.
\end{equation}

\paragraph{ Holomorphic/real Lie groupoids, definition.}
There is no issue in defining Lie groupoids in both the real on
the complex case, see \cite{SX}.
 We could say that a holomorphic/real Lie groupoid is a
small category where objects, arrows and operations
are holomorphic/smooth. We prefer
to give a long but down-to-the-earth definition.

\begin{defn}
A {\em holomorphic/real Lie groupoid} $\Gamma \toto M$ consists of two complex/real
manifolds $\Gamma$ and $M$, together with two holomorphic/smooth
surjective submersions $s,t:\Gamma\to M$, called the source and the
target maps, and a holomorphic/smooth inclusion $\epsilon: M\to
\Gamma$, which admits
a group law, that it to say, such that there exists {\em (i)} a
holomorphic/smooth inverse map $\Gamma \to \Gamma $ denoted $\gamma
\to \gamma^{-1} $  permuting the source and target maps and {\em
(ii)} a holomorphic/smooth map
 $\Gamma \times_{t,M,s} \Gamma  \to \Gamma$  denoted
$(\gamma_1,\gamma_2) \to \gamma_1 \cdot \gamma_2 $
satisfying
$$   \begin{array}{rcll} s(\gamma_1 \cdot \gamma_2)&= & s(\gamma_1) &
  \forall  (\gamma_1,\gamma_2) \in \Gamma \times_{t,M,s} \Gamma \\
t(\gamma_1 \cdot \gamma_2)&= & t(\gamma_2) &
  \forall  (\gamma_1,\gamma_2) \in \Gamma \times_{t,M,s} \Gamma \\
(\gamma_1 \cdot \gamma_2) \cdot \gamma_3 & = & \gamma_1 \cdot
  (\gamma_2 \cdot \gamma_3 ) & \forall
  ( \gamma_1,\gamma_2,\gamma_3 )
  \in \Gamma \times_{t,M,s} \Gamma \times_{t,M,s}  \Gamma \\
   \epsilon ( s(\gamma)) \cdot \gamma &=&\gamma & \forall \gamma \in \Gamma \\
   \gamma \cdot   \epsilon ( t(\gamma))          &=&\gamma  &
 \forall \gamma \in \Gamma   \\
\gamma \cdot \gamma^{-1} &=& \epsilon(s(\gamma))& \forall \gamma \in \Gamma\\
\gamma^{-1} \cdot \gamma &=& \epsilon(t(\gamma))& \forall \gamma \in \Gamma \\
   (\gamma_1 \cdot \gamma_2)^{-1} & = & \gamma_2^{-1} \cdot
  \gamma_1^{-1} &\forall (\gamma_1 ,\gamma_2 \in \Gamma )  \in \Gamma
  \times_{t,M,s} \Gamma \\
 \end{array} $$
\end{defn}

Of course, the previous list of requirements is redundant, and
several of the previous axioms could be erased.

 To any holomorphic/real Lie groupoid, there is an holomorphic/real algebroid
associated with, (see \cite{McKenzie} Section 3.5 for the real case,
the complex case is similar). 
This Lie algebroid is defined by the vector bundle $A :=\coprod_{x \in M} A_x$
where $A_x  = {\rm ker} (\diff_x s)\subset T_{\epsilon (x)} \Gamma $.
%
%
%
%
%
%
%
The converse is not true: a Lie algebroid may not be the Lie
algebroid of a Lie groupoid, see \cite{CF}. When such a Lie
groupoid exists, we say that the Lie groupoid is question
{\em integrates} the algebroid. 
%
%
%
Recall from \cite{SX} that, when a holomorphic Lie algebroid,
when considered as a real Lie algebroid by forgetting the complex
structure, integrates to a smooth
source-simply-connected Lie groupoid $\Gamma \toto M $, then 
$\Gamma \toto M $ inherits an unique natural complex structure that turns it
into a holomorphic Lie groupoid. 

We use the following notations. For any $L, L' \subset M$, one introduces $\Gamma_L := s^{-1}(L) $,
$\Gamma^{L'} := t^{-1}(L') $ and $\Gamma_L^{L'}= \Gamma_L \cap \Gamma^{L'}$.
Given a point $x  \in M$, we use the shorthands $\Gamma_x , \Gamma^x$ for
$\Gamma_{\{x\}} $ and $ \Gamma^{\{x\}}$ respectively.

\paragraph{ Left- and right-invariant vector fields.}
When the Lie algebroid integrates to a Lie groupoid,
sections of $A \to M $ are in one-to-one correspondence with 
right-invariant (resp. left-invariant) vector fields on $ \Gamma \toto
M$.
Given a section $X \in \Gamma (A) $, we denote by $
\overrightarrow{X} $ (resp. $ \overleftarrow{X} $) the right-invariant 
(resp. left) vector fields corresponding to it.

We extend this convention to sections of $\wedge^\bullet A $. Namely, following \cite{ILX},
we define,  for all $k \geq 1$ and all (maybe local) 
sections $e_1,\dots,e_k \in \Gamma (A \to M)$, a right-invariant
$k$-vector field on $\Gamma \toto M$,
by $\overrightarrow{e_1 \wedge\dots \wedge e_k}  = \overrightarrow{e_1}
 \wedge \dots \wedge \overrightarrow{e_k}  $   (resp.
$\overleftarrow{e_1 \wedge\dots \wedge e_k}  = \overleftarrow{e_1}
 \wedge \dots \wedge \overleftarrow{e_k}  $  ) and extend this
 correspondence by multilinearity. For $k=0$, sections of
$\wedge^0 A $ are simply functions: we define then $
\overrightarrow{f}=s^* f $ (resp. $
 \overleftarrow{X} =t^* f$). 

The map from $\Gamma(\wedge^{\bullet} A \to M) \to {\mathcal X}^{\bullet} (\Gamma)
$ given by $X \to    \overrightarrow{X} $ is a morphism of
Gerstenhaber algebras,
where ${\mathcal X}^{\bullet} (\Gamma) $ stands for the Gerstenhaber
algebra of multivector fields.
%
%
When dealing with left-invariant vector fields, the same result
holds up to a sign, namely the map defined by
$X \to   (-1)^k \overleftarrow{X} $ for all $X \in \Gamma(\wedge^{k} A \to M) $ 
is a morphism of
Gerstenhaber algebras.

%

\paragraph{ Modules of a Lie groupoid.} 
A left action of a Lie groupoid $\Gamma
\toto M $ on a pair  $(X,\phi)$, with $X$ a manifold
and 
$\phi:X \to M $, is a holomorphic/smooth map from $\Gamma
\times_{t,M,\phi} X $ (which is a manifold because the target
map is a submersion) to $X$, called the {\em action map} and denoted
by $(\gamma,x) \to \gamma \cdot x $, which satisfies the following
axioms (see \cite{McKenzie} Section 1.6):
  $$  \left\{ \begin{array}{rcll}   \gamma \cdot x & \in &  \phi^{-1}(s(x)) & \forall (\gamma,x) \in \Gamma \times_{x,M,\phi} X  \\  \epsilon (\phi(x)) \cdot x  &=& x & \forall x
  \in X \\  \gamma_1 \cdot (\gamma_2 \cdot x) &= & (\gamma_1 \cdot \gamma_2) \cdot x
   & \forall (\gamma_1,\gamma_2,x) \in \Gamma \times_{t,M,s} \Gamma
  \times_{t,M,\phi}  X  \end{array} \right. $$
Right action are defined in the same way. We say then that $(X,\phi)$
is a {\em (left or right) $\Gamma$-module}.
Given a left (resp. right) action of $\Gamma \toto M$ on $X \to M$, we call {\em quotient
space} and denote by $\Gamma \backslash X  $ (resp. $X/ \Gamma $) the space  of orbits of
the action, id est, the space $X/ \sim $  where $\sim $ is the
equivalence relation that identifies  $ x$ and $\gamma \cdot x $ (resp.  $x \cdot \gamma$) for
any  $(\gamma,x) \in \Gamma \times_{t,M,\phi} X  $ (resp. $X \times_{\phi,M,s} \Gamma $).

A groupoid $ \Gamma \toto M$ acts on 
$ ( M, {\rm Id}) $ by
 $$ \gamma \cdot m = t(\gamma)  \, \, \, \mbox{if $s(\gamma)=m $} .$$
The following result concatenates several results in \cite{McKenzie}.

\begin{lem} \label{lem:leaf0} Let $ \Gamma \toto M $ be a
  source-connected Lie groupoid integrating
$A \to M$. Any two points in $M$ are is the same Lie  algebroid leaf if and only
if they are in the same orbit ${\mathcal S} $ of the action of $ \Gamma \toto M$ on
$M$ (i.e. if and only if they are source and target of some $\gamma 
\in \Gamma$). Moreover, the map $\gamma \to (s(\gamma), t(\gamma)) $ 
restricts to a submersion from $\Gamma_{\mathcal S} = \Gamma_{\mathcal
  S}^{\mathcal S}= \Gamma^{\mathcal S}$ to ${\mathcal S} \times {\mathcal S} $.
\end{lem}
 We
say that a subset $X \subset M $  is {\em $\Gamma $-connected} if and only
if, for any two $x,y \in X$ there exists a finite sequence
$x=x_0,x_1,\cdots, x_l =y $ such that, for all $i \in \{0,\cdots,
n_1\}$, either $x_i, x_{i+1}$ are the same connected component of $X$, or
$x_i$ and $x_{i+1} $ are the source and target of some $\gamma \in \Gamma
$.

\paragraph{Holomorphic Poisson manifold and integration.}
We recall that for any Poisson manifold $(M,\pi_M)$,
the cotangent bundle $T^*M  \to M$ is endowed with a natural
Lie algebroid structure $(T^*M,\pi_M^\# , [\cdot,\cdot]_{\pi_M})$.

 A  {\em real/holomorphic Poisson manifold $(M,\pi_M) $} is a complex manifold
$ M$ endowed with a smooth/holomorphic bivector field $\pi_M$ such that $\schouten{\pi_M,\pi_M}_{TM}=0$. 
Recall from \cite{SX} that the real part $\re (\pi_M)$ and
the imaginary part $\im (\pi_M)$ of $\pi_M $ are compatible smooth Poisson
structures.

We refer to \cite{CDW} for an overview of the theory of integration of
Poisson manifolds to symplectic groupoid in the real case, and we give
a brief overview.
A {\em symplectic Lie groupoid} is a pair $(\Gamma \toto
 M,\omega_\Gamma )$ where
 $\Gamma \toto M $ is a Lie groupoid and
$   \omega_\Gamma $ a symplectic structure on $\Gamma $
which is compatible with respect to the multiplication.
 A symplectic Lie groupoid $(\Gamma \toto M, \omega_{\Gamma}) $
is said to {\em integrate} $(M,\pi_M) $  if $s_* \pi_\Gamma =\pi_M $.
Any source-simply-connected Lie groupoid
$\Gamma \toto M$ integrating the algebroid $(T^* M,\pi_M^{\#},[\cdot,\cdot]^{\pi_M},
) $  admits a unique symplectic structure
integrating $(M,\pi_M)$.

Assume that there exists a source-simply-connected symplectic Lie groupoid $\Gamma
\toto M $ integrating the algebroid  $(T^*M,(\re(\pi_M))^{\#} , [\cdot,\cdot]_{\re(\pi_M)})$. 
According to \cite{SX}, the Lie groupoid $\Gamma \toto M$ is indeed a holomorphic symplectic Lie
groupoid with respect to holomorphic symplectic form $\omega$ whose real part is  $ \frac{1}{4}\omega_R$.

\section{Resolution of the closure of an algebroid leaf.}
\label{sec:R1}

\subsection{Definition of a resolution in real and complex geometries}
 \label{sec:R1_D1}

In the context of algebraic geometry, a {\em resolution} of a variety
$W$ is pair $(Z,\phi)$ where $Z$ is a smooth (= without singularities) variety 
and $\phi:Z \to W$ is a proper regular birational map onto $W$.
In particular,  the restriction of $\phi   $ to $\phi^{-1}(W_{reg }) \to
W_{reg}$ is biregular (where $W_{reg}$ stands for the regular part
of $W$).

 One difficulty arises when we try to reformulate this definition in the context
of complex or real geometry:  what kind of ``singular'' varieties
should we consider ? We avoid this difficulty by restricting ourself
to the following situation. All the singular smooth/complex
``varieties'' that we are going to study are of the form
$\bar{{\mathcal S}}$ where
 ${\mathcal S}$ is a (not closed in general)
locally closed ($=$ embedded) submanifold
 of a manifold $M$, and where the closure is with
respect to the usual topology of $M$.


In this context, by a resolution, we mean the following:

\begin{defn} \label{def:desing}
Let $\bar{{\mathcal S}}$ be a the closure of an locally closed (=
embedded) submanifold ${\mathcal S} $ of a complex/real manifold $ M$.
 A {\em resolution of
$\bar{{\mathcal S}} $} is a pair $(Z,\phi)$ where $Z$ is a
complex/real manifold and $\phi:Z \to M$ is a holomorphic/smooth map such that
\begin{enumerate}
\item $\phi(Z) = {\bar{\mathcal S}}$,
\item $\phi^{-1}({\mathcal S}) $ is dense in $Z$,
\item the restricted map $\phi :\phi^{-1}({\mathcal S}) \to {\mathcal S}$ is an biholomorphism/diffeomorphism.
\end{enumerate}
When $\phi:\phi^{-1} ({\mathcal S})  \to {\mathcal S}$ is only an
\'etale map (i.e. a surjective local biholomorphism/diffeomorphism),
then we speak of an {\em \'etale resolution of $\bar{\mathcal S} $}.
When $\phi:\phi^{-1} ({\mathcal S})  \to {\mathcal S}$ is only a
covering (i.e. $\phi^{-1} ({\mathcal S})$ is  a connected set
and $\phi:\phi^{-1}({\mathcal S})  \to {\mathcal S} $ is an \'etale
map), then we speak of a {\em covering resolution of $\bar{\mathcal S} $}.

For \'etale or covering resolutions, by saying that the typical
fiber is isomorphic to some discrete set $Z $, we mean that
$\phi^{-1}(x ) \simeq Z $ for all $ x \in {\mathcal S}$.%
\end{defn}

\begin{rmk}
Let $M$ be a nonsingular variety over ${\mathbb C} $ and $W \subset N $ an irreducible subvariety.
Then the regular part $W_{reg} $ of $W$ is a locally closed complex submanifold of $M$ and 
$\overline{W_{reg}} = W$ (where the closure is with respect to the usual topology of $M$).
A resolution (in the sense of algebraic geometry) of $W $ is also a holomorphic resolution
(in the sense of Definition \ref{def:desing}) of $\overline{W_{reg}} $.
\end{rmk}

\begin{rmk}
Since ${{\mathcal S}} $  is dense in $ \overline{{\mathcal S}} $,
Condition 1 in Definition \ref{def:desing} implies Condition 2 in the
complex case. 
\end{rmk}

\begin{rmk}
It deserves to be noted that $\bar{{\mathcal S}}$ may be itself a manifold with no
singularity: take, for instance, $M= {\mathbb C}$ and ${\mathcal
S}={\mathbb C}^*$. As a consequence,  in the case where
$\bar{{\mathcal S}}$ happens to be an affine subvariety of ${\mathbb
C}^n$, then ${\mathcal S}$, which is always included in
$(\bar{{\mathcal S}})_{reg}$, may be strictly included into
$(\bar{{\mathcal S}})_{reg}$. So that, in the particular case where $\bar{{\mathcal
S}}$ is itself an affine variety,
 a holomorphic resolution $(Z,\phi) $ of $\bar{{\mathcal S}}$ is the sense of 
Definition \ref{def:desing} may not be a resolution in the sense of
algebraic geometry,
even when $(Z,\phi)$ is in the category of algebraic varieties.
\end{rmk}


A {\em morphism between two \'etale resolutions $(Z_i,\phi_i) $,
$i= 1 , 2 $}
is a holomorphic/smooth
map $\Psi: Z_1 \to Z_2 $ such that $\phi_2 \circ \Psi = \phi_1  $, i.e.
such that the following diagram commutes

$$ \xymatrix{
Z_1 \ar[r]^{\Psi} \ar@<0.5ex>[dr]_{\phi_1}&   Z_2 \ar@<0.5ex>[d]^{\phi_2} \\
 &\bar{\mathcal S}    } $$

\subsection{Lie groupoids and resolutions of  the closure of an algebroid leaf.}
\label{sec:R1G}

We explain in this section how to build a resolution
 of the closure $\bar{\mathcal S} $
of a locally closed (= embedded) submanifold ${\mathcal S} $ of a
manifold $M$ which happens to be an algebroid leaf of an integrable algebroid.
%

%

\begin{lem}\label{lem:leaf}
Let ${\mathcal S} $ be a locally closed  leaf of a Lie algebroid $(A
\to M,\rho, [\cdot,\cdot])  $. Then $\bar{{\mathcal S}}$ is a
disjoint union of algebroid leaves.
%
\end{lem}
\begin{proof}
Though one could derive this result directly from general considerations
about integrable distributions, we give an easy proof that uses local integration.
Let $\Gamma \toto M$ be a source-connected local Lie groupoid integrating 
the algebroid $ (A \to M, \rho, [\cdot,\cdot])$.
Choose a point $m \in   \bar{{\mathcal S}}$ and an element $\gamma \in
\Gamma_m$.
In any neighbourhood $W$ of $\gamma \in \Gamma  $,
there exists, since the source map is a submersion, at least one
element $\gamma' \in W $
which is mapped by the source map to an element in ${\mathcal S}$.
The target  $t(\gamma')$ of such an element belongs to
$ {\mathcal S} $. Since $W$ can be chosen arbitrary small,
  $ t(\gamma)$ belongs to $\bar{{\mathcal S}}$. 
By Lemma \ref{lem:leaf0}, $\bar{{\mathcal S}} $ needs to be an union
of algebroid leaves. 
\end{proof} 
%
%
 
%
%

\begin{rmk} 
In view of Lemma \ref{lem:leaf}, we may suggest that the notion of stratified space 
would be more relevant here and  $\bar{{\mathcal S}} $
should be called a ``strate'', however, it seems to be customary to speak about ``algebroids leaves'' 
(and in particular of  ``symplectic foliation'' in the case of Poisson
manifold) so we  adhere to this traditional vocabulary.
\end{rmk} 

The construction of \'etale resolutions of $ \bar{\mathcal S}$
requires the following object.

\begin{defn}\label{def:algcro} Let  $(A
\to M,\rho,[\cdot,\cdot]) $ be a Lie algebroid and ${\mathcal S} $
a locally closed (= embedded) algebroid leaf.
We say that a submanifold $L $ of $M$ is an {\em algebroid crossing of
$\bar{\mathcal S} $} if and only if the following conditions are
satisfied
\begin{enumerate}
\item $L$ is a submanifold of $M $,
\item $L \subset \bar{\mathcal S} $,
\item $L$ intersects all the algebroid leaves contained in $\bar{\mathcal S} $,
\item  $L \cap {\mathcal S} $ is dense in $L $, 
\item the vector bundle  $ B \to L \cap {\mathcal S}$ defined over $L \cap {\mathcal S} $ 
 by  $$  B_x = \rho^{-1}(T_xL) \, \, \, \, \, \,  \forall x \in  L \cap {\mathcal S} $$
extends to a (holomorphic/smooth) vector bundle over $L$, denoted again by
 $B \to L$ and called
{\em  normalisation of the algebroid crossing $L$}. 
\end{enumerate}
\end{defn}

\begin{rmk}
Note that for a manifold $L$ that satisfies Conditions (1)-(4), the
normalisation, if any, is unique. note also that Condition 52) is a consequence of Condition (4).
 \end{rmk}

\begin{example} \label{ex:lightcone}
The reader may ask how  it can be that $\overline{\mathcal S} $ is singular 
while $L$ is smooth, especially taking under account that Definition \ref{def:algcro}(3)
forces $L$ to go through the singularities of $\overline{\mathcal S} $.
These conditions are, however, perfectly compatible.
Assume, for instance, that $M  $ is a vector space, and that ${\mathcal S} $
is a cone, minus the origin. Then $\overline{\mathcal S} = {\mathcal S} \cap \{0\}$ 
is the cone itself. Any straight line $L$ contained in that cone  is a smooth submanifold of $M$
and admits a non-empty intersection with both $ {\mathcal S}$ and  $ \{0\}$. 
Here is an example for which this situation occurs.
 The Lie group ${\rm SO}_n ({\mathbb C}) $ acts on  ${\mathbb C}^n$, and one can form the action 
 groupoid ${\rm SO}_n ({\mathbb C}) \times {\mathbb C}^n \toto 
{\mathbb C}^n $. 
The following set is an algebroid leaf of the Lie algebroid of the previous groupoid:
 $$ {\mathcal S} = \left\{ (z_1, \cdots, z_n)  \, \in \, {\mathbb C}^n - \{0\} \, \, | \,  \, z_1^2 + \cdots + z_n^2 =0 \right\}.$$
Any straight line $L$ through ${\mathcal S} $ is an algebroid crossing with
normalisation $ {\rm Lie}{(H)} \times L  \to L $, where $H$ is the stabiliser of $L$ in $ {\rm SO}_n ({\mathbb C}) $.  
 \end{example}

By density of $L \cap {\mathcal S} $ in $L $, 
the inclusion $\rho (B_x) \subset  T_x L$ holds for any $x \in L $.
In particular:

\begin{lem}\label{lem:subalg}
The normalisation $B \to L $ of an algebroid crossing $L$ of $\overline{\mathcal S} $,
admits a unique  algebroid structure such that the inclusion map into $A
\to M, $ is an algebroid morphism. 
(In other words, the normalisation $B \to L $ of an algebroid crossing
$L$ is a subalgebroid of  $(A
\to M,\rho,[\cdot,\cdot]) $ ).
\end{lem}

The following Proposition is the main result of this section.
It makes it possible to construct an \'etale resolution of $\overline{\mathcal S} $
out of an algebroid crossing, and helps us to decide whether this
\'etale resolution turns
to be a covering resolution or a resolution.

\begin{prop}\label{prop:desing+groupoid}
Let $(A \to M,\rho,[\cdot,\cdot]) $ be a Lie algebroid, ${\mathcal
S} $  a locally closed orbit of this algebroid, and $L$ an algebroid
 crossing of $\bar{\mathcal S}$ with normalisation $B \to L$. If
\begin{enumerate}
\item there exists a source-connected Hausdorff Lie groupoid $\Gamma
 \toto M $ integrating
the Lie algebroid $(A \to M,\rho,[\cdot,\cdot])  $ and
\item  there exists a  sub-Lie groupoid $R \toto L $ of $ \Gamma \toto M$,
  closed as  a subset of
  $\Gamma_L^L $, integrating
the subalgebroid $B \to L  $,
\end{enumerate}
then
 \begin{enumerate}
\item  $(Z(R),\phi)  $ is an \'etale resolution of $\bar{\mathcal S} $,
where  \begin{enumerate}
 \item $Z(R)=  R \backslash \Gamma_L$ and,
\item  $\phi: Z(R) \to M $
is the unique holomorphic/smooth map
such that  the following diagram commutes
\begin{equation}\label{eq:comdia1}
 \xymatrix{
\Gamma_L   \ar[r]^{p} \ar[rd]^{t}  & Z(R) \ar@<0.5ex>[d]^{\phi} \\
  & M}
\end{equation}
where $p:\Gamma_L \to Z(R)= R\backslash \Gamma_{L} $ is the natural projection.
\end{enumerate}
\item When  $L \cap {\mathcal S} $ is a $R$-connected set,
this \'etale resolution is a covering resolution with typical fiber
$\frac{\pi_0(I_x(\Gamma) ) }{\pi_0(I_x(R)) } $, where $x \in
{\mathcal S} $ is an arbitrary point, and $I_x(\Gamma) $ (resp.
$I_x(R) $) stands for the isotropy group of $ \Gamma \toto M$ (resp.
of $R\toto L$) at the point $x$.

\item This \'etale resolution is a resolution if
and only if $R$ contains $\Gamma_{L \cap {\mathcal S}}^{L
    \cap {\mathcal S}} $. In this case, we have
     $R=\overline{ \Gamma_{L \cap {\mathcal S}}^{L \cap {\mathcal S}}} \cap \Gamma_L^L$.

\item When $L \cap {\mathcal S} $ is a connected set
and $ R \toto L $ is a source-connected sub-Lie groupoid of $\Gamma
\toto M $, then the typical fiber is   $\frac{\pi_1( {\mathcal S} )
}{j(\pi_1( L \cap {\mathcal S}  ))} $, where $j$ is the map induced
at the fundamental group level by the inclusion of $L \cap {\mathcal
S} $ into ${\mathcal S} $.
\end{enumerate}
\end{prop}
\begin{example}
For the algebroid crossing given in Example \ref{ex:lightcone}, the sub-Lie
groupoid $R = H \times L \toto L $ satisfies the conditions of Proposition \ref{prop:desing+groupoid}(3).
Let us describe the obtained resolution. The manifold
$Z(R)$ is equal to the quotient of $ {\rm{ SO}}_n ({\mathbb C}) \times L $,
 with respect to the left $H$-action of $H$ given by 
$ h \cdot (O,l) = (hO, h(l)) $  for all $h \in h, O \in  {\rm{ SO}}_n ({\mathbb C}), l \in L$. 
The map $\phi: Z(R) \to \overline{{\mathcal S}} $ is given by $\phi([g,l]) = g^{-1} (l)$,
where $[g,l] $ stands for the class in $Z(R) $ of an element $(g,l) \in 
{\rm{ SO}}_n ({\mathbb C}) \times L $.
\end{example}
\begin{rmk}
If $L$ itself is a closed subset of $M$, then $R$ is closed as  a subset of
  $\Gamma_L^L $ if and only if it is a closed subset of $\Gamma$.
  \end{rmk}
\begin{proof}
1) Since $ L $ is a submanifold of $M$,
and the source map $s: \Gamma \to M$ is a surjective submersion,
 $\Gamma_L =s^{-1}(L) $ is a submanifold of $\Gamma$.
 The Lie groupoid  $R \toto L$ acts on $\Gamma_L $ by left multiplication.
This action is  free and the action map , i.e. the map,
$$ \begin{array}{rcl}  \mu:R \times_{t,L,s} \Gamma_L    & \to &
\Gamma_L  \times \Gamma_L \\ \mu (r,\gamma) &=  &(r \cdot \gamma, \gamma)  \end{array}
  $$
is proper. We need to prove this last claim precisely. Let $K$
be a compact subset of $\Gamma_L \times \Gamma_L$, and
$(r_n,\gamma_n)_{n \in {\mathbb N}} $ a sequence in $R
\times_{t,M,s} \Gamma_L$ with $\mu(r_n,\gamma_n) = (r_n \cdot
\gamma_n ,\gamma_n)\in K$. By compactness of $K$, one can extract a
subsequence $(r_{\sigma(n)} \cdot \gamma_{\sigma(n)},
\gamma_{\sigma(n)}) $ that converges to $ (g_1,g_2) \in K $. Then $
g_1 $ and $g_2^{-1}$ are composable and $(r_{\sigma(n)})_{n \in
{\mathbb N}} $ converges to $ r= g_1 \cdot g_2^{-1}$. Since $R$ is
closed as a subset of $\Gamma_L^L $, $r$ is an element in $R$. By
construction $(r,g_2)$ is an element in $  \mu^{-1}(K) $ which is
therefore a compact subset of $ R \times_{t,M,s} \Gamma_L $.
This justifies the claim.

The action of the Lie groupoid $R \toto L $ on $\Gamma_L $
 being a free and proper action,
the quotient space  $ Z(R) := R \backslash \Gamma_L$ is
 a Hausdorff manifold
and the projection $p: \Gamma_L  \to Z(R)$ is a surjective submersion.

Since $t(r \cdot \gamma) = t(\gamma) $ for any
$(r ,\gamma) \in  R \times_{t,L,s} \Gamma_L $,
there exists an unique map $\phi: Z(R) \to M $
such that the diagram (\ref{eq:comdia1}) commutes.
Moreover, the identity
 $ \phi(Z(R)) = t(\Gamma_L)$ holds by construction of $\phi$.
But, by definition of an algebroid crossing, $L$ has a non-empty
intersection with all the symplectic leaves contained in
$\bar{{\mathcal S}}$, so that,  by Lemma \ref{lem:leaf0}, any
element of $\bar{{\mathcal S}}$ is the target of at least one element in
$\Gamma_L$, and we obtain the identity $\phi(Z(L)) = t(\Gamma_L) = \bar{\mathcal
S}$.

The identity
 $  t(\Gamma_{L \cap {\mathcal S}}) ={\mathcal S}$
 holds by Lemma \ref{lem:leaf0}.
By construction of $\phi $,
 the identity  $ p(\Gamma_{L \cap {\mathcal S}})
= \phi^{-1}({\mathcal S})$ holds as well and $\phi:p( \Gamma_{L \cap
{\mathcal S}} ) \to {\mathcal S}  $ is a surjective submersion
because  $t: \Gamma_{L \cap {\mathcal S}} \to {\mathcal S}$ itself a
surjective submersion  by Lemma \ref{lem:leaf0} again.
Now, the dimension of $Z(R)  $ is given
by
 $$ dim( Z(R)) = dim(\Gamma_L) - rk (B) = dim(L)  + rk(A) - rk (B)$$
But since $\rho(B_x) = T_xL $ for all  $x \in L \cap {\mathcal S} $
while $ \rho(A_x) = T_x {\mathcal S} $, the ranks of the algebroids
$A$ and $ B$ are given by
 $ rk(A) = dim({\mathcal S})+ dim(ker(\rho) ) $
and $rk(B) = dim(L) +  dim(ker(\rho)) $
Hence $rk(A) - rk(B) = dim({\mathcal S}) - dim(L)$ and we obtain
$$ dim( Z(R)) =  dim(L) - dim(L) + dim({\mathcal S})= dim({\mathcal S}) .$$
The dimensions of the manifolds $Z(R) $ and ${\mathcal S} $ being
equal and  $\phi $ being a surjective submersion from the open
subset $\phi^{-1}({\mathcal S}) \subset Z(R)$  onto ${\mathcal S} $,
the map  $\phi: \phi^{-1}({\mathcal S}) \to {\mathcal S} $ is an \'etale map. 

It remains to prove that $\phi^{-1}({\mathcal S})  $ is open and dense in $Z(R)
$. But $ L \cap { \mathcal S} $ is open in $L $ by Lemma \ref{lem:leaf} and dense in $L$ by assumption.
Hence $\Gamma_{ L \cap { \mathcal S}} $ is open and dense in $\Gamma_L $.
Since  $p$ is a surjective
submersion (and in particular an open map), $\phi^{-1}({\mathcal S})= p(\Gamma_{L \cap {\mathcal S}})$
is open and dense in $Z(R) $. This achieves the proof of~1).

2)  By Lemma \ref{lem:leaf0}, the Lie groupoid $\Gamma \toto M$ acts transitively on ${\mathcal S}$, hence for any $x \in {\mathcal S}$
the target map induces a biholomorphism/diffeomorphism
$$ I_x(\Gamma)  \backslash \Gamma_x \simeq {\mathcal S} ,$$
 where $I_x(\Gamma)$ stands for the isotropy group at $x \in {\mathcal S}$ of the Lie groupoid $\Gamma \toto M$.

For any $x \in  L \cap {\mathcal S} $, the anchor map of the
algebroid $B \to L $, is onto so that connected components of $L
\cap {\mathcal S} $ are algebroid leaves for the algebroid $B \to L
$, and, by Lemma \ref{lem:leaf0} again, the action of $R \toto L $
on $  L \cap {\mathcal S} $ acts transitively on each connected
component. Since $L \cap {\mathcal S} $ is $R$-connected, for any
two connected components, there exists $r \in R $ having its source
in the first one and its target in the second one. Hence the action
of $R \toto L $ on $  L \cap {\mathcal S} $ is transitive, we obtain
a  biholomorphism/diffeomorphism
  $$  \phi^{-1} ({\mathcal S}) \simeq  R \backslash \Gamma_{L \cap {\mathcal S}} \simeq
  I_x\big( R\big)  \backslash \Gamma_x $$
 where  $I_x\big( R\big)$
  stands for the isotropy group at the point $x \in L \cap {\mathcal S}$
of the Lie groupoid $R \toto L$. Since $\Gamma_x $ is a connected
set, so is $\phi^{-1} ({\mathcal S}) =  p (\Gamma_x) $, and the
\'etale resolution $(Z,\phi) $ is a covering resolution. Moreover,
the typical fiber of  $\phi:\phi^{-1} ({\mathcal S}) \to {\mathcal
 S}$ can be identified with the
quotient
$$\frac{  I_x \big( \Gamma\big) }{ I_x\big(  R\big)  } .$$
We study this quotient in detail. For any $x \in L \cap {\mathcal
S}$,
 the kernel of the anchor maps of the algebroids
$B \to L$ and $A \to M$ coincide, therefore the isotropy Lie
algebras of both Lie groupoids at an arbitrary point
 $x \in L \cap {\mathcal S}$ coincide; therefore their isotropy Lie
groups at $x \in   L \cap {\mathcal S}$ have the same connected
component of the
 identity. Denoting, for an arbitrary
 Lie group $H$, by $\pi_0(H) $ the discrete group of connected
 components
of $H$, there is an isomorphism
 $$  \frac{  I_x\big( \Gamma\big) }{ I_x\big(
  R\big) } \simeq \frac{  \pi_0 ( I_x\big(
  \Gamma \big)) }{ \pi_0  ( I_x\big(
  R\big)   ) }  . $$
This achieves the proof of 2)

3) The restriction to $\phi^{-1}({\mathcal S}) \subset Z(R) $ of
$\phi
   $
is an biholomorphism/diffeomorphism if and only if for any two
$\gamma_1,\gamma_2 $ in $\Gamma_{L \cap {\mathcal S}} $ with
$t(\gamma_1) =t(\gamma_2) $, $\gamma_1 \gamma_2^{-1}   $ is an
element of $R$. In other words, the restriction to
$\phi^{-1}({\mathcal S}) \subset Z(R) $ of $\phi
   $
is a biholomorphism/diffeomorphism if and only if $\Gamma_{L \cap
{\mathcal S}}^{L \cap
   {\mathcal S}} \subset R$. Since
$R $ is a closed subset of $\Gamma_L^L$ by assumption,
 this last requirement is equivalent to
the requirement $ \overline{\Gamma_{L \cap {\mathcal S}}^{L \cap
   {\mathcal S}}} \cap \Gamma_L^L \subset R$.

But, $L \cap {\mathcal S}$ being by assumption open and dense in $L
$,  $R \cap \Gamma_{L \cap {\mathcal S}}^{L \cap
   {\mathcal S}} $ is dense in $R$. In other words, $R \subset \overline{\Gamma_{L \cap {\mathcal S}}^{L \cap
   {\mathcal S}}} $. Since the inclusion $R \subset \Gamma_L^L $ holds
also, the inclusion $R \subset  \overline{\Gamma_{L \cap {\mathcal
S}}^{L \cap
   {\mathcal S}}} \cap \Gamma_L^L $ holds.

In conclusion the inclusion $ \overline{\Gamma_{L \cap {\mathcal
S}}^{L \cap
   {\mathcal S}}}\cap \Gamma_L^L \subset R$ holds if and only if
the equality $ \overline{\Gamma_{L \cap {\mathcal S}}^{L \cap
   {\mathcal S}}} \cap \Gamma_L^L= R $ holds. This achieves the proof of 3).

4) Fix $x \in L \cap {\mathcal S}$. We have a commutative diagram of fiber bundles:
$$   \xymatrix{ I_x(R) \ar[r]  \ar[d]    &    R_x   \ar[r]^{t}  \ar[d]
&   L \cap {\mathcal S}\ar[d]  \\  I_x(\Gamma)  \ar[r] & \Gamma_x
  \ar[r]^{t} &
  {\mathcal S}   }   $$
  where vertical maps are inclusion maps.
 The usual long exact sequence in
homotopy gives the following commutative diagram, where horizontal
lines are exact sequences:
 $$   \xymatrix{ \cdots  \ar[r]  &  \pi_1(  R_x)  \ar[r]  \ar[d]    &
 \pi_1 (  L \cap {\mathcal S})   \ar[r]  \ar[d]^j
&  \pi_0 (I_x(R)) \ar[r]  \ar[d] & \pi_0(R_x  ) \ar[d] \ar[r] &\cdots\\
 \cdots \ar[r]  &
  \pi_1(  \Gamma_x) \ar[r]     & \pi_1 (   {\mathcal S})   \ar[r]
&   \pi_0 (I_x(\Gamma))  \ar[r] & \pi_0(\Gamma_x  )  \ar[r] &\cdots }$$
By assumption, $\Gamma_x$ and $R_x $ are connected while
$\Gamma_x $ is simply-connected. Hence, in the central square of the
 previous diagram, that is to say, in the diagram
 $$ \xymatrix{  \pi_1 (  L \cap {\mathcal S})  \ar[r]  \ar[d]^j    &   \pi_0 (I_x(R))  \ar[d] \\
  \pi_1 (   {\mathcal S})     \ar[r]
&  \pi_0(  I_x(\Gamma))  } $$
 the  upper line is a surjective group morphism,
the lower line is a group isomorphism, while the right vertical line
is an injective group morphism (since it is induced by the
inclusion). The proof of 4) now follows by an easy diagram chasing.
\end{proof}

At first glance, it seems that  resolutions obtained  from the
 previous   construction depends strongly on the choice of an algebroid crossing.
We study this dependence using the notion of Morita equivalence 
of Lie groupoids. 
This notion of Morita equivalence goes back to
\cite{HS}, but we invite the reader to consult \cite{LTX1}
for an introduction to the notion of Morita equivalence of 
modules of Lie groupoids that matches the presentation below.

\begin{prop}\label{prop:repr_isom}
Let $(A \to M,\rho,[\cdot,\cdot]) $ be a Lie algebroid, ${\mathcal
S} $  a locally closed orbit of this algebroid, and $L_i, i=1,2$ two algebroid
 crossings of $\bar{\mathcal S}$ with normalisations $B_i \to L_i$. Assume that:
\begin{enumerate}
\item there exists a source-connected Hausdorff Lie groupoid $\Gamma
 \toto M $ integrating
the Lie algebroid $(A \to M,\rho,[\cdot,\cdot])  $ and
\item  there exists, for $i=1,2$,  sub-Lie groupoids $R_i \toto L_i $ of $ \Gamma \toto M$,
  closed as  a subset of
  $\Gamma_{L_i}^{L_i} $, integrating
the Lie algebroid $(B_i \to L_i,\rho,[\cdot,\cdot])  $,
and containing $\Gamma_{L_i \cap {\mathcal S}}^{L_i \cap {\mathcal S} }  $.
\end{enumerate}
 Let $(Z_i,\phi_i) $, $i=1,2 $  be the two corresponding  resolutions
as in Proposition  \ref{prop:desing+groupoid}(3).

The following  are then equivalent:
\begin{enumerate}
\item[(i)] the resolutions $(Z_1,\phi_1) $ and $(Z_2, \phi_2 )$ are isomorphic,
\item[(ii)]  $  \overline{ \Gamma_{L_1 \cap {\mathcal S}}^{L_2    \cap {\mathcal S}}} \cap \Gamma_{L_1}^{L_2} $
  is a submanifold of $ \Gamma$, and the restrictions to this submanifold of the source and the  target
  maps are surjective submersions onto $L_1 $ and $L_2 $ respectively,
 \item[(iii)] there exists a submanifold $I$ of $\Gamma $ 
that gives a Morita equivalence
 between the Lie groupoids $R_1 \toto L_1 $ and $R_2 \toto L_2 $:
  $$  \xymatrix{ \Gamma_{L_1} \ar[rd] & R_1 \dar[d]  & \ar[ld]^s I  \ar[rd]_t & R_2\dar[d] & \ar[ld] \Gamma_{L_2}  \\ & L_1 & & L_2 &   } $$ 
 In this case moreover,  the $R_1 $-module $\Gamma_{L_1} $
 corresponds to the $R_2 $-module $\Gamma_{L_2} $ with respect to the Morita equivalence $I$. 
\end{enumerate}
\end{prop}

We need some preliminary results. Let $(Z(R)= R \backslash
\Gamma_L,\phi) $ be an \'etale resolution as in Proposition
\ref{prop:desing+groupoid}(1). There is a natural map $j$ from $L $ to
$Z(R)$ obtained by composing the restriction to $L$ of the unit map
$\epsilon: L \to \Gamma_L $ with the canonical projection $p :
\Gamma_L \to Z(R) $. By Equation~(\ref{eq:comdia1}), $\phi $ is a left
inverse of $j$, i.e. $\phi \circ j = {\rm id}_L $, so that $j(L) $
is a submanifold of $Z(R)$. In short:

\begin{lem} \label{lem:j(L)}
Let $(Z(R) = R \backslash \Gamma_L,\phi) $ be an \'etale resolution
as in Proposition \ref{prop:desing+groupoid} and $j: L \to Z(R)$ as
above. Then $j(L) $ is a submanifold of $Z(R) $ and the restriction
of  $ \phi$ to $j(L) $ is a biholomorphism/diffeomorphism onto $L$.
\end{lem}

We now prove Proposition \ref{prop:repr_isom}.

\begin{proof}
We say simply diffeomorphism rather than
biholomorphism/diffeomorphism.

{\it (i)} $\Longrightarrow$  {\it (ii)}. Denote by
$p_i, i=1,2 $ and $j_i, i=1,2 $ the projections from $\Gamma_{L_i}
\to Z_i $ and the inclusions of $L_i $ in $Z_i $ respectively.

Let $\Psi: Z_1 \to Z_2 $ be a diffeomorphism of resolutions from
 $(Z_1,\phi_1) $ to $(Z_2,\phi_2) $. Then $  \Psi \circ p_1 :
\Gamma_{L_1} \to Z_2 $ is a surjective submersion. The inverse image
$ I_1 $ of the submanifold $j_2(L_2) $ of $ Z_2$ by $  \Psi \circ
p_1$ is a submanifold of $\Gamma_{L_1} $, which is closed as a
subset of $\Gamma_{L_1}^{L_2}$, and, by Lemma \ref{lem:j(L)},
 $t=\phi_2 \circ \Psi \circ p_1 :
I_1 \to L_2 $ is a surjective submersion.
 In particular, $I_1$ contains $ \Gamma_{L_1 \cap {\mathcal S}}^{L_2 \cap {\mathcal S}} $ as a dense
subset. Hence $I_1=\overline{\Gamma_{L_1 \cap {\mathcal S}}^{L_2
\cap {\mathcal S}}} \cap \Gamma_{L_1}^{L_2}  $ and the latter is a submanifold of $\Gamma $.

Similarly, since $\Psi^{-1}: Z_2 \to Z_1 $ is a diffeomorphism, $
\Psi^{-1} \circ p_2 : \Gamma_{L_2} \to Z_1$ is a surjective
submersion. The inverse image $ I_2 $ of the submanifold $j_1(L_1) $
of $ Z_1$ by $  \Psi^{-1} \circ p_2$ is equal to $
\overline{\Gamma_{L_2 \cap {\mathcal S}}^{L_1 \cap
    {\mathcal S}}} \cap \Gamma_{L_2}^{L_1}  $ and $t:I_2 \to L_1 $
is surjective submersion.

In particular, we obtain $I_2 =I_1^{-1} $ and,
since the inverse map intertwines sources and targets,
 $s:I_1 \to  L_1$  is also a surjective submersion. This completes the proof
 of {\em (i)} $\Longrightarrow $ {\em (ii)}.

 {\it (ii)}  $\Longrightarrow$  {\it (iii)}. Let $ I= \overline{\Gamma_{L_1 \cap {\mathcal S}}^{L_2 \cap {\mathcal S}}} \cap \Gamma_{L_1}^{L_2}  $.
For any $m \in L_1 \cap {\mathcal S}$, the Lie groupoid $R_2 \toto L_2$ 
acts transitively on the fiber of $s: I \to L_1$ over $m$, since this fiber is precisely equal to $\Gamma_m^{L_2 \cap {\mathcal S}} $
while $R_2 $ contains $\Gamma^{L_2 \cap {\mathcal S}}_{L_2 \cap {\mathcal S}} $.
Let us show that this fact remains true for all $m \in L_1 $. Let $c,c' \in I$ be two points in $I$
with $s(c) = m = s (c') $. There exist sequences $(c_n)_{n \in {\mathbb N}}, (c_n')_{n \in {\mathbb N}}$ 
of elements of $I$,
converging to $c $ and $c'$ respectively, and which satisfy $s(c_n )= s(c_n') $ 
for all $n \in {\mathbb N} $. The sequence $ (c_n')^{-1}c_n $ is a sequence of $R_2 $ 
that converges to $ (c')^{-1} c$. Since $R_2 $ is closed in $\Gamma_{L_2}^{L_2} $,
we have $(c')^{-1} c  \in R_2$, and $R_2 \toto L_2 $ acts transitively on the fibers of $ s: I \to L_1$.
This action is also, of course, free.

   Similarly, the fibers of $t: I \to
L_1$ are precisely the fibers of the $R_1 $-action. The left and right actions of $R_1 \toto L_1 $ and $R_2 \toto L_2 $
is free and proper, with $I_1 \simeq I / R_2 $ and $ R_1 \backslash  I \simeq l_2$. 
Hence $ I$ is a Morita bimodule that gives a Morita equivalence between the Lie groupoids $R_1 \toto L_1 $ and $R_2 \toto L_2 $.

We now prove that,
in this case moreover,  the $R_1 $-module $\Gamma_{L_1} $
 corresponds to the $R_2 $-module $\Gamma_{L_2} $ with respect to the Morita equivalence $I$. 
W e recall some general facts about Morita equivalences (see \cite{LTX1} or \cite{Xu:2004}). A Morita bimodule between
two Lie groupoids induces a one-to-one correspondence between their left modules that we now describe. Let $(X_2,\chi_2) $ be a left $R_2 $-module.
Then the diagonal action of $R_2 \toto L_2 $ on  $I \times_{s,L_2,\chi_2} X $ is free and
proper, so that the quotient space 
$$ X_1:=\frac{I \times_{s,L_2,\chi_2} X_2}{R_2} = \frac{\left\{ (c,x) \in I \times X_2 \, | \, t(c) = \chi_2(x) \right\}  }{ (c,x)  \sim  (c r^{-1}, r x) } $$
 is a manifold again. The map $t: I \to L_2 $ factorizes to yield a surjective submersion 
$\chi_1: X_1  \to L_1 $. The left action of the Lie groupoid $R_1 \toto L_1 $ on $I $ induces a left action of  $R_1 \toto L_1 $
 on $(X_1,\chi_1) $ which gives the desired structure of left $R_1$-module. 

In the particular case where $(X_2,\chi_2) = (  \Gamma_{L_2},s   ) $ and the Morita equivalence is with respect 
to $I$,  we easily check that $X_1  $ is isomorphic to $Z_1 = R_1 \backslash \Gamma_{L_1} $, 
the isomorphism in question being simply given by  
$$    [c,x_2]  \to p_1(c x_2)   $$
where $(c,x_2) $ is an element of $I \times_{s,L_2,s} \Gamma_{L_2}$, $[c,x_2] $ its class in $\frac{I \times_{s,L_2,\chi_2} X_2}{R_2}
$ and $p_1: \Gamma_{L_1} \to Z_1 = R_1 \backslash \Gamma_{L_1} $ the canonical projection. This proves 3).

 {\it (iii)}  $\Longrightarrow$  {\it (i)}.
Under the correspondence between left $R_1 $- and $R_2 $-modules just described above,     
when $X_2 = R_2 \backslash \Gamma_{L_2} $ is a manifold, then so is $X_1 = R_1\backslash \Gamma_{L_1} $
and  both manifolds are canonically isomorphic (see  \cite{Xu:2004}). The isomorphism  $\Psi$ in question is given by $[x_2] \to [(i,x_2)] $,
where $x_2 \in X_2$ is an arbitrary element, $[x_2] $ its class in  the quotient space $R_2 \backslash X_2 $, $i \in I $ is any element with $t(i) =
\phi_2(x_2)$, and $[i,x_2]$ is the class of $(i,x_2) $ in the quotient space $ R_1 \backslash X_1 $. In the particular case $X_2 =
\Gamma_{L_2} $ and $\chi_2 = s $, $\Psi$  maps $R_2 \backslash \Gamma_{L_2} $ to $R_1 \backslash \Gamma_{L_1} $.
 Moreover,  $\Psi $ satisfies by construction the relation $  \phi_1 \circ \Psi = \phi_2 $  and is therefore an isomorphism of resolutions.
\end{proof}

\section{Resolution compatible with a multi-vector field.}
\label{sec:R2}

\subsection{Definition of a resolution compatible with a multi-vector field}
\label{sec:R2_D2}

By an abuse of language (justified by its use in algebraic geometry),
 a $k$-vector field
 $\pi_M$ on the manifold $M$ is said to be {\em
  tangent to  $\bar{\mathcal S} $} if it is tangent to ${\mathcal S}
$, i.e.  if $\pi_M[f,f_2,\cdots,f_{k}]    $ vanishes on ${\mathcal
S} $ for any local function $f$ that vanishes on ${\mathcal S} $
(equivalently on $\overline{\mathcal S}$),
and any local functions $f_2,\cdots,f_k \in {\mathcal F}(M)$.

Let us now explain what we mean by a resolution compatible with a
$k$-vector field $\pi_M$ tangent to $\bar{\mathcal S} $.

\begin{defn} \label{def:comp_desing}
Let $\bar{\mathcal S}$ be the closure of a complex/real locally
closed (= embedded) submanifold ${\mathcal S} $ of $M$ and $\pi_M $
a $k$-vector field tangent to $\bar{\mathcal S} $, a {\em resolution
compatible with the $k$-vector field $\pi_M $ } is a pair
$(Z,\phi)$ where
\begin{enumerate}
\item $(Z,\phi) $ is a resolution of $\bar{\mathcal S} $,
\item the $k$-vector field $\pi_Z $ defined on  $\phi^{-1}({\mathcal
  S}) $ by $\phi_* \pi_{Z} = \pi_M $ extends to a holomorphic/smooth
  $k$-vector field on $Z$.
\end{enumerate}
When $(Z,\phi)$ is only an \'etale/covering resolution, then we speak of an
 {\em \'etale/covering  resolution compatible with the $k$-vector
 field $\pi_M$}. If moreover, $ \pi_M$ is a Poisson bivector field,
we speak of a {\em Poisson resolution}.

We denote by $\pi_Z $ again the extension to $Z$ of $\pi_Z $. For
convenience,
we often denote a compatible resolution by a triple
$(Z,\phi,\pi_Z ) $. But we invite the reader to keep in mind that, for a compatible \'etale resolution,
$\pi_Z$ is indeed determined by $(Z,\phi)$.
 \end{defn}

\begin{rmk}
For any resolution compatible with $\pi_M $,
the relation $ \phi_* \pi_{Z} = \pi_M  $  holds indeed  on $Z$
by density of $\phi^{-1}({\mathcal  S}) $ in $Z$.
In particular, when $\pi_M $ is  a Poisson bivector field, so is $\pi_Z $.
\end{rmk}

For clarity, we list all the conditions required in order to
have an \'etale resolution $(Z,\pi_Z,\phi)$  of the closure $\bar{\mathcal S}
$ of a locally closed submanifold of $M$ which is compatible with a
$k$-vector field $\pi_M $ tangent to $\bar{\mathcal S}$.

\begin{enumerate}
\item $Z$ is a manifold and $\pi_Z $ a $k$-vector field,
\item  $\phi:Z \to M$ is a  holomorphic/smooth map from the manifold
$Z $ to the manifold $M$,
\item $\phi(Z)=\bar{{\mathcal S}}$,
\item $\phi^{-1}({\mathcal S})$ is open and dense in $Z$ (note that density is
  automatically satisfied  in the
  complex case),
\item  the restriction of $\phi$ to a map from $\phi^{-1}({\mathcal
 S})$ to ${\mathcal S}$  is a biholomorphism/diffeomorphism,
\item $\phi: Z \to M$  maps $\pi_Z $ to $\pi_M $.
\end{enumerate}

For \'etale/covering symplectic resolutions, the fourth point needs
to be replaced by ``the restriction of $\phi$ to
$\phi^{-1}({\mathcal S})$ is an \'etale/covering map over
${\mathcal S}$''.

The following Lemma will be useful in the next section.
Notations are as before.

\begin{lem}\label{lem:etalemorph}
Let $(Z_i,\phi_{Z_i})$, with $i=1,2 $,
be \'etale resolutions of $\bar{\mathcal S} $,
and $\Psi: Z_1 \to Z_2 $  a local diffeomorphism 
which is a  morphism of resolutions of $\bar{\mathcal S} $
$$ \xymatrix{
Z_1 \ar[r]^{\Psi} \ar@<0.5ex>[rd]_{\phi_{Z_1}}  & Z_2 \ar@<0.5ex>[d]^{\phi_{Z_2}} \\
  & \bar{\mathcal S} }$$
Then, if $(Z_1,\phi_{Z_1})$ is an \'etale
 resolution compatible with  $\pi_M $, then  $(Z_2,\phi_{Z_2}) $ 
is also an \'etale resolution compatible with  $\pi_M $. 
\end{lem}
\begin{proof}
We denote by $\pi_{Z_1} $ the unique $k$-vector field on $Z_1$ satisfying
 $(\phi_{Z_1})_* \pi_{Z_1} = \pi_M$.
 For any two points $x,y  \in Z_1$ with $\Psi (x) = \Psi (y) $,
there exists a  local biholomorphism/diffeomorphism $\Phi $
of $Z_1 $ over the identity of   $Z_2$,
(id est  $\Psi \circ \Phi =\Psi $) defined on
some open subset $U \subset Z_1$ with $\Phi (x) = y $.
The relation $ \Psi_*( \pi_{Z_1} )_{|_x}=
  (\pi_{Z_1})_{|_{y}}  $
holds for any $ x  \in U \cap \phi_1^{-1}({\mathcal S})$ since
   $$ (\phi_{Z_1})_* (\pi_{Z_1})_{|_x} = (\pi_M)_{|_{\phi_{Z_1}( x)}} =
(\phi_{Z_1})_* (\pi_{Z_1})_{|_{y}}  $$ and since $(\phi_{Z_1})_* : T_m
Z_1: T_{\phi_{Z_1}(m)} {\mathcal S}$
is invertible for all  $x \in U \cap \phi_1^{-1}({\mathcal S})$. 
The open subset $U \cap \phi_{Z_1}^{-1}({\mathcal S})$ being
dense in $U$, the relation $ \Phi_*( \pi_{Z_1} )_{|_x}=
  (\pi_{Z_1})_{|_{y}}  $ holds for any $x \in U$.
This amounts to the fact that  $\pi_{Z_1} $ goes to the quotient and defines
a $k$-vector field $\pi_{Z_2} $ on $Z_2 $. By construction, $\Psi_* \pi_{Z_1}
=\pi_{Z_2} $. It is immediate that
$$(\phi_{Z_2})_* \,  \pi_{Z_2} = (\phi_{Z_2})_* \circ \Psi_*  \, \pi_{Z_1}=
 (\phi_{Z_1})_* \, \pi_{Z_1}=\pi_M .$$
The $k$-vector field $\pi_{Z_2} $ is defined everywhere on
 $Z_2 $ and projects on $\pi_M $ through $(\phi_2)_* $, so that
 the \'etale resolution $(Z_2,\phi_2) $ is compatible with $\pi_M $.
\end{proof}

\subsection{Lie groupoids, multiplicative multivector fields 
 and compatible resolutions of the closure of an algebroid leaf}
\label{sec:R2G}

 Recall (see Definition 2.6 in \cite{ILX}) that a
$k$-vector field on a Lie groupoid $\Gamma \toto M $ is said to be
multiplicative
when the graph of the multiplication of the groupoid, that is to say the
submanifold of $ \Gamma^3$ given by
$$ {\rm Gr(\Gamma)} = \{(\gamma_1,\gamma_2,\gamma_1\gamma_2) \in \Gamma^3 |
t(\gamma_1)= s(\gamma_2)\} ,$$
is coisotropic with respect to $\pi_{\Gamma} \oplus \pi_{\Gamma}\oplus
(-1)^{k+1}\pi_{\Gamma}$.

Multiplicative vector fields on Lie groupoids have a very rich geometry.
We invite the reader to read Section 2 in \cite{ILX} to get an
overview of that difficult matter, and recall two points of fundamental
importance for the present purpose.

We recall (see Remark 2.4 in \cite{ILX}) that a submanifold $N$ 
of $M$ is {\em coisotropic with respect to a $k$-vector field $\pi_M
  \in {\mathcal X}^k(M)$}
if and only if, for any local functions
$f_1,\dots,f_k \in {\mathcal F}(M)$ that vanish on $N$, the function
$ \pi_M [f_1,\dots , f_k] $ vanishes on $N$. 

The purpose of the present section is to 
\begin{enumerate}[{(}i{)}]
\item  recall from \cite{ILX} 
how this multiplicative  $k$-vector field $\pi_{\Gamma} $ induces a $k$-vector
  field $\pi_M $ on $M$, and explain 
why, for $k\geq 2 $, this $k$-vector field is tangent to $\bar{\mathcal S} $,
\item show that any resolution $(Z(R) = R \backslash \Gamma_L,\phi) $ 
of an algebroid leaf $\overline{\mathcal S} $ constructed with the
help of an algebroid crossing $L$ as in  Proposition
\ref{prop:desing+groupoid}
is compatible with
$\pi_M $  under the condition that $L$ is coisotropic with respect to $ \pi_M$.
\end{enumerate}

(i) {\bf From multiplicative $k$-vector field to $k$-vector field on
  ${\mathcal S}$.}

We prove in this section that a multiplicative $k$-vector field on $\Gamma \toto M $
induces a $k$-vector field on $M$ which is tangent to all
algebroid leaves for $k \geq 2 $, see Proposition \ref{prop:tangency}
below. We characterise, in term
of $k$-differentials on the Lie algebroid, algebroid crossing which are coisotropic with
respect to this induced $k$-vector field, see Proposition
\ref{prop:cois-ckvector} below.

To start with, according to  Proposition 2.21 in \cite{ILX}, there exists a unique
$k$-vector field $\pi_M $ on $M$ such that 
   \begin{equation}
\label{eq:proje} s_* \pi_{\Gamma} =\pi_M   \mbox{ and }  t_* \pi_{\Gamma}
   =(-1)^{k+1}\pi_M.\end{equation}

We now prove the following Proposition.

\begin{prop}  \label{prop:tangency}
For any multiplicative $k$-vector field $\pi_{\Gamma} $ on 
$\Gamma \toto M $ with $ k \geq 2$, the $k$-vector field $\pi_M  = s_* \pi_{\Gamma}$ 
is tangent to all the algebroid leaves.
\end{prop}
\begin{proof}
According to Propositions 2.17,2.18 and Corollary 2.20
in \cite{ILX}, there exists a $k$-differential $\delta:\Gamma(\wedge^\bullet A \to M)
 \to \Gamma(\wedge^{\bullet+k-1}A \to M ) $ such that, 
 \begin{equation}
\label{eq:rightmult}
\left\{
\begin{array}{rcl}
\schouten{ \pi_{\Gamma}, \overrightarrow{a}}_{T\Gamma} & = &  \overrightarrow{\delta (a)} \\
\schouten{  \pi_{\Gamma}, \overleftarrow{a} }_{T\Gamma}  & =& \overleftarrow{\delta (a)} \\
\end{array}\right.
\end{equation}
where $\overrightarrow{a} $ (resp. $\overleftarrow{a}  $) is the
right-invariant (resp. left-invariant) multi-vector field 
on $\Gamma $ corresponding to $a \in \Gamma(\wedge^l A \to M)$, with the understanding that
$\overrightarrow{f}= s^* f $ (resp. $\overleftarrow{f}= t^* f $) when
$f$ is a function on $M$ (see Section \ref{Sec:Notations}
and especially Equation~(\ref{k-dif-prop}) for the definition of a $k$-differential).

 Note that \cite{ILX} deals with real
Lie groupoids only, but these results extend to the complex setting
without any difficulty, with the understanding, however,
that $\Gamma(\wedge^\bullet A \to M)$
 stands for the sheaf of local sections.

For any local functions $g_1,\cdots,g_{k-1},f $ on $M$,
with $f $ vanishing on ${\mathcal S}$, we have, in view of Equation~(\ref{eq:proje})
  $$ s^* (\pi_M[g_1,\cdots,g_{k-1},f])=  \pi_\Gamma [s^*g_1,\dots ,s^*
  g_{k-1}, s^* f]  =  \schouten{ \pi_\Gamma, s^* g_1  }_{T\Gamma} [ s^* g_2,\cdots,s^* g_{k-1},s^* f  ]  .$$
Equation~(\ref{eq:rightmult}) yields
  $$  s^* (\pi_M[g_1,\cdots,g_{k-1},f])=  \overrightarrow{\delta(g_1)}
  [ s^* g_2,\cdots,s^*g_{k-1},s^*f  ] .$$
But, in turn, we have
 $$ \overrightarrow{\delta(g_1)}
  [ s^* g_2,\cdots,s^*g_{k-1},s^*f  ] =   s^* \left( \, \rho (  \delta(g_1)) [g_2,\cdots,g_{k-1},f
  ] \, \right) .$$ Since $f$ vanishes on  $ {\mathcal S}  $, and since
  $\rho(\delta(f_1)) $ is tangent to all the algebroid leaves, 
$\rho (  \delta(g_1)) [g_2,\cdots,g_{k-1},f]$  vanishes on ${\mathcal
  S} $ as well. In conclusion, the function
 $$  \pi_M[g_1,\cdots,g_{k-1},f]  $$
  vanishes on ${\mathcal S} $ for any
 local functions  $g_1,\cdots,g_{k-1} $ provided that $f$ vanishes on ${\mathcal S} $.
Hence $\pi_M $ is tangent to ${\mathcal S} $. 
\end{proof}

Assume now that we are given an algebroid crossing $L$
of $\overline{{\mathcal S} }$ with normalisation $B \to L $, where
${\mathcal S} $ is a locally closed (= embedded) 
algebroid leaf.
We now explore the picture for $L$
 coisotropic with respect to 
 $\pi_M $.

We say that a $k$-differential $\delta : \Gamma(\wedge^\bullet A_{|_U} \to U) \to
\Gamma(\wedge^{\bullet+k - 1} A_{|_U} \to U)$ is {\em compatible with an algebroid
crossing $L$ with normalization $B \to L$}
 if and only if $\delta(\left\langle B\right\rangle_U) \subset \left\langle B\right\rangle_U $ 
for all open subset $U \subset M $, where
 $\left\langle B\right\rangle_U $ is  
the ideal in $  (\Gamma(\wedge^\bullet A_{|_U} \to U),\wedge)  $
generated by functions vanishing on $L$ and sections of $A$
whose restriction to $L$ is a section of $B \to L$.

In other words, a  local section $ X $ of $\Gamma(\wedge^p A_{|_U} \to U)$
belongs to $\left\langle B\right\rangle_U $ if, and only if, for all $m \in L \cap U$,
$X_{|_m} $ is a linear combination of terms of the form
$ b \wedge a_1 \wedge \cdots \wedge a_{p-1}$, with $b \in B_m $,
$ a_1 , \cdots , a_{p-1} \in A_m $.

The need of the following Proposition will appear in the sequel.

\begin{prop}\label{prop:cois-ckvector}
Let $k \geq 2 $ be an integer and
$\Gamma \toto M$ a Lie groupoid
endowed with a multiplicative $k$-vector field
 $\pi_\Gamma $. 
Denote by $\pi_M $ the $k$-vector field 
 on $ M$ induced by $\pi_\Gamma$ by Equation  (\ref{eq:proje}), and, by $\delta$ 
 the $k$-differential $\delta $ 
induced by Equation~(\ref{eq:rightmult}).

Let  ${\mathcal S} $ be a locally closed 
algebroid leaf. An algebroid crossing $ L $ of $\overline{{\mathcal S} }$, with normalisation $B \to L$,
is compatible with  the $k $-differential $\delta $ if and only if $L $ is coisotropic  with respect to $\pi_M $.
\end{prop}

The following Lemma is an immediate consequence of the density of
$ L  \cap {\mathcal S}  $ in $L$ and from the defining relation
$B_m = \rho^{-1}(T_m L) $ for all $m \in  L  \cap {\mathcal S}$.

\begin{lem}
\label{lem:coiso}
Let $U \subset M  $ be an open subset.
\begin{enumerate}
\item  For any  homogeneous element $X \in \left\langle B\right\rangle_U $ of degree $\geq 1$, 
$L \cap U$ is coisotropic with respect to  $\rho(X) $.
\item Conversely, any section $X \in \Gamma (\wedge^i A_{|_U} \to U) $ 
with $i \geq 1 $ such that $L \cap U$ is coisotropic with respect to $\rho(X) $ is a section
of $\left\langle B\right\rangle_U $. 
\end{enumerate}
\end{lem}

Now, we can turn our attention to the proof of Proposition \ref{prop:cois-ckvector}.

\begin{proof}
 Let $U \subset M $ be a trivial open subset of $M$.
Assume that $ B \to L$ is compatible with the $k$-differential
$\delta $.
 To start with, we recall from Lemma 2.32 in \cite{ILX}
  that the following
identity holds
\begin{equation}\label{eq:2.32}
 \pi_M[f_1 ,\dots, f_k]= (-1)^{k+1}
 \rho (\delta(f_1))[f_2 ,\dots, f_k]  \end{equation}
for any local functions $f_1,\cdots,f_k \in {\mathcal F}(U)$.
Now, let $f_1,\dots,f_k$ be functions that vanish on $L \cap U $.
The $k$-differential 
$\delta $ being compatible with the algebroid crossing $B \to L $,
we have $\delta(f_1)  \in \left\langle B\right\rangle_U$. By Lemma \ref{lem:coiso}-(1),
$L \cap U$ is coisotropic 
with respect to  $\rho (\delta(f_1)) $, so that the function
$ \pi_M[f_1 ,\dots, f_k]= (-1)^{k+1} \rho (\delta(f_1))[f_2 ,\dots, f_k]  $
vanishes on  $L \cap U $. In conclusion, $L$ is coisotropic with respect to
$\pi_M $.

Conversely, assume that $L$ is coisotropic with respect to $\pi_M$.
 Let $f \in {\mathcal F}(U)$ be a
 local function that vanishes on $L \cap U $.
 It is immediate that $L$ is also coisotropic with respect to  $
 \schouten{\pi_M, f}_{TM}$.
By Eq. (\ref{eq:2.32}), the identity $  \schouten{\pi_M, f}_{TM} = (-1)^{k+1}\rho
\big( \delta(f) \big) $ holds, so that, according to Lemma
 \ref{lem:coiso} (2), we have $ \delta (f) \in \left\langle B\right\rangle_U $.

At this point, to show that $\delta(\left\langle B\right\rangle_U)  \subset \left\langle B\right\rangle_U $, it suffices to 
show that $ \delta (b) \in \left\langle B\right\rangle_U$ for any section $b \in \left\langle B\right\rangle_U $.
According to Lemma \ref{lem:coiso} (2), it suffices indeed
to prove that $ \rho( \delta(b)) [f_1,\cdots,f_{k}]  $
vanishes on $L \cap U $ if $f_1,\cdots,f_{k} \in {\mathcal F}(U)$ 
vanish on $L \cap U $. Let us prove this point
  $$\begin{array}{rcll} \rho \left( \delta(b) \right) [f_1,\cdots,f_{k}] &= &
 \schouten{\rho(\delta(b)), f_1}_{TM}  [f_2,\cdots,f_{k}] & 
 \\ &= & \rho \left( \schouten{\delta (b),f_1}_A  \right) [f_2 , \cdots , f_{k}]
   & \\ &  =&    \rho \left( \delta(\schouten{b,f_1}_A ) \right) [f_2 , \cdots , f_{k}] &\\
   & & -    \rho \left( \schouten{b,\delta(f_1)}_A  \right) [f_2 , \cdots , f_{k}]   & \mbox{ by Eq. (\ref{k-dif-prop}) }  \\
   \\
\end{array} $$
Now,  $ \schouten{b,f_1}_A = \rho(b) [f_1]$ is a function that vanishes on $ L\cap U$. Hence $\delta(\schouten{b,f_1}_A)$ is a section in $\left\langle B\right\rangle_U $,
so that $L$ is coisotropic with respect to $ \rho \left(
\delta(\schouten{b,f_1}_A ) \right) $ by Lemma \ref{lem:coiso}, and the function
  $$  \rho \left( \delta(\schouten{b,f_1}_A ) \right) [f_2 , \cdots , f_{k}]$$
 vanishes on $L \cap  U$. Also, we have $\rho \left( \schouten{b,\delta(f_1)}_A \right) =  \schouten{\rho(b),\rho(\delta(f_1))}_{TM} $;
 which implies, since the Schouten bracket of coisotropic multi-vector fields is again coisotropic,
  that $L $ is coisotropic with respect to $\rho \left( \schouten{b,\delta(f_1)}_A \right)$, and the function
 $$   \rho \left( \schouten{b,\delta(f_1)}_A  \right) [f_2 , \cdots , f_{k}] $$
vanishes on $L \cap  U$.
These two last equations imply that the function  $\rho \left( \delta(b) \right) [f_1,\cdots,f_{k}]$
vanishes on $L \cap  U$. This completes the proof. 
\end{proof}

\smallskip

 (ii) {\bf From a multiplicative $k$-vector field on $\Gamma \toto M $ to a $k$-vector field on $Z(R) $.}

We prove in this section that, given a multiplicative $k$-vector field $\pi_\Gamma $ on 
the Lie groupoid $\Gamma \toto M $, the \'etale resolutions constructed in Proposition
\ref{prop:desing+groupoid} are automatically compatible with $\pi_M = s_* 
\pi_\Gamma$, provided that the chosen algebroid crossing, with respect to which they are constructed, is coisotropic with respect to $\pi_M$. Roughly speaking, the idea  is to construct explicitly the $k$-vector field on the resolution
$(Z(R), \phi)$ out of $\pi_\Gamma$ by the kind of reduction procedure described in Lemma \ref{lem:technicalities} below.

The following Theorem summarises what 
this section adds to Proposition \ref{prop:desing+groupoid}.

\begin{theo}\label{theo:compatible_resolutions}
Let $(A \to M,\rho,[\cdot,\cdot]) $ be a Lie algebroid,
${\mathcal S} $ a locally closed algebroid leaf and $L$
an algebroid crossing of $\overline{\mathcal S} $ with normalisation $B \to L $.
If
\begin{enumerate}
\item there exists a source-connected Hausdorff
 Lie groupoid $\Gamma \toto M $ integrating $(A \to M,\rho,[\cdot,\cdot])
 $, and
\item there exists a multiplicative $k$-vector field
$ \pi_{\Gamma}$, with $k \geq 2$,
 on $\Gamma \toto M$ such that $L$ is coisotropic with respect
  to $s_* \pi_\Gamma = \pi_M $, and
\item  there exists a  sub-Lie groupoid $R \toto L $ of $ \Gamma \toto N$,
  closed as  a subset of
  $\Gamma_L^L $, integrating
the algebroid $(B \to M,\rho,[\cdot,\cdot])  $,
\end{enumerate}
then
 \begin{enumerate}
\item $\pi_M $ is well-defined and tangent to ${\mathcal S} $.
\item $(Z(R),\phi,\pi_{Z(R)} )$ is an \'etale resolution of
  $\bar{\mathcal S} $
compatible with $ \pi_M$,
where \begin{enumerate}
 \item  $Z(R)=  R \backslash \Gamma_L$ and,
\item  $\phi: Z(R) \to M $
is the unique holomorphic/smooth map
such that  the following diagram commutes
\begin{equation}\label{eq:comdia1'}
 \xymatrix{
\Gamma_L   \ar[r]^{p} \ar[rd]^{t}  & Z(R) \ar@<0.5ex>[d]^{\phi} \\
  & M}
\end{equation}
where $p:\Gamma_L \to Z(R)= R\backslash \Gamma_{L} $ is the natural
projection, and
\item $\pi_{Z(R)} $ is the unique $k $-vector field on $M$ such that
for all local functions $\tilde{f}_1,\cdots, \tilde{f}_k \in {\mathcal
  F} (\Gamma) 
 $ and $f_1, \cdots,f_k \in {\mathcal
F}(Z(R)) $   the relations
 $$ \left\{ \begin{array}{ccc} 
p^* f_1 & = & \imath^*  \tilde{f}_1 \\
  &\vdots & \\ 
 p^* f_k &= &\imath^* \tilde{f}_k \\ 
\end{array} \right. $$
imply
\begin{equation}
  \label{eq:feedback} 
     p^* \pi_{Z(R)} [f_1,\cdots ,f_k] = \imath^* \pi_\Gamma 
[\tilde{f_1},\cdots, \tilde{f_n}]
\end{equation}
(wherever these identities make sense). Here $\imath  $ stands for the
inclusion map $\Gamma_L \hookrightarrow \Gamma $.
%
\item When $ (\Gamma \toto M , \pi_\Gamma)$ is a Poisson Lie groupoid,
id est $k=2$ and $\pi_\Gamma  $ is a multiplicative
Poisson bivector field on  $ \Gamma \toto M$, then $\pi_{Z(R)} $ is a Poisson
  bivector field and the resolution $(Z(R),\phi,\pi_{Z(R)}) $ 
is a Poisson resolution.
 \end{enumerate}
\item When  $L \cap {\mathcal S} $ is a $R$-connected set,
this \'etale resolution is a covering resolution with typical fiber
$\frac{\pi_0(I_x(\Gamma) ) }{\pi_0(I_x(R)) } $, where $x \in
{\mathcal S} $ is an arbitrary point, and $I_x(\Gamma) $ (resp.
$I_x(R) $) stands for the isotropy group of $ \Gamma \toto M$ (resp.
of $R\toto L$) at the point $x$.

\item This \'etale resolution is a resolution if
and only if  $R=\overline{ \Gamma_{L \cap {\mathcal S}}^{L
    \cap {\mathcal S}}} \cap \Gamma_L^L$.

\item When $L \cap {\mathcal S} $ is a connected set
and $ R \toto L $ is a source-connected sub-Lie groupoid of $\Gamma
\toto M $, then the typical fiber is isomorphic to  $\frac{\pi_1(
{\mathcal S} ) }{j(\pi_1( L \cap {\mathcal S}  ))} $, where $j$ is
the map induced at the fundamental group level by the inclusion of
$L \cap {\mathcal S} $ into ${\mathcal S} $.
\end{enumerate}
\end{theo}

We postpone until the end of this section the proof of  Lemma \ref{lem:technicalities} below,
which gives us the general frame to operate  reduction of multi-vector fields.
For a submanifold $N $ of a manifold $ Q $ and a point $x \in N$,
 $T_xN^{\perp} \subset T_x^* Q $
stands for $\{\alpha \in T_x^* Q \, |  \, T_x N \subset {\rm ker}(\alpha)\}$.

\begin{lem}\label{lem:technicalities}
Let $Q $ be a manifold,  $\imath : N \hookrightarrow Q$ a submanifold of $Q$
and  $\Phi:N \to P $ a surjective submersion with connected fibers.
Let $\pi_Q $ be a $k$-vector field on $M$.
If
\begin{enumerate}
\item  $(\pi_Q)_{|_x} \left[ \wedge^{k-1} ({\rm ker} (\diff_x \Phi))^{\perp} \wedge
  T_x N^{\perp} \right] =0$ and,
\item for any $n \in N $ and $u \in T_n N $, there exists a vector
  field $X$ on $Q$ tangent to 
$N$, whose value at $n$ is $u$, and such that 
 $(L_X \pi_Q)_{|_x} [ \wedge^{k}  ({\rm ker} (\diff_x \Phi))^{\perp}] =0 $.
\end{enumerate}
then there exists a unique $k $-vector field 
$\pi_P $ on $P$ such that, for any local functions $\tilde{f}_1,\cdots, \tilde{f}_k $ in 
${\mathcal F}(Q)$ on $N$ and $ f_1,\cdots ,f_k  $ in 
${\mathcal F}(P)$, the relations
 $$ \left\{ \begin{array}{ccc} 
\Phi^* f_1 & = & \imath^*  \tilde{f}_1 \\
  &\vdots & \\ 
 \Phi^* f_k &= &\imath^* \tilde{f}_k \\ 
\end{array} \right. $$
imply
\begin{equation}
  \label{eq:feedback2} 
     \Phi^* \pi_P (f_1,\cdots ,f_k) = (-1)^{k+1} \imath^* \pi_Q [\tilde{f_1},\cdots, \tilde{f_n}]
\end{equation}
(wherever these identities make sense).
\end{lem}

We now prove Theorem \ref{theo:compatible_resolutions}.
\begin{proof}
Theorem \ref{theo:compatible_resolutions}(1) is nothing but the
statement of Proposition \ref{prop:tangency}.
Points (3)-(5) are a consequence of Theorem \ref{theo:compatible_resolutions}(2)
and Proposition \ref{prop:desing+groupoid}(2)-(4). 
It remains us, therefore, the task of proving Theorem \ref{theo:compatible_resolutions}(2).
Now, when the assumptions of Theorem \ref{theo:compatible_resolutions}(2) are satisfied
for a sub-Lie groupoid $R \toto L $ of $\Gamma \toto M $ integrating $B \to L $,
all the assumptions of Proposition \ref{prop:desing+groupoid}(1) are satisfied as well.
Therefore, it only remains us the task  of proving that
the multivector field $\pi_{Z (R)} $ that appears in Eq. (\ref{eq:feedback}) can 
indeed be constructed.
  
To start with, we consider $R_0 \toto L $  the connected component of the
identities of the groupoid $R \toto L $. Again, $R_0 \toto L $
is a sub-Lie groupoid of $\Gamma \toto M $ that integrates $B \to L $
and which is closed in $ \Gamma_L^L$. 

Let us check that the assumptions of Lemma \ref{lem:technicalities} are satisfied
with $Q := \Gamma, N := \Gamma_L, P:= Z(R_0), \pi_Q:= \pi_{\Gamma}$
 and $ \Phi:= p$ the natural
projection $\Gamma_L \to Z(R_0) =R_0 \backslash \Gamma_L $.
First,  since $R_0 \toto L$ is source-connected,
 $ p$ is a surjective submersion with connected fibers.
 The two points below show that the assumptions 1 and 2 in Lemma \ref{lem:technicalities}
 are satisfied.

\begin{enumerate}
\item The identity 
 \begin{equation}\label{eq:coiso} (\pi_\Gamma)_{|_\gamma} \left[
   \wedge^{k-1}  ({\rm ker } (\diff_\gamma p))^{\perp}  \wedge
   (T_{\gamma}\Gamma_L)^{\perp} \right]=0 \end{equation}
holds for any $\gamma \in \Gamma_{L } $.
Since $s^*: T_m L^{\perp} \to T_\gamma \Gamma_L^{\perp} $
is one-to-one, for any $\alpha \in (T_\gamma \Gamma_L)^{\perp} $,
 there exists a local function $f  $ on $M$, defined in a neighbourhood
of $s(\gamma) $ and vanishing on $L$,
 such that $\diff_\gamma s^* f = s^* \diff_{s(\gamma)} f =
\alpha $. Now, Eq. (\ref{eq:rightmult}) gives
  $$   \schouten{\pi_\Gamma , s^* f }_{T\Gamma} =  \schouten{ \pi_\Gamma , \overrightarrow{f} }_{T\Gamma}
  =   \overrightarrow{\delta(f) }. $$
By assumption, $L$ is coisotropic with respect to $\pi_M $,
and it follows from  Proposition \ref{prop:cois-ckvector}
that $  \delta(f) $ is a section of $\left\langle B\right\rangle_U $.
In particular, $\overrightarrow{\delta(f) }_{\gamma} $ lies in the
  ideal of $\wedge^\bullet T_\gamma \Gamma $ 
(with respect to the wedge product) generated by the space of elements of the form
$ \overrightarrow{b}_{|_\gamma}$ with $b \in B_{s(\gamma)} $. 
But this space is, by definition, the foliation given by the  left action of  $R_0 \toto L $ on $\Gamma_L $,
and coincides therefore with the kernel of $ {\rm d}_\gamma p$.
This amounts to the fact that
   $ \overrightarrow{\delta(f) }_{|_\gamma} (\beta) =0$ for any $\beta
  \in \wedge^{k-1} ({\rm ker } (\diff_\gamma p))^{\perp}  $.
Hence $(\pi_\Gamma)_{|_\gamma} (\alpha \wedge \beta) =0 $.
This proves Eq. (\ref{eq:coiso})
and Condition 1 in Lemma \ref{lem:technicalities} is satisfied.
\item For any $\gamma \in T_\gamma \Gamma_L $  and any 
$u \in T_\gamma \Gamma_L $, there exists a section of $b \in \left\langle B\right\rangle_U $,
defined on a neighbourhood $U$ of $s(\gamma) $ such that $u =
\overrightarrow{b}_{|_\gamma} $. 
Now, Eq. (\ref{eq:rightmult}) gives
  $$   L_{\overrightarrow{b}} \pi_\Gamma = \schouten{\pi_\Gamma,\overrightarrow{b} }_{T\Gamma} = (-1)^k\overrightarrow{\delta(b)}  
   $$
Since $  {\delta(b) }$ is a section of $\left\langle B\right\rangle_U $
  according to Proposition \ref{prop:cois-ckvector}, we have again
   $ \overrightarrow{\delta(b) }_{|_\gamma}(\beta ) =0 
$  for any $\beta  \in   \wedge^{k} ({\rm ker
  } (\diff_\gamma p))^{\perp} $ and Condition 2 in Lemma \ref{lem:technicalities} is satisfied.
\end{enumerate}

By Lemma \ref{lem:technicalities} therefore, the quotient space
$Z(R_0) = R_0
\backslash \Gamma_L $ inherits a $k$-vector field $\pi_{Z(R_0)} $ that satisfies
Eq. (\ref{eq:feedback}). 
Our next goal is to show that the triple $( Z(R_0)  , 
\phi_0 , \pi_{Z(R_0)} )  $ is a resolution compatible  with $\pi_M $. In short, we have to show 
that 
 \begin{equation}\label{eq:itisok} (\phi_0)_*  \pi_{Z(R_0)} = \pi_M .\end{equation}

Recall, from \ref{eq:proje},  that $(-1)^{k-1}  t_* \pi_\Gamma = \pi_M $, so that
$ (-1)^{k-1} \pi_\Gamma [t^*f_1,\cdots, t^*f_k] = t^* \pi_M
[f_1,\cdots,f_k]$.
Now, 
 Let $f_1 \cdots f_k \in
{\mathcal F}(M) $ be local functions.
For $i =1,\dots,k $, the restriction to $\Gamma_L $  of the
function $ t^* f_i $ is equal  to the pull-back
through
$p $ of the function $   \phi_0^* f_i  $. In equation:
$$ \left\{ \begin{array}{ccc} 
p^*  (\phi_0^* f_1)& = & \imath^*  (t^ * f_1) \\
  &\vdots & \\ 
 p^*  (\phi_0^* f_k) &= &\imath^* (t^ * f_k)  \\ 
\end{array} \right. $$
Hence Lemma \ref{lem:technicalities} (more precisely Eq.
(\ref{eq:feedback2}))
 yields
 \begin{equation}\label{eq:target_m} (-1)^{k-1}  p^* \pi_{Z(R_0)}[  \phi_0^* f_1,\cdots,
\phi_0^* f_k ] = \pi_{\Gamma}\big[ t^* f_1 ,\cdots, t^*
f_k \big]_{|_{\Gamma_L}} \end{equation}
According to Eq. (\ref{eq:proje}),
 the relation $ \pi_{\Gamma}\big[ t^* f_1 ,\cdots, t^* f_k
\big] |_{\Gamma_L} = (-1)^{k-1} t^* (\pi_M [f_1,\cdots,f_k]) $
holds. Together with Eq. (\ref{eq:target_m}), this
 implies in turn:
\begin{eqnarray*} (-1)^{k-1} p^* \pi_{Z(R_0)}[  \phi_0^* f_1,\cdots,\phi_0^* f_k ] &=&
 (-1)^{k-1} t^* \pi_M [f_1,\cdots,f_k] \\ &=&
 (-1)^{k-1} p^* (\phi_0^* \pi_M [f_1,\cdots,f_k] ).\end{eqnarray*}
Since $p$ is a surjective submersion, this amounts to
$$ \pi_{Z(R_0)}[  \phi_0^* f_1,\cdots,\phi_0^* f_k ] = p^* \circ \phi_0^* \pi_M [f_1,\cdots,f_k] ,$$
and this completes the proof of the relation (\ref{eq:itisok}). In conclusion $(Z(R_0),\phi_0) $ is a \'etale resolution
compatible with $\pi_M $.

Now, the natural inclusion $R_0 \subset R $, where $R \toto L$ is a sub-Lie groupoid
of $\Gamma \toto M $ integrating the normalisation $B \to L $ of the
algebroid crossing $L$ induces a morphism of \'etale resolution 
$$ Z(R_0) = R_0 \backslash \Gamma_L \to  R \backslash \Gamma_L = Z(R),$$
which, moreover, is a local diffeomorphism. 
According to Lemma \ref{lem:etalemorph}, the $k$-vector field $\pi_{Z(R_0)}$ goes 
to the quotient and defines a $k$-vector field $\pi_{Z(R)} $ on $Z(R) $ which
satisfies itself $\phi_* \pi_{Z(R)} = \pi_M $.
This completes the proof.
\end{proof}

\begin{rmk}
It is natural to ask what happens for $k=1$.
For  a multiplicative vector field $X $ on $\Gamma \toto M $, the
vector field $s_* X $ is well defined but
does not need to be tangent to ${\mathcal S} $. But if we assume that
it is tangent to ${\mathcal S}$,
then to require an algebroid crossing to be coisotropic with  respect
to $s_* X$ means that $s_* X $ must be also tangent to $L$.
Under these conditions,  one can define by the same method a vector field
on $Z(R) $ that projects onto $ s_* X $.

For $k=0 $, the situation is as follows.
A multiplicative function on $\Gamma \toto M $ is simply a function that satisfies
$f(\gamma_1 \gamma_2) = f(\gamma_1)+ f(\gamma_2) $ for all compatible
$\gamma_1,\gamma_2  \in \Gamma$.
Such a function does not induce any particular object on $M$.
However, when $f$ vanishes on $R $,  it goes to the quotient 
to yield a function on $Z(R) $.
\end{rmk}

We finish this section with a proof of Lemma \ref{lem:technicalities}.

\begin{proof}
Choose arbitrary local functions ${f}_1,\cdots, {f}_k $ defined 
on an open subset $U \subset P$. In a neighbourhood $V\subset Q$ of any $ m \in \phi^{-1}(m)$,
there exists functions  $\tilde{f}_1,\cdots, \tilde{f}_k $ 
such that:
 $$ \left\{ \begin{array}{ccc} 
\phi^* f_1 &  =    & \imath^*  \tilde{f}_1 \\
           &\vdots & \\ 
\phi^* f_k &    =  &\imath^* \tilde{f}_k \\ 
\end{array} \right. $$
We want to define a $k$-vector field $\pi_P$  on $P$ by 
  $$ (\pi_P)_{|_p} [f_1,\cdots,f_k]  =   (\pi_{Q})_{|_x} [ \tilde{f}_1,\cdots, \tilde{f}_k ] , $$
  where $x \in N$ is any point with $p(x) = p$.
We have to check that this definition makes sense:
for this purpose it suffices to check that {\em (i)} the  restriction to $N $ of  $\pi_{Q} [ \tilde{f}_1,\cdots,
  \tilde{f}_k ]$ 
does not depend on the choice of the local functions
  $\tilde{f}_1,\cdots,\tilde{f}_k$,
and {\em (ii)}  that  the
  restriction to $N $ of  $\pi_{Q} [ \tilde{f}_1,\cdots,
  \tilde{f}_k ]$ is constant along the fibers of $\Phi: N \to P $ (i.e. does not depend on $x$).

First, we check  {\em (i)}. If we assume that
 $ \tilde{f}_k$ vanishes on $N \cap V$, then we have
  $$ \diff_x \tilde{f}_1 \wedge \cdots \wedge \diff_x \tilde{f}_k  \in 
  \wedge^{k-1} ({\rm ker} (\diff_x \Phi))^{\perp} \wedge
  T_x N^{\perp}.$$
As a consequence, the restriction to $ N$ of the 
function $ \pi_{Q} [ \tilde{f}_1,\cdots,   \tilde{f}_k ]   $
vanishes; this  proves {\em (i)}.

Now, we prove {\em (ii)}. The fibers of $\Phi: N \to P $ being connected
subsets, it suffices to check that 
for any   $x \in V \cap p^{-1}(U)  $ and $u \in T_{x} N $, 
 we have   $ u [\pi_{Q} [ \tilde{f}_1,\cdots, \tilde{f}_k ]] = 0 $.
By assumption, there exists a local vector field $X$ through $u $
tangent to $N$
and which satisfies  $(L_X \pi_Q)_{|_x}  \in   \wedge^{k} ({\rm ker} (\diff_x \Phi))^{\perp} \wedge
  T_x N^{\perp} $. Hence
 $$ \begin{array}{rcl}  u \left[\pi_{Q} \left[ \tilde{f}_1,\cdots, \tilde{f}_k
 \right] \right]_{|_x} &=&     X \left[ \pi_{Q} \left[ \tilde{f}_1,\cdots, \tilde{f}_k
 \right] \right]_{ |_x} \\ & =&  (L_X \pi_{Q})_{ |_x} \left[ \tilde{f}_1,\cdots,\tilde{f}_k \right]
  - \sum_{i=1}^k  (\pi_{Q})_{|_x} \left[ \tilde{f}_1,\cdots, X[\tilde{f}_i] , \cdots , \tilde{f}_k \right]
  . \\  
  \end{array}   $$   
But the function $ (L_X\pi_{Q}) [ \tilde{f}_1,\cdots,
 \tilde{f}_k]  $ vanishes at the point $x $ by assumption since
$$ {\diff}_x \tilde{f}_1,\cdots,
 {\diff}_x \tilde{f}_k \in ({\rm ker}{\rm d}_x \Phi)^{\perp}.$$
Moreover, the restriction to $N $ of $X[\tilde{f}_i] $ vanishes,
hence, according to {\em (i)}, we have 
 $$(\pi_{Q})_{|_x} \left[ \tilde{f}_1,\cdots, X[\tilde{f}_i] ,
 \cdots \wedge \tilde{f}_k\right] =0,  \quad \forall i =1,\cdots, k .$$
Hence $u \big[\pi_{Q} [ \tilde{f}_1,\cdots, \tilde{f}_k
 ]\big] =0 $ (at the point $x$). This completes the proof of {\em (ii)}, and the proof of
 Lemma \ref{lem:technicalities} as well.
 \end{proof}

\section{Symplectic groupoids and symplectic resolutions.}
\label{sec:R3}

\subsection{Definition of a symplectic resolution.}
\label{sec:R3D3}

Let $(M,\pi_M) $ be a Poisson manifold, and $(T^*M \to M,
[\cdot,\cdot]^{\pi_M},\pi^{\#}_M) $ the  Lie algebroid associated with.
By classical theory of Poisson manifold, the leaves of this
algebroid are symplectic manifolds, and the inclusion maps are
Poisson maps. For a given leaf ${\mathcal S} $, we denote by
$\omega_{\mathcal S} $ its symplectic form. By construction,
 the bivector field $\pi_M$ on the manifold $M$ is tangent to  all the symplectic leaves.

We explain what we mean by a symplectic resolution.

\begin{defn} \label{def:comp_desing_sympl}
Let ${\mathcal S }$  be a locally closed symplectic leaf, and let
$\bar{\mathcal S}$ be its  closure, a {\em symplectic resolution of $\bar{\mathcal S} $ }
is a triple $ (Z,\omega_Z,\phi)$ where
\begin{enumerate}
\item  $(Z,\phi) $ is a resolution of $\bar{\mathcal S} $,
\item the $2$-form $\omega_Z$  defined on $\phi^{-1}({\mathcal S}) $
by $\omega_Z := \phi^* \omega_{\mathcal S} $ extends to a holomorphic/smooth
symplectic $2$-form  on $Z$.
\end{enumerate}
When $(Z,\phi) $ is only an \'etale/covering resolution, then we speak of an {\em \'etale/covering symplectic resolution}.
For convenience, we often denote  a symplectic resolution as a triple $(Z,\phi \omega_Z) $.
\end{defn}

\begin{rmk} \label{rmk:symp_comp}
 We leave it to the reader to check that an \'etale symplectic resolution
 resolution $(Z,\phi, \omega_Z) $ is in particular 
 an \'etale resolution compatible with $ \pi_M$, and, indeed, an \'etale Poisson resolution.
 \end{rmk}

Let us enumerate for clarity all the conditions required in order to
have a symplectic resolution $(Z,\omega_Z,\phi)$ of the closure of a
symplectic leaf ${\mathcal S} $.

\begin{enumerate}
\item $Z$ is a manifold and $\omega_Z $ is a symplectic form with
  Poisson bivector $ \pi_Z$,
\item  $\phi:Z \to M$ is a holomorphic/smooth map from the manifold
$Z $ to the manifold $M$,
\item $\phi(Z)=\bar{{\mathcal S}}$,
\item $\phi^{-1}({\mathcal S})$ is dense in $Z$, 
\item  the restriction of $\phi$ to a map from $\phi^{-1}({\mathcal S})$ to ${\mathcal S}$  is a biholomorphism/diffeomorphism,
\item $\phi: Z \to M$  is a Poisson map.
\end{enumerate}

For \'etale/covering symplectic resolutions, the fifth point above needs
to be replaced by  ``the restriction of $\phi$ to
$\phi^{-1}({\mathcal S})$ is  an \'etale/covering map over
${\mathcal S}$''.

We finish this introductory section by defining complete symplectic resolutions.
Let $(M_i,\pi_i)$, $i=1,2$ be Poisson manifolds.
 Recall that a Poisson map $\phi: M_1 \to M_2$ is said to be {\em
   complete} if, for all open subset $U \subset M_2 $ and all $f \in
 {\mathcal F}(U) $,  
the flow starting at $m \in M_1 $ of the Hamiltonian vector
 field $ {\mathcal X}_{\phi^* f} $  is defined at
 the time $t=1$  as soon as
 the flow starting at $\phi(m) \in M_2 $  of the Hamiltonian vector
 field $ {\mathcal X}_{f} $  is defined at the time
$t=1 $. 


\begin{defn}
Let $(M,\pi_M)$ be a Poisson manifold and ${\mathcal S}$ a symplectic leaf.
An \'etale symplectic resolution $(Z,\phi,\omega_Z)$ of $\overline{\mathcal S}$ is said to be complete
if the map $\phi $ is a complete Poisson map from $(Z, \omega_Z^{-1} ) $ to $(M,\pi_M) $.
\end{defn}


\subsection{Symplectic groupoid and symplectic resolution of the closure of a symplectic leaf.}
\label{sec:R3G}

Let $(M,\pi_M) $ be a Poisson manifold and ${\mathcal S} $ a locally
closed symplectic leaf.
We introduce in this section the notion of Lagrangian crossing of
$ \overline{\mathcal S}$, which
is a particular case of algebroid crossing adapted to the construction of
symplectic resolutions. We then specialise 
Theorem \ref{theo:compatible_resolutions} to the case of Lagrangian
crossing of Poisson manifolds.

\begin{defn}\label{def:Lag_cross}
Let ${\mathcal S}$ be a symplectic leaf of a Poisson manifold $(M,\pi_M)$. 
A {\em Lagrangian crossing} of $\overline{\mathcal S}$ is a
submanifold $L$ of $M$ such that
\begin{enumerate}
\item $L \cap {\mathcal S}$ is  dense in $L$ and is a Lagrangian submanifold of ${\mathcal S}$
\item $L$ has a non-empty intersection
with all the symplectic leaves contained in $\bar{{\mathcal S}}$.
\end{enumerate}
\end{defn}

We immediately connect this notion to the notion of algebroid
crossing.

\begin{prop} \label{prop:LCtoCoiso}
Any  Lagrangian crossing $L$  is an algebroid crossing
of $(T^*M \to M,\pi^{\#}_M ,[\cdot,\cdot]^{\pi_M})$
with normalisation $TL^{\perp} \to L $ 
and is coisotropic with respect to the bivector field $\pi_M $.
\end{prop}
\begin{proof}
First, $L \cap {\mathcal S}$ being a Lagrangian submanifold of
  ${\mathcal S}$, the identity $\pi_M^{\#}(T^*_x L^{\perp}) = T_xL$
holds for any $x \in L \cap {\mathcal S}$. By density of $L \cap
{\mathcal S}$ in $L$, the inclusion $\pi_M^{\#}(T^*_x L^{\perp})
  \subset T_xL$ holds for any $x \in L$. As a consequence,
 $L$ is coisotropic with respect to $ \pi_M $.
The vector bundle  $TL^{\perp}\to L $ is a
is the normalisation of $L$.
 
 According to  Definition \ref{def:Lag_cross} (2),
$L$ intersects all the algebroid leaves contained in  $\bar{\mathcal
   S} $, since these leaves are precisely the symplectic leaves. As a consequence, 
   first, $L$   is an algebroid crossing of $\bar{\mathcal S}$, and, second, its normalisation is $TL^{\perp} \to L $.
\end{proof}

%

In view of Proposition \ref{prop:LCtoCoiso}, one can apply 
Theorem \ref{theo:compatible_resolutions} 
to the particular cases where algebroid crossings with normalisation are Lagrangian crossing
of Poisson manifolds.


\begin{theo} \label{theo:sympl_resol}
Let $(M,\pi_M) $ be a  Poisson manifold, ${\mathcal S} $  be a locally 
closed symplectic leaf, and $L$ a Lagrangian crossing of $\bar{\mathcal S} $. If,
\begin{enumerate}
\item there exists a symplectic Hausdorff Lie groupoid
 $(\Gamma \toto M ,\omega_{\Gamma})$ integrating the Poisson manifold $M$, and
\item  there exists a closed sub-Lie groupoid $R \toto L $ of $\Gamma \toto N$,  closed as  a subset of
  $\Gamma_L^L $, integrating $TL^{\perp} \to L  $,
\end{enumerate}
then
\begin{enumerate}
\item
$(Z(R),\omega_{Z(R)},\phi)  $ is an \'etale complete symplectic resolution of $\bar{\mathcal S} $, where  \begin{enumerate}
\item $Z(R)=  R \backslash \Gamma_L$,  and,
\item  $\phi: Z(R) \to M $
is the unique holomorphic/smooth map such that  the following diagram commutes
\begin{equation}\label{eq:comdia2}
{\xymatrix{
    \Gamma_L \ar[r]^{p} \ar[rd]^{t} &
  Z(R)  \ar[d]^{\phi}  \\  & M   }}
\end{equation}
where $p:\Gamma_L \to Z(L)= R\backslash \Gamma_{L} $ is the natural projection, and
\item  $\omega_{Z(R)} $  is the symplectic form defined by
   \begin{equation}\label{eq:consomega}  p^* \omega_{Z(R)} =  -\imath^* \omega_\Gamma . \end{equation}
\end{enumerate}
\item When   $L \cap {\mathcal S} $ is a $R$-connected set,
this \'etale complete symplectic resolution is a covering complete symplectic resolution with typical fiber  
$\frac{\pi_0( I_x ( \Gamma ) ) }{\pi_0(I_x(R)) } $, where $x \in {\mathcal S} $ is an arbitrary
point, and $I_x(\Gamma) $ (resp. $I_x(R) $) stands for the isotropy
group of $ \Gamma \toto M$ (resp. of $ R \toto L  $) at the point $x$.

\item This \'etale complete symplectic  resolution is a complete symplectic  resolution if
and only if
$R$ contains $ \Gamma_{L \cap {\mathcal S}}^{L
    \cap {\mathcal S}} $. In this case, we have
$ R = \overline{\Gamma_{L \cap {\mathcal S}}^{L
    \cap {\mathcal S}}} \cap \Gamma_L^L $. 

\item When  $L \cap {\mathcal S} $ is a connected set
and $ R \toto L $ is the source-connected sub-Lie groupoid of
$\Gamma \toto M $ with Lie algebroid $TL^{\perp} \to L $, then the
typical fiber is isomorphic to  $\frac{\pi_1( {\mathcal S} )
}{j(\pi_1( L \cap {\mathcal S}  ))} $, where $j:\pi_1( L \cap
{\mathcal
  S}) \to \pi_1( {\mathcal S})$ is the map induced at
the fundamental group level by the inclusion of $L \cap {\mathcal S}
$ into ${\mathcal S} $.
\end{enumerate}
\end{theo}

We need, first, a few lemmas.

\begin{lem} \label{label:complete}
For any  symplectic groupoid
$(\Gamma \toto M,\omega_\Gamma) $
that integrates a  Poisson manifold
$(M,\pi_M) $, the source map 
(resp. the target map) is a complete
map from $\Gamma =\pi_\Gamma =\omega_\Gamma^{-1} $ to $(m,\pi_M) $
(resp. $(M,-\pi_M) $).
\end{lem}
\begin{proof}
For the real case, we refer to \cite{CDW}.
In the complex case, note, to start with, that for any holomorphic
Poisson manifold $(N,\pi_N) $ where $\pi_N$ is a holomorphic Poisson
structure with real part $\pi^R_N $, and any holomorphic local function $f$,
the Hamiltonian  vector field ${\mathcal X}_f $ is twice the
Hamiltonian vector field (with respect $\pi_N^R $) of the real part of
$f$. As  a consequence a Poisson map between holomorphic Poisson
manifolds is complete if the induced Poisson map between their real
parts is complete.

Now, according to \cite{SX}, the real part of
$ \pi_\Gamma$ is precisely the Poisson bivector field 
associated to the symplectic structure  
integrating the real part of $\pi_M $. The source and target
maps are then complete maps since Lemma \ref{label:complete}
holds true in the real case.  
\end{proof}

\begin{lem} \label{lem:coiso_Poisson}
Let $P_1,P_2$ be two Poisson manifolds and $\phi:P_2 \to P_1$ be a
Poisson or anti-Poisson map which is a surjective submersion.
If $L$ is a coisotropic submanifold in $P_1$, then $\phi^{-1}(L)$ is coisotropic in $P_2$.
\end{lem}
\begin{proof}
We assume that $\phi$ is a Poisson map,
the anti-Poisson case being similar, up to a sign.
Since $ \phi$ is a submersion,
for any $x \in P_2 $, the dual $ ({\rm ker}(\diff_x \phi))^{\perp}$
of the kernel of $ d_x \phi$ is generated by
covectors of the forms $\diff_{x_2} (\phi^* f) $, with $ f \in {\mathcal
  F}(P_1)$. Since $\phi: P_2 \to P_1 $ is a Poisson map
 $$ \phi_* (\pi_{P_2}^{\#} (\diff_{x_2} \phi^* f)) =  \phi_* ({\mathcal
 X}_{\phi^* f}(x_2))  =   {\mathcal
 X}_{f}(\phi(x_2))  $$
Since $ L$  is coisotropic, the relation $  {\mathcal
 X}_{f}(\phi(x_2))  \in   T_{\phi(x_1)}L $ holds. Since $\phi $ is a
 submersion, we obtain
 $ \pi_{P_2}^{\#} (\diff_{x_2} \phi^* f )\in T_{x_2} (\phi^{-1}(L ))  $.
This completes the proof.
\end{proof}





Now we turn our attention to the proof 
of Theorem \ref{theo:sympl_resol}.

\begin{proof}
First, Conditions 1 and 2 in Theorem \ref{theo:sympl_resol} imply that
Conditions 1, 2 and 3 in Theorem \ref{theo:compatible_resolutions} are satisfied.
According to Theorem \ref{theo:compatible_resolutions} therefore,
there exists a bivector field $\pi_{Z(R)}  $ on $Z(R) $
such that  $(Z(R),\phi,\pi_{Z(R)}) $ is an \'etale resolution of $
\bar{\mathcal S} $ which is compatible with the bivector field
$\pi_M $. In view of Definition \ref{def:comp_desing_sympl},
 what remains to prove is that, indeed:
{\em (i)} the bivector field $\pi_{Z(R)} $
 is the Poisson bivector field associated to a symplectic structure
$\omega_{Z(R)} $ on $Z(R)$ that satisfies Eq. (\ref{eq:consomega}), and {\em (ii)} that $\phi $ is a complete map.

We prove {\em (i)}. Denote by $\pi_{\Gamma}$ the multiplicative bivector field
on $\Gamma $ associated with the symplectic structure $\omega_\Gamma $. 
The submanifold $L$ being coisotropic in $M $ and the source map being an
Poisson map, $\Gamma_L $ is a coisotropic submanifold of the
symplectic manifold $(\Gamma,\omega_{\Gamma})$.

Since the $2$-differential associated with $\pi_\Gamma $ is the de Rham differential
 (see \cite{ILX} for instance),  Eq. (\ref{eq:rightmult}) amounts to the following relation  
  \begin{equation}\label{eq:david} \overrightarrow{\diff f} = \schouten{\pi_{\Gamma},s^* f}_{T\Gamma}= {\mathcal X}_{s^* f} = \pi_\Gamma^{\#} (s^ * \diff f) \end{equation}
for any local
function $f \in {\mathcal F}(M)$ (where ${\mathcal X}_f$ stands for
the Hamiltonian vector field of $f$).
One can immediately rewrite Eq. (\ref{eq:david}) as
 \begin{equation}\label{eq:david'} \overrightarrow{\alpha}  =  (\pi_\Gamma^{\#})_\gamma (s^* \alpha)  \end{equation}
for all $\alpha \in T_m \Gamma $ and $\gamma \in \Gamma_m $.
Since $s^* $ is a diffeomorphism from $ (T_mL)^{\perp} $ to $ (T_\gamma\Gamma_L)^{\perp} $, 
Eq. (\ref{eq:david'}) gives:
   $$   \left\{ \overrightarrow{\alpha}_{|_\gamma}, \alpha \in 
       T_m L^{\perp}  \right\}= ((\pi_{Z(R)}^{\#})_{|_\gamma} T_\gamma\Gamma_L)^{\perp} $$
Now, the kernel of the projection map $p: \Gamma_L \to Z(R) = R \backslash \Gamma_L $,
at a point $\gamma \in \Gamma_L $, consists precisely 
of the left-hand term in the previous expression, so that we have:
 \begin{equation}\label{lastf} ((\pi_\Gamma^{\#})_{|_\gamma})^{-1}  ( {\rm ker}( \diff_\gamma p) ) =  (T_\gamma\Gamma_L)^{\perp} \end{equation}

Now, for any $z \in Z(R)$ and any $\alpha  \in T^*_z Z(R)$, it
follows by classical bilinear algebra from Eq. (\ref{eq:feedback})
that the tangent vector $(\pi_{Z(R)}^{\#})_{|_z}(\alpha)$ is by construction
given by
\begin{equation}\label{eq:conspi} (\pi_{Z(R)}^{\#})_{|_z}(\alpha) = -({\rm d}_\gamma p)\big( (\pi^{\#}_{\Gamma})_{|_\gamma} (\widetilde{ p^*
  \alpha})\big), \end{equation}
where   $\gamma$ is a point in  $ p^{-1}(z)$  and,  $\widetilde{ p^*  \alpha} \in T_{\gamma}^* \Gamma $ is 
 a covector whose restriction to $T_{\gamma} \Gamma_L $ coincides with $  p^*  \alpha  $. 
In particular, by Eq. (\ref{lastf}), if $(\pi_{Z(R)}^{\#})_{|_z}(\alpha) =0 $, then $(\widetilde{ p^* \alpha}) $
belongs to $ (T_\gamma \Gamma_L)^\perp $, and $ \alpha$ needs to vanish.
In other words, $(\pi_{Z(R)}^\#)_{|_z} $ is an injective map, which implies that 
$ \pi_{Z(R)}$ is the Poisson bivector field of a symplectic structure $\omega_{Z(R)} $. 
Eq. (\ref{eq:conspi}) amounts to Eq. (\ref{eq:consomega}). 
This completes the proof of {\em (i)}.

Next, we prove {\em (ii)}. Let $U \subset M$ be a open subset, $f \in {\mathcal F}(U)$
a function such that the flow $\Phi^M_\tau $ starting at $ m \in U$ is defined for the time $\tau =1$,
 and let $z \in \phi^{-1}(m) $ be a point.
The target map $t$ from $( \Gamma,-\pi_{\Gamma}) $ to $( M ,\pi_M) $ is a complete Poisson map,
so that the flow $\Phi_\tau^\Gamma $ of $t^* f $ starting at $\gamma $ is defined for $ \tau=1$,
where $\gamma \in \Gamma $ is any point such that $p(\gamma) =z $.
Since the Hamiltonian vector field ${\mathcal X}_{t^* f}$ is tangent to $\Gamma_L$ and since $p_* {\mathcal X}_{t^* f} = {\mathcal X}_{\phi^* f} $,
the flow starting at $z$  of $ {\mathcal X}_{\phi^* f}$ is equal to $ p \circ \Phi^\Gamma_\tau $.
In particular, it is defined for $\tau =1 $. This completes the proof of {\em (ii)}.
 \end{proof}


\begin{rmk}
We are redevable to Jiang-Hua Lu for the following remark. If $R \toto L$ is a source-connected Lie groupoid, then, by
construction, the procedure  that we have used to build the symplectic structure on $Z(R) $ out of the symplectic structure of
$\Gamma$ matches exactly the procedure called symplectic reduction \cite{OR}
with respect to the coisotropic submanifold~$\Gamma_L$. 
\end{rmk}




Proposition \ref{prop:repr_isom} can then be adapted easily.

\begin{prop}\label{prop:repr_isom2}
 Let $L_i $, $ i=1,2$ be two Lagrangian crossing of $\bar{\mathcal S} $ such that $R_i = \overline{ \Gamma_{L_i \cap
    {\mathcal S}}^{L_i    \cap {\mathcal S}}} \cap \Gamma_{L_i}^{L_i}$ is a sub-Lie groupoid
of $\Gamma \toto M $. Let $Z_i $, $i=1,2 $  be the  resolutions
corresponding to by Theorem  \ref{theo:sympl_resol}(3), i.e. $Z_i = R_i \backslash \Gamma_{L_i} $.
The following  are equivalent:
\begin{enumerate}
\item[(i)] the symplectic resolutions $(Z_1,\phi_1,\omega_{Z_1}) $ and $(Z_2,\phi_2,\omega_{Z_2}) $  are isomorphic,
\item[(ii)]  $  \overline{ \Gamma_{L_1 \cap {\mathcal S}}^{L_2    \cap {\mathcal
  S}}} \cap \Gamma_{L_1}^{L_2} $
  is a Lagrangian submanifold of $ \Gamma$, and the restrictions
to this submanifold of the source and the  target
  maps are surjective submersions onto $L_1 $ and $L_2 $ respectively,
 \item[(iii)] there exists a submanifold $I$ of $\Gamma $ that gives a
     Morita equivalence
 between the Lie groupoids $R_1 \toto L_1 $ and $R_2 \toto L_2 $:
  $$  \xymatrix{ \Gamma_{L_1} \ar[rd] & R_1 \dar[d]  & \ar[ld]^s I
  \ar[rd]_t & R_2\dar[d] & \ar[ld] \Gamma_{L_2}  \\ & L_1 & & L_2 &
  } $$ 
  In this case moreover,  the $R_1 $-module $\Gamma_{L_1} $
 corresponds to the $R_2 $-module $\Gamma_{L_2} $ with respect to the Morita equivalence $I$. 
Also, the submanifold $I \subset \Gamma $ is Lagrangian in $\Gamma $. 
\end{enumerate}
\end{prop}
\begin{proof}
Two symplectic resolutions isomorphic as resolutions are isomorphic as symplectic resolutions.  
The result is just then an immediate consequence of Proposition \ref{prop:repr_isom}.

The only difficulty is to prove the last assertion.
The manifold $\Gamma_{L_1 \cap
    {\mathcal S}}^{L_2    \cap {\mathcal S}} $ is coisotropic by Lemma
\ref{lem:coiso_Poisson} and is therefore, indeed, 
Lagrangian because its dimension is half  the dimension of $\Gamma$.
 Hence $ I = \overline{ \Gamma_{L_1 \cap
    {\mathcal S}}^{L_2    \cap {\mathcal S}}} \cap \Gamma_{L_1}^{L_2}$
    is a Lagrangian submanifold.
\end{proof}
%


\section{Characterisation of symplectic resolutions of the previous  type.}
\label{sec:char}

We restrict our attention, in this section, to the most interesting case, id est, the case of symplectic resolutions.
The aim of the present section is to characterise  proper symplectic
resolutions   of  $\overline{\mathcal S} $ 
isomorphic to the symplectic  resolutions
of the form  $(Z(R)=R\backslash \Gamma_L,\phi,\omega_{Z(R)})
$ constructed out of a sub-Lie groupoid $R \toto L $ 
integrating a Lie algebroid crossing $L $ 
as in Theorem \ref{theo:sympl_resol} (3).
Assume that we have a Lagrangian crossing $L$ of the closure $\overline{\mathcal S} $ of a locally closed symplectic leaf 
${\mathcal S} $ of an integrable complex or real Poisson manifold $(M,\pi_M) $.
%

%
\begin{defn}\label{def:compa}
Let $L$ be a Lagrangian crossing of the closure $\bar{\mathcal S}$ of  a locally closed 
symplectic leaf ${\mathcal S} $ of a Poisson manifold $(M,\pi_M)$.
A symplectic resolution $(Z,\phi, \omega_Z) $ of $\bar{\mathcal S} $ is said to be {\em $L$-compatible} if there exists a submanifold $L_Z $ of $Z$ such that the restriction of $ \phi$  to $L_Z $ is a biholomorphism/diffeomorphism onto $L$.
 \end{defn}

\begin{example} According to Lemma \ref{lem:j(L)}, the symplectic resolution $(Z(R)=R \backslash \Gamma_L ,\phi_{Z(R)},\omega_{Z(R)})
$ constructed as in Theorem \ref{theo:sympl_resol} (3) is $L$-compatible. In this case, we have $L_Z= j(L)$.
\end{example}
\begin{rmk} 
For any $L$-compatible symplectic resolution,  we  have $L_Z= \overline{\phi_Z^{-1}(L \cap {\mathcal S})}$.
In particular, the manifold $L_Z$ that appears in Definition \ref{def:compa} is unique.
Indeed,  we could characterise $L$-compatible symplectic resolution as being those such that $\overline{\phi_Z^{-1}(L \cap {\mathcal S})}$ is a submanifold of $Z$ to which the restriction of $\phi $ is a biholomorphism/diffeomorphism onto $L $.
\end{rmk}
A resolution $(Z,\phi )$ is said to be {\em proper} if the map $\phi$ is a proper map.
Note that a proper symplectic resolution  is always complete.
\begin{example}
The symplectic resolution constructed in Proposition \ref{prop:x2y2} is not proper, while
the Springer resolution of a Richardson orbit is.
\end{example}
\begin{theo}\label{theo:Linverse}
Let $L $ be a Lagrangian crossing of $\overline{\mathcal S} $, where ${\mathcal S} $ 
is a locally closed symplectic leaf ${\mathcal S} $
of a  Poisson manifold $(M,\pi_M)$. 

Let $(Z,\phi_{Z},\omega_{Z})  $ be a $L$-compatible complete symplectic resolution
of $\overline{\mathcal S} $.

Let $\Gamma \toto M$ be a source-simply connected and source-connected
symplectic groupoid that integrates the Poisson manifold $(M,\pi_M) $.
Then:

\begin{enumerate}
\item
 there exists a sub-Lie groupoid $R \toto L $ of $\Gamma  \toto M $ integrating $TL^{\perp} \to L $ closed in $\Gamma_L^L $
  and containing $\Gamma_{L \cap {\mathcal S}}^{L \cap {\mathcal S}} $;
\item 
the symplectic resolution  
$(Z(R)= R \backslash \Gamma_L , \phi_{Z(R)} ,\omega_{Z(R)}) $ (whose existence is granted by
Theorem \ref{theo:sympl_resol} (3))  is isomorphic
(as a symplectic
resolution) to an open subset of $(Z,\phi_{Z},\omega_{Z}) $;
\item if, moreover, the symplectic resolution  
$(Z(R)= R \backslash \Gamma_L , \phi_{Z(R)} ,\omega_{Z(R)}) $ is proper,
then the symplectic resolutions
$(Z,\phi_{Z},\omega_{Z}) $ and $(Z(R)= R \backslash \Gamma_L , \phi_{Z(R)}
,\omega_{Z(R)}) $ are isomorphic (as symplectic resolutions).
\end{enumerate}
\end{theo}

Before proving Theorem \ref{theo:Linverse}, we have to adapt in our context a
 result from the proof of Theorem 8 in \cite{CF2} and to state that $(Z,\phi_Z) $
 is a right $\Gamma$-module. 
There are important differences between our case and the setting
of \cite{CF2}. The authors of \cite{CF2} work with symplectic realizations, id est, symplectic varieties $(S,\omega_S) $ endowed with
 a surjective submersion from $S$ to $M$ which is also a Poisson map, while we do not assume $\phi $ to be surjective in general.
Moreover, the holomorphic case is not considered in their work. Also, we prefer to work with right action, while \cite{CF2} works with
left action, but this last point makes of course no major difference.
However, the following fact, adapted from the proof of Theorem 8 in \cite{CF2}, remains valid.

\begin{prop}\label{prop:cfCF2} 
Let $(\Gamma \toto M \omega_\Gamma) $ be a source-connected and source-simply connected
symplectic groupoid integrating a Poisson manifold $(M,\pi_M) $. 
Let ${\mathcal  S} $ be a locally closed symplectic leaf, and 
 $(Z , \phi_Z ) $ an \'etale symplectic resolution of $\overline{\mathcal  S} $.
 There is a unique action of the Lie groupoid $\Gamma \toto M $ on $Z \stackrel{\phi_Z}{\to} M $
 whose restriction  to $ \phi_Z^{-1} ( {\mathcal S})$ is given by 
  \begin{equation} \label{eq:rest_act}    \gamma \cdot z = \phi_Z^{-1} ( t(\gamma )) \quad \forall z  \in \phi_Z^{-1} ( {\mathcal S}), 
  \gamma \in \Gamma_{s(z)}.  \end{equation}
\end{prop}
\begin{example}
When $(Z= R \backslash \Gamma_L,\phi_Z)$ is a symplectic resolution associated to a sub-Lie groupoid $R\toto M$ integrating
a Lagrangian crossing $L$ as in Theorem \ref{theo:sympl_resol}(3), then the unique action that satisfies Equation (\ref{eq:rest_act})
is the action induced by the right action of $\Gamma$ to itself.
More precisely, it is given by  $\gamma \cdot [g] = [g  \gamma ]$ for all $[g] \in  R \backslash \Gamma_L, \gamma \in \Gamma $
with $s(\gamma) = \phi_Z ([g])$, where 
$[g'] \in R \backslash \Gamma_L $ stands for  the class in  $  R \backslash \Gamma_L$ of a given element $g' \in \Gamma_L $. 
\end{example}
\begin{proof}
Since $\phi_Z$ is an isomorphism from $\phi_Z^{-1}({\mathcal S})  $ to ${\mathcal S} $,
there is at most one Lie groupoid action of $\Gamma \toto M $ on $ Z \stackrel{\phi_Z}{\to} M$
that satisfies Eq. (\ref{eq:rest_act}) by density of $ \phi_Z^{-1}({\mathcal S}) $ in $Z$. 
This proves uniqueness.

We sketch the argument of the existence and prove in detail only those which
differ from \cite{CF2}.

Recall that a {\em cotangent path} is a map $a(u) $ of class ${\mathcal C}^1$
from $[0,1]$ to $T^*M $ satisfying
 $$ \frac{\diff m(u)}{\diff u}  = \pi_M^{\#} \big(a(u)\big)  \quad \forall u \in [0,1]$$
where $m(u)  $ is the base path of $a(u)$, id est, the projection of
$a(u)$ onto $M$ (through the canonical projection $T^*M \to M $).
Also, one assumes that $u \mapsto a(u) $ is equal to zero
for all $u$ in  neighbourhoods of $0$ and $1$, so that one can concatenate two $A$-paths 
when the end point of the first one coincides with the starting point of the second one.
 There is a notion of homotopy of cotangent paths, (see  \cite{CF2}), and it is now a classical result that if the
Poisson manifold
$(M,\pi_M)$ integrates to a source-connected and source-simply connected groupoid 
$\Gamma \toto M $, there is an isomorphism
  \begin{equation}\label{eq:Wein} \Gamma  \simeq  \frac{\mbox{cotangent paths}}{\mbox{homotopy}} .\end{equation}
The source and targets of the element  $\gamma \in \Gamma $
corresponding to a cotangent path are  the starting and end points of
its base path respectively. Product in $\Gamma $ corresponds to
concatenation of cotangent paths.
%

The argument used in the step 2 of the proof of theorem 8
in \cite{CF2} does not use the assumption that
what the map denoted by $\mu $ in \cite{CF2} (and which is our $\phi_Z$)
is a surjective submersion.
It therefore remains valid and yield the following result:
given a cotangent path $a(u) $ with base path $m(u) $,
  and some $z \in Z $ with $\phi_Z(z)=m(1) $, there exists a unique curve $z(u) $ on $Z$ 
with starting point  $z(0) =z $ and which satisfies
the differential equation
  \begin{equation} \label{eq:intcotg}
\left\{
  \begin{array}{lcrr}\frac{\diff z (u)}{\diff u} &=& \pi_{Z}^{\#} \big( \phi_{\Sigma}^* (a(u)\big) & \\ 
        \phi_{\Sigma} \big( \sigma(u)\big) &=&  m(u) & \, \, \, \forall u \in [0,1]  \\
\end{array}\right.
\end{equation}

Lemma 2 in  \cite{CF2} does not use the fact that
the map called $\mu $ in  \cite{CF2} is a surjective submersion
and remains valid. Its conclusion is that
  $z(1)$ does not depend on the class of homotopy
of the cotangent path $a(u) $. This second point, together 
with the identification given in Eq. (\ref{eq:Wein}),  
yields the existence of a map $Z \times_{\phi_Z,M,s}\Gamma  \to Z$.
This map defines a (right)-groupoid action (as shown in the end of Step 2
in \cite{CF2}, up to the fact that \cite{CF2} considers left action).
Eq. (\ref{eq:rest_act})  is  automatically satisfied.

%

The only delicate point is to show that this action is indeed
holomorphic in the complex case.
But it follows immediately from Eq. (\ref{eq:rest_act}) that the restriction to
$\phi_Z^{-1}({\mathcal S})\times_{\phi_Z,M,s}\Gamma_{ \mathcal S}  $ of the  action map $Z \times_{\phi_Z,M,s}\Gamma  \to Z $ is holomorphic.
Since $\phi_Z^{-1}({\mathcal S}) $ is dense is  $Z$, the action map is holomorphic.
\end{proof}
%
%

Now, we can turn our attention to the proof of Theorem 
\ref{theo:Linverse}.

\begin{proof}
1) We use the shorthand $j=(\phi_Z)_{ |_{Z}}$ to denote the restriction
of $\phi_Z $ to a biholomorphism/diffeomorphism from $L_Z $ onto $L$.
The action of $\Gamma \toto M $  on $Z \stackrel{\phi_Z}{\to}  M $ 
restricts and yields a map, that we denote  by $\Xi $, 
from   $    L_Z \times_{j,L,s} \Gamma_L $ to $Z  $.
But  $L_Z \times_{j,L,s} \Gamma_L  $ is simply isomorphic to $\Gamma_L
$ (the isomorphism being simply the projection onto the second
component),  so that $\Xi $ can be considered as a map from $\Gamma_L
$  to $Z $,
%

%
The relation (\ref{eq:rest_act}) can be rewritten as $\phi_Z \circ \Xi= t$. Under this form, it implies that the restriction
of $\Xi $ to $\Gamma_{L \cap {\mathcal S}} $ is given by $ \phi_Z^{-1} \circ t $. This fact has several consequences.
\begin{enumerate}
\item First, $\phi_Z^{-1}({\mathcal S}) \subset  \Xi(\Gamma_L) $.
\item Let  $\omega_{\Gamma_L} $ be the restriction of $\omega_\Gamma $ to $\Gamma_L $. 
The following relation holds:
  \begin{equation}\label{eq:pbreg}  \Xi^* \omega_Z = -\omega_{\Gamma_L} \end{equation}
 Let us prove this relation. Since the target map is an anti-Poisson map, we have
 the relation $t^* \omega_{\mathcal S} = (-\omega_{\Gamma}) |_{\Gamma^{\mathcal S}}  $.
 But $\Gamma_{L \cap {\mathcal S}} \subset  \Gamma^{\mathcal S}$, and we can conclude that
 $t^* \omega_{\mathcal S} =  -\omega_{\Gamma_L} $.
 This relation, together with the relation $ \Xi^* \circ \phi_Z^* = t^*$ and the fact that
 $phi_Z$ is a symplectomorphism from $(\phi_Z^{-1}({\mathcal S}), \omega_Z) $ to $({\mathcal S} , \omega{\mathcal S}) $
 implies that Eq. (\ref{eq:pbreg}) holds on $\phi_Z^{-1}({\mathcal S}) $, hence on $Z$
 by density. 

%
\item Let $ \epsilon(L) \subset \Gamma_L  $ be the image of $L$ through
  the unit map. The restriction of $ \Xi  $ to $\epsilon (L) $
takes its values in $L_Z $. Since both $t: \epsilon(L) \to L $
and $j:  $ are biholomorphisms/diffeomorphisms, the restriction of
  $\Xi $ to $\epsilon (L) $ is a  biholomorphism/diffeomorphism onto $L_Z $.
\end{enumerate}

 Let us show that $\Xi $ is a submersion onto its image.
 Recall that the rank of $\omega_{\Gamma_L} $ is
$dim({\mathcal S})$,
as well as the rank of $\omega_Z $.   Eq. (\ref{eq:pbreg})
implies then that
$\diff_\gamma \Xi $ is surjective at every point $\gamma  \in \Gamma_L $, id. est,
that $  \Xi $ is a submersion onto its image.

The inverse image of $L_Z $ through $\Xi$ is therefore a sub-manifold of $\Gamma_L $ that we denote by $R$.
Let us describe all the properties of $R$. First, by construction,  $R \subset \Gamma $ is the set of points in $r \in \Gamma $
such that $\Xi(r,m) \in L_Z $ where $ m = j^{-1}(s(r))$.  
With the help of that characterisation,
we leave it to the reader to check that $\Xi$ is stable by inverse.
Since $\Xi$ is a submersion onto its image, the restriction of $\Xi $ to $R$ is a
submersion onto its image $L_Z $, and the restriction of the target map $t $ to $R$ 
is also a submersion onto its image $L$. Since the inverse map intertwines the source and the target maps,
the restriction of the source map $s $ to $R$ 
is also a submersion onto ts image $L$.
Now, one sees easily that the product of two compatible elements in $R$is in $R$ also.
Therefore $R \toto L$ is
 a groupoid. Since it obviously contains
$\Gamma_{L \cap {\mathcal S}}^{L \cap {\mathcal S}} $
 as a dense open subset, it  proves (1).

2) According to Theorem \ref{theo:sympl_resol} (3), $ (Z(R) =R \backslash \Gamma_L, \phi_{Z(R)}  ,\omega_{Z(R)})  $
is a symplectic resolution of $\bar{\mathcal S}$. The map $\Xi $ goes to the quotient and yields a map, that we denote 
$\Psi $ from $ Z(R)  $ to $Z $. By construction, the restriction of $\Psi$ to the dense open subset
$\phi_{Z(R)}^{-1}({\mathcal S})  $ is given by $\Xi = \phi_Z^{-1}
\circ \phi_{Z(R)} $. In particular, it is a one-to-one map and it is a symplectomorphism.
By density, it implies that $\Psi$ is a symplectomorphism, and is therefore an open map, and 	a local diffeomorphism
onto its image.
But a local diffeomorphism which is a diffeomorphism on a dense open subset is a diffeomorphism.
This completes the proof.

%

3) It remains to prove that  $\Psi$ is onto when $\phi_{Z(R)} $ is proper. 
Let $U$ be an open subset of $Z $ contained in $\psi_Z^{-1}({\mathcal S}) $ 
whose closure $\overline{U} $ (w.r.t. the topology of 
the manifold $Z$) is compact. Let $V = {\Psi}^{-1}(U)$.
It is elementary that $V = \phi_{Z(R)}^{-1}( \phi_{Z} (U)) $,
so that $V \subset  \phi_{Z(R)}^{-1}( \phi_Z (\overline{U})  )$.
But  $\phi_Z (\overline{U}) $ is compact since $\phi_Z$
 is continuous and $ \phi_{Z(R)}^{-1}( \phi_Z (\overline{U})  ) $
is compact by properness of $ \phi_{Z(R)} $.
Now, ${\Psi}(\overline{V}) $ is compact and contains $U$, it
 therefore contains $\overline{U} $. Since any point in $Z$ lies inside
 the closure of a relatively compact open set contained in the dense
 open subset  $\psi_Z^{-1}({\mathcal S}) $,  ${\Psi} $
is surjective. This completes the proof.
\end{proof}
%
%


%
\section{Examples}
\label{sec:Ex}
\subsection{Poisson brackets on ${\mathbb R}^2 $ }
\label{sec:Ex_2D}

In this section,
we work in the setting of real
differential geometry.
We present an example for which
{ the space we are working on is a regular manifold, but the Poisson structure has  singularities.}

Let $\kappa (x,y)$ be a smooth non-negative function on ${\mathbb R}^2$.
We assume that the set of zeroes $(z_i)_{i \in I}$ is a discrete
subset of ${\mathbb R}^2$;
so that we can assume $I \subset {\mathbb N}$.
Define a Poisson bracket on ${\mathbb R}^2$ by
  $$ \{x,y\}= \kappa(x,y) $$

We focus on  the symplectic leaf ${\mathcal S} = {\mathbb R}^2 -
\{z_i,i \in I\} $. Note that in this case $\bar{\mathcal S} =
{\mathbb R}^2$  is a smooth manifold.

Any smooth submanifold $L$ of dimension  $1$, i.e. any smooth
embedded curve, is a coisotropic submanifold, whose intersection
with ${\mathcal S }$ is Lagrangian. So that any smooth embedded
curve $L$ that goes through all the points $z_i,i \in I $ is a
Lagrangian crossing. Such a curve always exists.

According to Corollary 5 in \cite{CF2}, this Poisson manifold is
integrable and integrates to a Lie groupoid $\Gamma \toto {\mathbb
R}^2 $. Theorem \ref{theo:sympl_resol}  allows us to built
symplectic \'etale resolutions whenever  $TL^{\perp} \to L $
integrates to a sub-Lie groupoid of $\Gamma \toto M $ closed in $\Gamma_L^L$.
In a future work, we shall see that this is always the case: more precisely there exists
a source-connected sub-Lie groupoid of $\Gamma \toto M $  as well as a
sub-Lie groupoid containing $\Gamma_{L \cap {\mathcal S} }^{L \cap {\mathcal S} } $
integrating $ TL^{\perp} \to L$ and closed in $\Gamma_L^L $.
However, we restrict ourself here to the most basic example.

\smallskip
{\bf Example: the case of the bracket $\{x,y\}=x^2 + y^2 $.}
\smallskip

We introduce complex coordinates $z =x+ iy $, $\bar{z}= x-iy $,
and we study the case $\kappa (x,y)= x^2 + y^2 = z \bar{z}  $.
According to \cite{Dazord}, the symplectic Lie groupoid is given in
this case by
\begin{enumerate}
\item $\Gamma := {\mathbb C} \times {\mathbb C} $,
\item the source map $ s(Z,z)= z$ and the target map $t(Z,z)= e^{Z
  \bar{z}} z $,
\item the product $(Z_1,z_1) \cdot (Z_2, e^{Z_1 \bar{z_1}} z_1 ) =
(Z_1 + e^{\bar{Z_1} z_1 } Z_2  , z_1 )  $
\item the symplectic structure
\begin{equation}\label{eq:omegaC2} \omega_{\Gamma} = z \bar{z} dZ \wedge d\bar{Z}+ 2 {\mathcal
R}{\rm e} (\overline{zZ}  dZ \wedge dz )
 +Z \bar{Z}   d\bar{z} \wedge dz  + 2 {\mathcal  R}{\rm e} (  dZ \wedge d\bar{z}
 ).\end{equation}
\end{enumerate}
(Indeed, the explicit structures of \cite{Dazord} have been slightly modify
in order to match our previous conventions).
The real axis is a Lagrangian crossing that we denote by $L$.
By construction,
 $\Gamma_L = \{(Z,\lambda) | Z \in {\mathbb C} , \lambda \in {\mathbb R}\} $.
The sub-Lie groupoid  $ R^{(0)} = {\mathbb R}^2 \subset \Gamma$
 integrates
$TL^{\perp} \to L $. So that all the assumptions of Theorem
\ref{theo:sympl_resol} are satisfied. We now describe explicitly the
symplectic \'etale resolution obtained by this procedure.

The quotient space
$Z(R^{(0)}) = R^{(0)} \backslash \Gamma_L $ is given
by
   $$ Z(R^{(0)}) \simeq \frac{ {\mathbb C} \times {\mathbb R }}{ \sim } $$
where $(Z,\lambda) \sim (Ze^{M \nu} + M ,\nu )      $ whenever
$\lambda = \nu e^{M \nu} $. 
In any class of the equivalence relation $\sim $, there is one and
   only one element of the form $ (i a,b) $ with $a,b  \in  {\mathbb
     R}$. Let us justify this point. Uniqueness is straightforward. 
      Now, we set 
$M=-({\rm Im}(Z))e^{-Z \lambda}$ and $ \nu = \lambda e^{-Z \lambda } $.
 One  checks easily that $t(M,\nu) = s(Z,\lambda)  $ and that
   $(M,\nu)(Z,\lambda) $ is of the requested form (i.e. the real part of the first component
vanishes).
 As a consequence, $Z(R^{(0)}) $ is simply isomorphic to ${\mathbb R}^2 $.
   The target map $t : \Gamma_L \to {\mathbb R}^2 \simeq {\mathbb C}$
factorizes   to yield the map
  $$  \phi (a,b):= b e^{iab} =  (b \,{\rm  cos}(ab), b \,{\rm  sin}(ab)) $$
  and $(Z(R^{(0)}),\phi)$ is an \'etale resolution of $\overline{\mathcal S}  \simeq {\mathbb R}^2$.
According to Theorem \ref{theo:sympl_resol}, $Z(R^{(0)}) $ endows a
symplectic structure $\omega_{Z(R^{(0)})} $ such that $\phi$ is a
Poisson map. We could compute this structure from the one on $\Gamma
$ as given by (\ref{eq:omegaC2}), but the computation is quite
tedious. It is much easier, since we know by Theorem
\ref{theo:sympl_resol}  that this structure has to exist, to deduce
directly from the fact that $\phi $ is Poisson map the explicit form
of the Poisson bracket $\{\cdot,\cdot\}_{Z(R^{(0)})} $ corresponding
to $\omega_{Z(R^{(0)})} $. We proceed as follows. The definition of a
Poisson map yields the relation
 \begin{equation}\label{eq:Poissonmap}   \{b \, {\rm  cos}(ab),   b \,{\rm  sin}(ab)  \}_{Z(R^{(0)})}
  =  \phi^* \{ x,y \}_{{\mathbb R}^2}  \end{equation}
where $x$, $y$ are the  coordinate functions of $M \simeq {\mathbb
  R}^2 $,
and where the relations $\phi^* x =b \, {\rm  cos}(ab)  $
and $\phi^* y =b \, {\rm  sin}(ab)  $ have been used.
 On the one hand, by a direct computation, we obtain
 \begin{equation}\label{eq:Poissonmap1}    \phi^* \{ x,y \}_{{\mathbb R}^2} = \phi^*  (x^2 +y^2) = b^2 ({\rm
   cos}^2(ab) +{\rm
   sin}^2(ab) ) = b^2   \end{equation}
and on the other hand the Leibniz rule gives 
  \begin{equation}\label{eq:Poissonmap2}   \{b \, {\rm  cos}(ab),   b \,{\rm  sin}(ab)  \}_{Z(R^{(0)})} = b^2 \{a,b  \}_{Z(R^{(0)})}  . \end{equation}
According to Equations. (\ref{eq:Poissonmap}-\ref{eq:Poissonmap1}-\ref{eq:Poissonmap2}),  the  induced structure on $
Z(R^{(0)})$
 is  simply given by $\{a,b\}_{Z(R^{(0)})}=1 $, which corresponds to the symplectic
 $2$-form
 $$ \omega_{Z(R^{(0)})} = {\rm d} a \wedge {\rm d} b $$


In order to find ``better'' resolutions, one has to replace $R^{(0)}
\toto L $ by  a bigger sub-Lie groupoid. This can be done as
follows. Fix $k \in {\mathbb N} $ and let
 \begin{eqnarray*} R^{(k)} &:=&\{ (\nu +  n \frac{i k \pi}{\lambda} ,\lambda) \, | \, \nu \in
 {\mathbb R}, \lambda \in   {\mathbb R}^*, n \in
 {\mathbb Z} \} \cup \{(\nu, 0) \,  | \, \nu \in {\mathbb R}\} \\ &=&  R^{(0)} \cup (\cup_{n \in {\mathbb Z}} \{ (\nu +  n \frac{i k \pi}{\lambda} ,\lambda) \, | \, \nu \in
 {\mathbb R}, \lambda \in   {\mathbb R}^* , n \in {\mathbb Z}^* \}   )\\ \end{eqnarray*}
For all $ k \in {\mathbb Z}$, $R^{(k)} $ is a closed sub-Lie groupoid of
$\Gamma \toto M $ with Lie algebroid $TL^{\perp} \to L$. For $k=0 $,
we recover the previous case, for $k=1 $, we have obviously
 $$ R^{(1)} = \overline{\Gamma_{L \cap {\mathcal S}}^{L \cap {\mathcal
 S}}} \cap \Gamma_L^L  .$$
By  Theorem \ref{theo:sympl_resol} (3) therefore, $ Z(R^{(1)}) = R^{(1)}
 \backslash \Gamma_L$
is a symplectic resolution of $\bar{\mathcal S} $. The quotient
space $Z(R^{(k)}) = R^{(k)} \backslash \Gamma_L $ is given by
   $$ Z(R^{(k)}) := \frac{ {\mathbb C} \times {\mathbb R }}{ \sim_k } $$
where $(Z,\lambda)
   \sim_k \big((-1)^{kn} Z e^{M \nu} + M + \frac{i k n\pi}{\nu},\nu \big)$
 whenever $\lambda = \nu (-1)^{nk} e^{M \nu} $.
 In any class of the relation $\sim_k $, there is at least
one element of the form $ (i a,b) $ with $a,b  \in  {\mathbb
     R}$, and a short computation shows that
 $Z(R^{(k)}) $ is indeed isomorphic to ${\mathbb R}^2/\sim_k' $,
where $\sim_k'$ is the equivalence relation that identifies $(a,b) $
with $((-1)^{kn} a+ \frac{nk \pi}{b},(-1)^{kn} b )  $ for all $n \in
{\mathbb Z} $ whenever $b \neq 0$.

The symplectic form $ \omega_{Z(R^{(0)})} = {\rm d} a \wedge {\rm d} b$
goes to the quotient through the equivalence relation $\sim_{k}' $,
as one easily checks. This achieves the construction of a symplectic
resolution of $\bar{\mathcal S} $.

We summarise this construction as follows, that presents the extremal cases $k=0$ and $k=1 $.

\begin{prop}\label{prop:x2y2}
\begin{enumerate}
\item
 The triple $(Z(R^{(0)}),\phi , \omega_{Z(R^{(0)})}) $ is an \'etale
symplectic resolution of $\bar{\mathcal S} $, more precisely:
 \begin{enumerate}
  \item $Z(R^{(0)}) \simeq {\mathbb R}^2$, we denote $(a,b)$ the canonical
  coordinates
  \item  $\omega_{Z(R^{(0)})} = {\rm d} a \wedge {\rm d} b$
 \item  $\phi (a,b):= b e^{iab} =  \big( b \,{\rm  cos}(ab), b \,{\rm  sin}(ab) \big) $
 \end{enumerate}
\item The triple $(Z(R^{(1)}),\phi_1 , \omega_{Z(R^{(1)})}) $ is a  symplectic
resolution of $\bar{\mathcal S} $, more precisely:
  \begin{enumerate}
  \item $Z(R^{(1)}) \simeq {\mathbb R}^2/\sim_1'$, where
 $ \sim_1'  $ is the equivalence relation that identifies
$(a,b) $ with $((-1)^{n} a+ \frac{n  \pi}{b},(-1)^{n} b )  $ for all
$n \in {\mathbb N} $ whenever $b \neq 0$. We denote by
$\overline{(a,b)} \in Z(R^{(1)}) $ the class of $(a,b) \in Z(R^{(0)}) $
  \item  $\omega_{Z(R^{(1)})} $ is the unique symplectic form satisfying
 $\Pi^* \omega_{Z(R^{(1)})}= \omega_{Z(R^{(0)})}$ where $\Pi: Z(R^{(0)}) \to
 Z(R^{(1)}) $
is the canonical projection
 \item  $\phi_1$ is (well)-defined by
 $\phi_1 \big(\overline{(a,b)}\big)=  (b \,{\rm  cos}(ab), b \,{\rm  sin}(ab)) $.
 \end{enumerate}
\item The projection $\Pi $ is a morphism of \'etale
symplectic resolutions.
\end{enumerate}
\end{prop}
Of course, the previous constructions could be done for any straight line through the origin, since
the Poisson structure  is invariant under a rotation centred at the origin.
Two such resolutions  are in general not isomorphic.
\begin{prop}
Let $L_1$ and $L_2$ be two straight lines through the origin.
Let $(Z_1,\phi_1)$ and $(Z_2,\phi_2)$ be the two symplectic resolutions associated with
as in Proposition \ref{prop:x2y2}(2). 
The resolutions $(Z_1,\phi_1)$ and $(Z_2,\phi_2)$ are isomorphic if and only if $L_1 =L_2 $.
\end{prop}
\begin{proof}
By symmetry, one can assume that $L_1$ is the real axis.
Then, by definition of the groupoid structure on $\Gamma \toto M $ 
 \begin{eqnarray*} \Gamma_{L_1 \cap {\mathcal S}}^{L_2 \cap {\mathcal S}} &=& \,  \{ \, (Z,\lambda ) \, | \, \lambda \in {\mathbb R}^*, \lambda e^{Z \lambda} \in L_2 \, \}  \\ & = & \,  \{ \, (Z,\lambda ) \, | \, \lambda \in {\mathbb R}^*,   e^{Z \lambda} \in L_2 \, \} \\
 &=&   \,  \{ \, (Z,\lambda ) \, | \,   \lambda {\rm Im}( Z)  - \alpha \, \in  \, \pi {\mathbb Z}  \,  \},  \\
 \end{eqnarray*} 
where   $\alpha \in {\mathbb R}$ is the angle with the horizontal line.
We leave it to the reader to check that $ \overline{\Gamma_{L_1 \cap {\mathcal S}}^{L_2 \cap {\mathcal S}}}  $
has an empty intersection with  $s^{-1}(0) $ when $\alpha $ does not belong to $  \pi {\mathbb Z}$.
(In other words, the latter means that there is no  sequences $(\lambda_k)_{k \in {\mathbb N}}, (Z_k)_{k \in {\mathbb N}}  $
with $\lambda_k  $ converging to $0$ and with $ Z_k$ convergent so that 
$ \lambda_k {\rm Im}( Z_k)  - \alpha \, \in  \, \pi {\mathbb Z}$ for all $k \in {\mathbb N}$).
In conclusion, when $L_2$ is not the horizontal axis, the restriction
of the source map $s$ to $ \overline{\Gamma_{L_1 \cap {\mathcal S}}^{L_2 \cap {\mathcal S}}}$ 
is not a surjective map onto $L_1$.   According to Proposition \ref{prop:repr_isom2} therefore,
the resolutions $(Z_1,\phi_1)$ and $(Z_2,\phi_2)$ can not be isomorphic.
\end{proof}

\subsection{Nilpotent orbits of a semi-simple Lie algebra: the
  Springer resolution}

\label{sec:Ex_Springer} We construct,  as a particular case of the
symplectic resolutions previously built,
 a very classical resolution called the {\em Springer
resolution}.

In the present case,  the space we are working on is a singular
  variety, but the Poisson structure is "as regular as possible"  in
  the sense that it is symplectic  at   regular points.

Let ${\mathfrak g}$ be a complex semi-simple Lie algebra and $G$  a
connected Lie group integrating ${\mathfrak
  g}$.
The Lie algebra ${\mathfrak g}$ is identified with its dual ${\mathfrak g}^* $
with the help of the Killing form.

We identify the linear Poisson manifold ${\mathfrak g}^*$ (endowed with the linear
Poisson structure) with
${\mathfrak g}$ with the help of the Killing form. This Poisson
manifold is integrable. More precisely, it integrates
 to the transformation Lie groupoid
$T^*G \simeq G \times {\mathfrak g} \toto {\mathfrak g}$ where
\begin{enumerate}
\item the source map, target map and product are as follows
 $$  \left\{ \begin{array}{ccc}  s(g,u) & =&u \\t(g,u) &= &  Ad_{g^{-1}} (u)
 \\ (g,u  ) \cdot (h, Ad_{g^{-1}}u) &= &( gh, u) \\
 \end{array} \right. $$
\item the symplectic structure is the canonical symplectic structure on
  a cotangent bundle.
\end{enumerate}
See \cite{CDW} or \cite{DeSilvaWeinstein} for more details.

Let ${\mathcal S}$ a nilpotent orbit in ${\mathfrak g}$. The
fundamental assumption that we have to make in order to
 construct an \'etale resolution of $\bar{\mathcal S} $
 is the following.
Assume that ${\mathcal S}$ is a {\em Richardson orbit}, i.e. that
there exists a parabolic subalgebra ${\mathfrak P} \subset
{\mathfrak g}$ whose nilradical ${\mathfrak U}$ satisfies the
property that ${\mathfrak U} \cap {\mathcal S}$ is dense in
${\mathfrak U}$.

\begin{lem}\label{lem:nilradical}
The nilradical ${\mathcal U}$ is a Lagrangian crossing of
$\bar{\mathcal S}$.
\end{lem}

I am strongly grateful to P. Tauvel and R. Yu for the following
proof.

\begin{proof}
According to the Theorem  of Richardson, (see \cite{TY}), for any $x \in
{\mathfrak U} \cap {\mathcal S}$, the intersection  of the $G$-orbit
of $x$ with ${\mathfrak P}$ is the $P$-orbit of $x $. As a
consequence, for any $v_1,v_2 \in T_x {\mathfrak U} \simeq
 {\mathfrak U}$, there exists
$p_1,p_2 \in  {\mathfrak P}$ such that $ v_i = [p_i,x]$, $i=1,2$. By
the definition of the symplectic structure $\omega_{\mathcal S} $ of
${\mathcal S} $, we have
$$\omega_{\mathcal S}(v_1,v_2)= \left\langle x,[p_1,p_2]\right\rangle  =\left\langle [x,p_1],p_2\right\rangle  .$$
The spaces  ${\mathfrak P} $ and  ${\mathfrak U}$ being dual to each
other w.r.t. the Killing form, this amounts to $ \omega_{\mathcal
S}(v_1,v_2)=  \left\langle [x,p_1],p_2\right\rangle  =0 $. In conclusion, ${\mathfrak U} $ is
a coisotropic submanifold of
 ${\mathfrak g} $.

Now, choose some $x \in {\mathfrak U} \cap {\mathcal S} $ and some
$w \in T_x{\mathcal S} $. There exists some $a \in {\mathfrak g} $
such that $w =[a,x] $.

Assume that $\omega_{\mathcal S}(v,w)=0 $ for all $v \in T_x
{\mathfrak U} $. Since $[{\mathfrak P},x] = {\mathfrak U}$, we have
  $$ \left\langle x,[p,a]\right\rangle  = \left\langle [x,p],a\right\rangle =0 \, \, \forall p \in {\mathfrak P} $$
Hence $a \in {\mathfrak U}^{\perp} = {\mathfrak P} $ and $v \in
T_x{\mathfrak U} $. Therefore ${\mathfrak U} \cap {\mathcal S} $ is
Lagrangian in ${\mathcal S} $.

The last delicate point is to check that ${\mathfrak U} $ intersects
all the symplectic leaves included of $\bar{\mathcal S} $.

Let $P$ be the connected parabolic subgroup that integrates ${\mathfrak P} $.
Since $G/P $ is projective, the projection $\Pi_2  $  onto the second component
 $ G/P \times {\mathfrak g}  \to  {\mathfrak g}$ is a closed map.
Now $G/P $ is the Grassmannian of all Lie subalgebras of  ${\mathfrak
  g}$  conjugate to ${\mathfrak P} $.
The set $S   $ of pairs $(\tilde{\mathfrak P} ,x  ) \in G/P \times {\mathfrak g} $ such that  $x $
belongs to the nilradical of $\tilde{\mathfrak P} $
 is a closed subset. Since $ \Pi_2$ is a closed map, $\Pi_2(S) $ is a
 closed subset of ${\mathfrak g} $.
But $\Pi_2(S) $ is precisely the union of all the adjoint orbits through
${\mathfrak U} $. This completes the proof.
\end{proof}

We have an identification ${\mathfrak P} \simeq {\mathfrak U}^{\perp}   $. As a
consequence, the sub-Lie groupoids of $G \times {\mathfrak g} \toto
{\mathfrak g}$ that integrate the subalgebroid $ {\mathfrak P}
\times {\mathfrak U} \to {\mathfrak U} $ are all the Lie groupoids
of the form $P \times {\mathfrak U} \toto {\mathfrak U} $ where $P$
is any parabolic Lie subgroup of $G$ with Lie algebra ${\mathfrak
P}$. In this case, one can identify $\Gamma_{\mathfrak U}$ with $G
\times {\mathfrak U} $ and the left action of $P \times {\mathfrak
U} \toto {\mathfrak U} $ corresponds precisely to the diagonal
action of $P$ given by
  $$ p \cdot (g,u) = (pg, Ad_{p} u) \ \ \ \ \forall g \in G, u \in {\mathfrak U}, p \in P .$$
so that $Z_{\mathfrak U}= \frac{G \times {\mathfrak U} }{P} $ and $ \phi (\overline{(g,u)}) = Ad_{g^{-1} } u $
where $\overline{(g,u)}$ stands for the class of $(g,u) \in G \times {\mathfrak U}$.

Since ${\mathfrak U}$ is a connected set, $ \frac{G \times {\mathfrak U}}{P},\phi)$ is
an \'etale symplectic resolution with typical fiber $\frac{
  \pi_0(Stab_G(x))}{\pi_0(Stab_P(x))}$. According to Theorem
\ref{theo:sympl_resol}, it is a symplectic resolution if and only if
$$ \pi_0(Stab_P(x)) =  \pi_0(Stab_G(x)).$$
We recover as a particular case of Theorem
\ref{theo:sympl_resol} the following Proposition,
which seems well-known (see
 Proposition 3.15 in \cite{Fu} for instance).

\begin{prop}\cite{Fu}
The Springer resolution  $(\frac{G \times {\mathfrak U}}{P} , \overline{(g,u)} \to
Ad_{g^{-1} u } )$ of the closure $\bar{\mathcal S} $ of a Richardson
orbit ${\mathcal S} $  is a covering symplectic resolution for any
parabolic subgroup $P$ such that ${\mathfrak U} \cap {\mathcal S}$ is
dense in $ {\mathcal S}$, where ${\mathfrak U} $ is the nilradical of
$Lie(P) $. It is a
symplectic resolution if and only if there exists a such a parabolic
subgroup $P$  such that $
\pi_0(Stab_P(x)) =  \pi_0(Stab_G(x))$.
\end{prop}

\begin{rmk}
For any $g \in G $, and any ${\mathfrak U}$ and $P$ as in Proposition
\cite{Fu} above, ${\rm Ad}_g {\mathfrak U} $
is again the nilradical of the Lie algebra of the parabolic subgroup $g P g^{-1}$.
In particular, ${\rm Ad}_g {\mathfrak U} $ is again a Lagrangian crossing of $\overline{\mathcal S} $
and ${\rm Ad}_g {\mathfrak U} $ is a Lie groupoid that integrates it.
In particular, one can form a second symplectic resolutions with $Z':=\frac{G \times {\rm Ad}_g ({\mathfrak U})}{ gP g^{-1} } $ 
and $\phi': Z' \mapsto \bar{\mathcal S} $  defined as before.
This second symplectic resolution is isomorphic to the first one.
The  Lagrangian closed submanifold $I$ of $\Gamma $ that gives, according to
Proposition \ref{prop:repr_isom2}, the  Morita equivalence between the Lie groupoids
$P \toto  {\mathfrak U}$ and $gP g^{-1} \toto {\rm Ad}_g {\mathfrak U} $ is 
$$ I:= g P \times {\mathfrak U} \subset G \times {\mathfrak g}.$$
\end{rmk}

\subsection{Exact multiplicative $k$-vector fields}

For any algebroid $(A \to M,\rho,[\cdot ,\cdot ]) $, and any Lie groupoid
$\Gamma \toto M $ that integrates it, and any $\Lambda \in \Gamma(\wedge^k A
\to M) $, the $k$-vector field $ \overrightarrow{\Lambda} - \overleftarrow{\Lambda}  $ is multiplicative.
The $k$-vector fields it defines on $M$ is simply $\rho (\Lambda) $,
and is tangent to all algebroid leaves.

We choose a locally closed algebroid leaf ${\mathcal S} $ and an algebroid crossing $L $ of 
$\overline{\mathcal S} $ coisotropic with respect to $\rho(\Lambda) $.

We leave it to the reader to check that the $k$-vector field on $\pi_{Z(R)}$
induced in this case on the \'etale resolution $Z(R)= R \backslash \Gamma_L $ (provided that it exists)
 is simply the infinitesimal $k$-vector field associated to $\Lambda $ through the right action of $\Gamma
\toto M $ on $Z(R) $.

\subsection{The Grothendieck resolution with its Evens-Lu Poisson structure.}

We show that the Grothendieck resolution is an example of resolution of a (holomorphic) algebroid leaf, and we use Poisson  groupoids
  to turn the Grothendieck resolution in a Poisson resolution.
This Poisson structure is precisely the one discovered by Sam Evens and Jiang-Hua Lu in \cite{LuEvens}.
Any Lie group $G$ acts on itself by conjugation, so that one can form the
action   groupoid $ G \times \underline{G} \toto \underline{G} $
(here $\underline{G} $ stands for the Lie group $G$ when it can be considered as a manifold acted
upon by conjugation, while we keep the notation  $G$ when $G$ is considered as a Lie group).
 Recall that the source map $s$ is given by $(g,h) \to h $, the target map $t$ is given by $(g,h) \to
g^{-1}h g $ and the product is given by $ (g_1,h_1)\cdot (g_2,h_2) = (g_1 g_2,
h_1)$ whenever $h_2 = t(g_1,h_1) $.

\smallskip
{\bf The Grothendieck resolution.}
\smallskip

Assume now that $G$ is a complex simple, connected and simply-connected  Lie group. Let $H$ be a Cartan subgroup,
$B$ a Borel subgroup containing $H$, and $U$ the unipotent radical of $B$.
We choose some $t \in H$ and denote by ${\mathcal S} $
the regular orbit (= conjugacy class) containing  $t$ in its closure,
see \cite{Hum}. The closure $\overline{\mathcal S}  $ of $ {\mathcal S}$ is what is called a {\em Steinberg fiber}
and denoted  by $F_t$ in \cite{LuEvens}.

\begin{lem}\label{lem:gro}
The submanifold $L=tU$ is an algebroid crossing of $\overline{\mathcal S}$
with normalisation ${\rm Lie }(R)  \to tU$.
\end{lem}
\begin{proof}
Any $g \in G$ belongs to a Borel subgroup, and all Borel subgroups are
conjugate. In particular, any $g \in \overline{\mathcal S} $ is conjugate
to an element in $t'U$, with $t' \in H$. Since the intersection of  two Steinberg fibers is empty
(see section 3.2 in \cite{LuEvens}),
one needs to have $t=t'$, and $tU $ intersects all the algebroid
leaves contained in $\overline{\mathcal S}$.

The Lie groupoid $R = B \times tU \toto tU $ is a closed sub-Lie groupoid
of the Lie groupoid $ \Gamma = G \times \underline{G} \toto
\underline{G}$, and:
\begin{equation}\label{eq:maxim}   \Gamma_{L \cap {\mathcal S}}^{L \cap {\mathcal S}}  \subset R .\end{equation}
Let ${\rm Lie }(R)  \to tU$ be the Lie algebroid of $ R \toto tU $.
For all $g \in  tU \cap {\mathcal S}$,  Eq. (\ref{eq:maxim}) amounts
to the fact that  ${\rm Lie} (R ) = \rho^{-1}(T_g tU) $, so that $     {\rm Lie} (R )
\to t U$ is the normaization of $tU$.
By construction $  L \cap {\mathcal S}$ is a dense open subset of $L=tU$. This completes the proof.
\end{proof}
By Lemma \ref{lem:gro} and (\ref{eq:maxim}), all the assumptions of Proposition \ref{prop:desing+groupoid}(3) are satisfied,
and one can construct a resolution of $\overline{\mathcal S}$, that we denote by $(X_t, \mu)$
in order to follow the notations of \cite{LuEvens} again. Let us give explicitly the construction of $ X_t$ and $\mu$.
To start with, we clearly have $\Gamma_L = G \times L = G \times tU $, while the quotient under the
left action of the Lie sub-groupoid $R \toto L  $ is the quotient  of $G \times tU$
through the action of the Lie group $B$ given by
 \begin{equation}\label{eq:baction} b\cdot (g, tu )   =(bg, btub^{-1}  ) \, \, \, \, \forall b \in B, g \in G, u \in U .\end{equation}
The map $\mu$ is the map  $\mu([g,tu]) \to Ad_g (tu) $ (where $[g,tu]  \in \frac{G \times {L}}{B} $ is the class
modulo the action of $B$ of the element $(g,tu) \in G \times L $).  

A short comparison with Section 3.2 in \cite{LuEvens}
shows that the resolution $(X_t, \mu)$  coincides with the Grothendieck resolution (also called Springer resolution).
In conclusion, we have the following proposition.

\begin{prop}\label{prop:gro}
Let $G$ be a simple, connected and simply-connected  complex Lie group, $H $ a Cartan
subgroup and $t \in H$. 
Let  ${\mathcal S} \subset G$ be a regular orbit with $t \in \overline{\mathcal S}$,
and $F_t = \overline{\mathcal S}$ be the Steinberg fiber. Let $(X_t,\mu)$ be as above.

The resolution $(X_t, \mu) $  is a resolution of  the Steinberg fiber $F_t = \overline{\mathcal S}$.
\end{prop} 
\smallskip
{\bf The Evens-Lu Poisson structure.}
\smallskip
Now, we endow the resolution $(X_t,\mu)$ of $F_t$  with a Poisson structure, following \cite{LuEvens} as a guideline again. The
construction below matches step by step the construction in  \cite{LuEvens},  only the interpretation claims to be new.
Let ${\mathfrak g} $ be the Lie algebra of $G$, and $  $.
The Lie algebra ${\mathfrak h}$ of $H$ is a Cartan subalgebra, and we can
choose a root decomposition $\Phi= \Phi^+ \cup \Phi^-$ and root vectors
$E_\alpha, E_{-\alpha}, \alpha \in \Phi_+ $ so that the space  
 $$ {\mathfrak n}= \oplus_{\alpha \in \Phi^+ } {\mathbb C} E_\alpha  $$
is the Lie algebra of $U$. Also we define  
 $$ {\mathfrak n}_- = \oplus_{\alpha \in \Phi^+ } {\mathbb C} E_{-\alpha}  .$$ 
We recall the construction of the standard Manin triple \cite{Korogodski}.
 Let $ {\mathfrak g}_\Delta $ be the diagonal of ${\mathfrak g} \otimes {\mathfrak g} $, and 
$$ {\mathfrak g}^*_{st} = \{(x+y, -y+z) \, | \, x \in {\mathfrak n}, z \in {\mathfrak n}_- , y \in {\mathfrak h} \} .$$ 
 Then $( {\mathfrak g} \oplus {\mathfrak g},{\mathfrak g}_\Delta,
{\mathfrak g}_{st}^* ) $ is a Manin triple.

According to \cite{ILX}, Section 4.5 (in particular Theorem 4.21), to any Manin triple $ ( {\mathfrak d},{\mathfrak g}_1, {\mathfrak g}_2 ) $,
is associated a natural multiplicative Poisson structure on the transformation groupoid $G_1 \times D/G_1 \toto D/G_1 $, where $D$ and $G_1 \subset D$ are  connected and simply-connected  Lie groups integrating $ {\mathfrak d}$ and ${\mathfrak g}_1$ respectively, and where $G_1$ acts on $D/G_1 $ by left multiplication.
In the present case, $D =G \times G $  and $G_1 $ is the diagonal $G_\Delta$ of $G \times G$, so that
$ D/G_1$ can be identified with $ G$ by mapping  $[g_1,g_2] \in D$ to $g_1 g_2^{-1} \in G$,
where $[g_1,g_2] $ stands for the class of $(g_1,g_2) $ modulo the action of $B$ given
by Equation (\ref{eq:baction}). 
Under this isomorphism the $G$-action on $D/G $ becomes the  conjugation of $G$ on $\underline{G} $,
 so that the Lie groupoid $G \times D/G \toto D/G$ can be identified with the Lie groupoid  $ G \times \underline{G} \toto \underline{G}$
 previously described.   
In conclusion, the Lie groupoid  $ G \times \underline{G} \toto \underline{G}$ can be endowed
with a Poisson structure $\pi_{G \times  \underline{G}} $ that turns it into a Poisson groupoid
(i.e. a Lie groupoid endowed with a multiplicative Poisson bivector field, see \cite{Xu:1995}).
A precise description of $\pi_{G \times  \underline{G}} $ is given by Equation (77)  in \cite{ILX}.
Since the pair $(G \times \underline{G} \toto \underline{G},\pi_{G \times  \underline{G}}  )$ is a Poisson Lie group,
 there exists a Poisson structure $\pi_{\underline{G}}$ on the base manifold $\underline{G}$ such that
 $ \pi_{\underline{G}}= s_*  (\pi_{G \times  \underline{G}}) $ as in Equation (\ref{eq:proje}).
We wish to compare these vector fields with the Poisson structure  on  $G \times G $  called  $\pi_D^+$
and the Poisson structure on $G$ called $\pi $ introduced  in \cite{LuEvens}, Section 2.2:

\begin{lem}\label{lem:comp:EL}
\begin{enumerate}
\item
The Poisson structures $\pi$ and $-\pi_{\underline{G}} $ coincide.
\item
The Poisson structure $\pi_{G \times  \underline{G}}$ is the image 
of the Poisson structure $-\pi_D^+$ through the map 
$(g_1, g_2) \to (g_2 , g_1 g_2^{-1}) $
\end{enumerate}
\end{lem}
\begin{proof}
We only prove the first point, since it is the only point really needed to prove Proposition \ref{prop:EL}.
The second one is just a cumbersome computation.

According to Proposition 4.19 in \cite{ILX}, the Poisson structure induced on $ \underline{G}\simeq G \times G / G_\Delta $ (which is the space denoted by $S$ in \cite{ILX}) is equal to:
  \begin{equation}\label{eq:davidencore} \pi_{\underline{G}}= - \sum_{i=1}^d (e_i)_S \wedge (\epsilon_i)_S \end{equation}
where $ X_S$ stands for the infinitesimal vector field induced by the action of 
an element $X \in  {\mathfrak d} = {\mathfrak g} \oplus {\mathfrak g}$ on $\underline{G}\simeq G\times G / G_\Delta $, and where $(e_i)_{i=1}^d, (\epsilon_i)_{i=1}^d $
are dual bases of ${\mathfrak g}_{st}^* $ and ${\mathfrak g}_\Delta$ respectively.
But, under the identification $G \times G / G_\Delta \simeq G $ mapping $[g_1, g_2]$ to 
$g_1^{-1}g_2 $, as previously described, we have $(X,Y)_S = X^L- Y^R $, where $X^L $ and $Y^R$
are the left and right actions of $X , Y \in {\mathfrak g}$ on $G$ respectively.
The construction  of $\pi_{\underline{G}} $ given in Eq. (\ref{eq:davidencore}) 
coincides then with the construction of $\pi$ described in Equation (2.10) in \cite{LuEvens}, up to a sign.
 \end{proof}
\begin{rmk}
To find again the explicit form of $ \pi$  given in Equation (2.10) in \cite{LuEvens},
one can choose the dual bases  of ${\mathfrak g}_{st}^* $ and ${\mathfrak g}_\Delta$ given by
 $$ \left\{(y_1,-y_1), \dots, (y_r,-y_r), (0, -E_{-\alpha}), (E_{\alpha}, 0 ) \, | \alpha \in \Phi^+ \right\} $$
 and
 $$ \left\{(y_1,y_1), \dots, (y_r,y_r), (E_\alpha, E_\alpha), (E_{-\alpha},E_{-\alpha}) \, | \alpha \in \Phi^+ \right\} .$$
 where $(y_i)_{i=1}^r $ is a base of  ${\mathfrak h} $ with $2\left\langle y_i,y_j\right\rangle  =  \delta_i^j$,
 and where we assume that $ \left\langle E_\alpha, E_{-\alpha} \right\rangle  =1 $ for all $\alpha \in \Phi^+ $.
 In the previous,  $ \left\langle \cdot, \cdot \right\rangle$ stands of course for the Killing form.
\end{rmk}

  %
 %
 %
 %
 %
%

%
According to Lemma 3.7 in \cite{LuEvens}, the submanifold $L = tU $
is coisotropic with respect to  $\pi_{\underline{G}} $ (which, according to Lemma \ref{lem:comp:EL}(1), coincides with the Poisson structure denoted by $\pi$ is \cite{LuEvens}).
All the conditions of Theorem  \ref{theo:compatible_resolutions}(4)  are therefore satisfied. 
We can  construct a Poisson structure on $ X_t = \frac{G \times L}{B}$ 
such that $\mu: X_t = \frac{G \times L}{B} \to \underline{G} $  is a Poisson map, id. est. $(X_t,\mu)$  is a resolution
compatible with the Poisson structure $ \pi_{ \underline{G} } $. In
conclusion, we have proved the follwing proposition:

\begin{prop}\label{prop:EL}  
Let $G$ be a simple connected and simply-connected  complex Lie group, and $H $ a Cartan
subgroup. Let  ${\mathcal S} \subset G$ be a regular orbit with $t \in \overline{\mathcal S}$.

The Poisson structure $\pi_{\underline{G}} $ is tangent to the Steinberg fiber $F_t = \overline{\mathcal S}$,
and the resolution $(X_t, \mu) $ of $F_t$ is a Poisson resolution.
\end{prop}

Since $\mu^{-1}({\mathcal S})$ is open and dense in $X_t$, the Poisson structure on $X_t$ such that $\mu$
is a Poisson map is unique, and the one constructed here needs to coincide with the one constructed
in \cite{LuEvens}. 
Proposition \ref{prop:EL} therefore  reproves Proposition 4.5 (2) is \cite{LuEvens}.
Indeed,  a step by step comparison, with the help of Lemma \ref{lem:comp:EL}(2), shows that the present construction
of the Poisson structure on $X_t$ 
as exposed in the proof of Theorem  \ref{theo:compatible_resolutions}(4), coincides precisely 
with the construction described in \cite{LuEvens}.
\subsection{Minimal resolutions of $ {\mathbb C}^2/ ({\mathbb Z}/ l {\mathbb Z}) $.}

    Let $l \in {\mathbb N}^*$ be an integer.
The group ${\mathbb Z}_l  ={\mathbb Z}/ l {\mathbb Z}$, seen as the group of $l^{th}$ roots of the unity,
acts  on $ {\mathbb C}^2$ by $\lambda \cdot (z_1, z_2) =(\lambda z_1, \overline{\lambda} z_2) $,
for all $\lambda \in {\mathbb Z}_l$.
The quotient space $W_l={\mathbb C}^2 / {\mathbb Z}_l  $ is an affine variety
that can be described as the zero locus in ${\mathbb C}^3$ of the function 
      $$ \chi_l (x,y,z) := xy - z^{l} .$$
The variety $W_l$ has only one singular point $O$; when seen as a subvariety
of ${\mathbb C}^3$, this singular point is the origin.

Let us recall several facts about its canonical Poisson structure.
The action of $G$ preserves the canonical symplectic structure on $ {\mathbb C}^2$, 
hence the canonical Poisson bracket on $  {\mathbb C}^2$ goes to the quotient and induces a  Poisson bracket
$  \{\cdot, \cdot\}_{W_l}$  on $W_l$. This Poisson structure $\pi_{W_l} $ is symplectic at all regular points of $W_l$.

Alternatively, it can be described as follows.
Consider the following Poisson bracket $ \pi_{{\mathbb C}^3}= \{ \cdot, \cdot \}_{{\mathbb C}^3}$ on $ {\mathbb C}^3$:
 \begin{equation}\label{eq:c3} \{x,y\}_{{\mathbb C}^3}=  \frac{\partial \chi_l}{\partial z} , 
 \, \{y,z\}_{{\mathbb C}^3}=  \frac{\partial \chi_l}{\partial x} , 
 \,  \{z,x\}_{{\mathbb C}^3}= \frac{\partial \chi_l}{\partial y} .\end{equation}
 Then $\chi_l$ is a Casimir function of the latest bracket, so that 
the Poisson structure it defines induces a Poisson structure on the zero locus $W_l $ of $\chi_l$.
 
 \begin{lem}
There exists  a algebraic Poisson variety $(N,\pi_N)$ such that
\begin{enumerate}
\item  $N$ is a nonsingular variety.
\item $(N, \pi_N)$ is integrable (when seen as a holomorphic Poisson manifold).
\item $(N,\pi_N)$ admits a symplectic leaf ${\mathcal S} $ whose closure $ \overline{\mathcal S}$
 is a subvariety of $N$ isomorphic to $W_l$ as an algebraic Poisson manifold. moreover,
 ${\mathcal S} $ is the regular part of $ \overline{\mathcal S}$.
 \end{enumerate} 
 \end{lem}
 \begin{proof}
 A natural candidate for $(N, \pi_N)$ would be $({\mathbb C}^3, \pi_{{\mathbb C}^3})$.
 But it is not clear that this structure is integrable. But $({\mathbb C}^3, \pi_{{\mathbb C}^3})$
 is, according to Theorem 5.5 in \cite{DSV}, a Poisson submanifold of a Poisson manifold
 $N$ (and also denoted $N$ in \cite{DSV}) that we now describe. 
 
Let $e \in {\rm sl}_l({\mathbb C}) $ be an element of the subregular nilpotent orbit,
and  ${\mathfrak n} \subset {\rm sl}_l({\mathbb C}) $ be a complement of the centraliser of $e$.
Then the affine space $N:= x+ {\mathfrak n}^{\perp} $  is of course a nonsingular submanifold, so that condition 1) is satisfied.

According to \cite{Sabourin}, one can choose ${\mathfrak n}$  so that 
the Poisson matrix of the linear Poisson structure of ${\rm sl}_{l}({\mathbb C}) $ is, for all $y \in N$, of the form 
  \begin{equation} \label{sabourin} \left( \begin{array}{cc} A(y) & B(y) \\-B^{\perp}(y)  &  C(y) \\ \end{array} \right)  \end{equation}
where $ C(y)$ is an invertible matrix (more precisely, a matrix of determinant $1$, see the proof of Theorem 2.3 in
 \cite{Sabourin}). 
 This amounts to the fact that $N$ is a Dirac submanifold  of ${\rm sl}_l({\mathbb C})  $.
 The classical procedure called Dirac reduction  in \cite{Sabourin} (and sometimes called Poisson-Dirac reduction)
 yields then a Poisson structure $\pi_N$ on $N $, which is polynomial by construction. 
 In conclusion $(N,\pi_N)$ is an algebraic Poisson variety.
 
 Let us show that it is integrable. First, $ {\rm sl}_l ({\mathbb C})$ is an integrable Poisson manifold,
 and the symplectic Lie groupoid that integrates this Poisson manifold
 is $ \Gamma = {\rm SL}_l({\mathbb C}) \times  {\rm sl}_l ({\mathbb
 C}) \toto  {\rm sl}_l ({\mathbb
 C})$. 
  Since the matrix $C(y)$ is invertible for all $y \in N$, $N$ is a
 cosymplectic submanifold of ${\rm sl}_l({\mathbb C})  $ (see
 \cite{CF2}, a cosymplectic manifold is a Poisson-Dirac manifold with
 a symplectic transverse structure). This implies that 
 $\Gamma_N^N  \toto N $ is a symplectic sub-Lie groupoid  of $\Gamma \toto  {\rm sl}_l ({\mathbb C})$,
 and the Poisson structure it integrates is the Dirac-Poisson structure on $N$, see \cite{CF2}. 
Hence $N$ is a Poisson-Dirac submanifold.
%
\end{proof}

For all $k \in \{1,\cdots, l-1\} $, the subvariety of ${\mathbb C}^3$
   $$ L_k := \{ (\lambda^k, \lambda^{l-k}, \lambda )  , \lambda \in {\mathbb C}   \}$$
is a subvariety of $W_l $. It can  also be described by the equations:
   \begin{equation}\label{eq:lkn} x= z^k \mbox{ and } y = z^{l-k}. \end{equation}

\begin{lem}
          For all $k \in \{1,\cdots, n-1\} $, the subvariety      $ L_k$
is a Lagrangian crossing of~${\overline{\mathcal S}} $. 
\end{lem}
\begin{proof}
The subvariety $L_k$ is nonsingular, since it is given by (\ref{eq:lkn}) when seen 
as a subvariety of ${\mathbb C}^3 $ via the inclusion $W_l \subset {\mathbb C}^3 $ described above. 
The intersection of $L_k $ with $ {\mathcal S}$ is of course Lagrangian since it has dimension $1$. 
Now, the identity $\overline{\mathcal S} = {\mathcal S} \cup \{O\} $ holds,  the subvariety $ L_k$ contains the point $O$,
and has a nonempty intersection with ${\mathcal S} $. 
All the assumptions are therefore satisfied.
\end{proof}

\vspace{0.5cm}

 The variety $W_l$  admits a minimal resolution $(\Sigma, \phi) $, see \cite{Sha} Section IV-4-3, and it is well-known
 that this resolution is symplectic, see Example 2.2 in \cite{FuSum}.
 
 Let us say a few words on this resolution.
 
 Explicitly, $(\Sigma, \phi)$ is constructed with a help of $l-1$ successive blowup, as follows.
 The blowup of $ W_l$ at $0$ is the projective variety of $P^3 ({\mathbb C})$ given, in the three canonical charts, by the equations
  $$   x = z^l y^{l-2} , \hspace{0.5cm} y = z^l x^{l-2}  , \hspace{0.5cm} \psi_{l-2}(x,y,z)=0 $$ 
 with the usual gluing relations. It is easy to check that the inverse image of $0$ consists of two copies of
 $P^1({\mathbb C}) $ that intersect transversally at a point if $l \neq 2$, and consists of one copy of $P^1({\mathbb C}) $ if $l=2$.
 Also, the two first components of the blowup $\tilde{W} $ are nonsingular, while the last one is isomorphic to $W_{l-2} $,
 to which the procedure can be applied recursively until the nonsingul	ar varieties $W_0$ or $W_1$ appear.  
 Applying successive blowup therefore, one gets a resolution $(\Sigma, \phi) $ of $W_l$ and a closed look at the construction amounts to the following properties.   
      
 \begin{lem}   \label{lem:techlk}   
 \begin{enumerate}
 \item $(\Sigma, \phi) $ is the minimal resolution.
 \item $(\Sigma, \phi) $ is a proper symplectic resolution which is compatible with the Lagrangian crossing
 $L_k $ for all $k=1,\dots, l-1$. Let $ \tilde{L_k}$ be the submanifold of $\Sigma $
 to which the restriction of $\tilde{L_k} $ is a biholomorphism onto $ L_k $. 
 \item The inverse image of $0$ consists of $l-1 $ projective curves $C_1, \cdots, C_{l-1} $, all
 isomorphic to $ { \mathbb P}^1({\mathbb C})$.
 \item   
 For all $k=1,\cdots, l-1 $, $\tilde{L_k}  $ intersects $C_k $ transversally
 at exactly one point, and $ \tilde{L_k} \cap C_i = \emptyset $ for $ k \neq i$.
 \item 
 For all $k=1,\cdots, l-2 $, $C_{k}  $ intersects $C_{k+1} $ at exactly one point $p_k $, and $ C_j \cap C_i = \emptyset $ for $| j -i| \geq 2 $.
 \item  The kernel of the differential of $\phi $ at a point $y \in C_k $ distinct from $p_k $ or $p_{k-1} $ is equal to 
 $ T_y C_k $.
 \item The differential of $ \phi$ vanishes at the points $p_1, \dots, p_{l-2}$.
 \end{enumerate}     
 \end{lem}
%
  
 
Let $\Gamma \toto M $ be the source-simply connected symplectic Lie groupoid that integrates $(M,\pi_M) $.
All the assumptions of Theorem \ref{theo:Linverse} are satisfied, so that there exists, for all $k=1, \dots, l-1 $,
a closed sub-Lie groupoid $R_k \toto L_k $ of $ \Gamma \toto M $ that integrates $ TL_k^{\perp} \to L_k $ and contains 
$ \Gamma_{L_k \cap {\mathcal S}}^{L_k \cap {\mathcal S}} $, and there exists an open subset $U_k $ of $\Sigma $
isomorphic to the symplectic resolution  associated to $R_k \toto L_k$
as in Theorem \ref{theo:sympl_resol} (3).
The following Proposition describes in a very explicit way this open subset:
\begin{prop}
For all $k=1, \cdots, l-1$, we have
  $$U_k = \Sigma - \cup_{i \neq k} C_i $$  
\end{prop}
\begin{proof}
We assume $k \neq 1$ and $k \neq l-1$ for simplicity.
The cases $k=1  $ and $k = l-1 $ can be dealt in the same way. 
%
%
%
Let $O_k \in \tilde{L_k}$ be the inverse image of $O$ through the restriction of
$\phi $ to a biholomorphism from $\tilde{L_k} $ onto $L_k $. Notice that $O_k $ belongs to $C_k$, according to 
Lemma \ref{lem:techlk}(4).
By construction, $ U_k$ is equal to
 $$ \Gamma \cdot (\tilde{L_k} \backslash \{O_k\} ) \cup \Gamma \cdot \{O_k\} ,$$
where $ \Gamma \cdot \{\tilde{L_k} \backslash O_k \} $ 
(resp. $\Gamma \cdot \{O_k\} $) stands for the $\Gamma$-orbit of $\tilde{L_k} \backslash \{O_k\} $ (resp. $O_k$)
with respect to the action of $\Gamma \toto M$ on $\Sigma$ defined in 
 Proposition \ref{prop:cfCF2}. Proposition \ref{prop:cfCF2} gives
 the identity $ \Gamma \cdot (\tilde{L_k} \backslash \{O_k\} ) =  \phi^{-1}(\mathcal S)$, since 
 the inclusion $ \tilde{L_k} \backslash \{O_k\} \subset \phi^{-1}({\mathcal S})$ holds by Lemma \ref{lem:techlk}(4). 
 We therefore have 
 $$ U_k =\phi^{-1}(\mathcal S) \cup  \Gamma \cdot \{O_k\}.$$ 
Now,  Proposition \ref{prop:cfCF2}, together with Equation (\ref{eq:intcotg})
give that a point $y \in \phi^{-1}(O) $ belongs to $U_k$ if and only if there exists a  path
$a(t) $ in $T^*_O M $ together with a smooth path  $\sigma (t) $ in $\Sigma $ (with $t \in [0,1]$) such that, first, $\sigma (1) = x  $ 
and $ \sigma (0) =O_k $, and second:
 \begin{equation}\label{eqn:sigma}   \frac{\diff \sigma(t)}{\diff t}  =  (\pi_\Sigma^\#)_{|_{\sigma (t)}}   \big( a (t) \circ   \diff_{\sigma (t)} \phi \big)  . \end{equation} 
Assume that  $y \in \cup_{i \neq k} C_i$, then $\sigma (t) $  has to go through one at least of the points $p_{k-1}$ or $p_k$ for some $t = t_0$. But then, there can not exist a smooth path $a(t) $ that satisfies (\ref{eqn:sigma}) since the right hand side of (\ref{eqn:sigma})
would vanish at $t =t_0 $ by Lemma \ref{lem:techlk}(7), so that $\sigma (t)$ would have to be a constant path.
Hence, $ y$ has to be a point in $C_k$.
For all $y \in C_k$, there exists  a smooth path  $\sigma (t) $, taking values in $ C_k \backslash \{p_k, p_{k+1}\}$, 
such that  $ \sigma (0) = O_k, \sigma (1) =y$.
Since the image of the dual map of $ \diff_{\sigma (t)} \phi $ has rank $1$ for all $t$ by Lemma \ref{lem:techlk}(6), it contains
the covector  $ \omega_\sigma (  \frac{\diff \sigma(t)}{\diff t}  , \cdot ) $.
There exists therefore a path $a(t) \in  T_O^* M $ which satisfies Eq. (\ref{eqn:sigma}), and $x \in \Gamma \cdot O_k$. this completes the proof.
\end{proof} 
\begin{rmk}
$ W_l$ is a particular type of Kleinian singularity. It is natural to ask whether these constructions could be done for other Kleinian singularities, those of type $D_l, l \geq 4$, or $E_6,E_7,E_8$. The answer is negative in general, due to the lack of Lagrangian crossing
in these cases.
\end{rmk}


\begin{thebibliography}{99}

\bibitem{Be2} A. Beauville,
 Symplectic singularities, {\em Invent. Math.} {\bf 139}, (2000) 541-549.

\bibitem{CDW}
A. Coste, P. Dazord and A. Weinstein,
Groupo{\"\i}des symplectiques, Publications du D{\'e}partement de
Math{\'e}matiques de l'Universit{\'e} de Lyon, {I}, {\bf 2/A} (1987) 1--65.

\bibitem{CF} M. Crainic and R.L. Fernandes, Integrability of
Lie brackets, {\em Ann. of Math.} {\bf 157} (2003) 575-620.

\bibitem{CF2} M. Crainic and R.L. Fernandes, Integrability of Poisson
brackets, {\em J. Differential Geometry}.

\bibitem{DSV} P. Damianou, H. Sabourin, P. Vanhaecke, Transverse Poisson structures
to adjoint orbits in semi-simple Lie algebras. To appear in Pacific Journal of Mathematics.


\bibitem{Dazord} F.A. Cuesta, P. Dazord, G. Hector, Sur
  l'int\'egration symplectique de la structure de Poisson singuli\`ere
$(x^2 + y^2)\frac{\partial}{  \partial x}\wedge \frac{\partial }{ \partial
  y} $
de ${\mathbb R}^2 $. Publicacions Matem\`atiques, Vol {\bf 33} (1989) 411-415.

\bibitem{DeSilvaWeinstein} A. Cannas da Silva; A. Weinstein,
 Geometric models for noncommutative algebras. Berkeley Mathematics
 Lecture Notes, AMS, {\bf 10} (1999).

\bibitem{DZ} J.P. Dufour; N.T. Zung, 
Poisson structures and their normal forms.
{\em Progress in Mathematics}, {\bf 242}, (2005).

\bibitem{FuSum}  B. Fu, A survey on symplectic singularities and resolutions, {\em math.AG/051034.}

\bibitem{Fu} B. Fu, Symplectic Resolutions for Nilpotent Orbits,
{\em Invent. Math.} {\bf 151} (2003), 167-186.

\bibitem{LuEvens} S. Evens, J.H. Lu, 
 Poisson geometry of the Grothendieck resolution of a complex semisimple group,
   math. QA/0610123.

\bibitem{GK} V. Ginzburg; D. Kaledin, Poisson deformations of
  symplectic quotient singularities.  
{\em Adv. Math.}  {\bf 186 } (2004),  1-57.

\bibitem{HS} M. Hilsum; G. Skandalis,
Stabilit� des $C\sp{*} $-alg�bres de feuilletages. 
{\em Ann. Inst. Fourier } {\bf 33} (1983),  201--208. 

\bibitem{Hum} J.E. Humphreys, Conjugacy classes in semisimple
 algebraic groups.
 Mathematical Surveys and Monographs, {\bf 43}. (1995).

\bibitem{ILX} D. Iglesias-Ponte, C. Laurent-Gengoux, P. Xu,
 Universal lifting theorem and quasi-Poisson groupoids, math.DG/0507396.

\bibitem{Korogodski} L. Korogodski, Y. Soibelman, Algebra of functions on quantum groups,
part I, AMS, {\em Mathematical surveys and monographs}, {\bf 56}, (1998).

\bibitem{SX} C. Laurent-Gengoux, M. Stienon, P. Xu, 
Holomorphic  Poisson structures and holomorphic groupoids, 
arXiv:0707.4253v1.

\bibitem{LTX1} C. Laurent-Gengoux, J.-L. Tu, P. Xu,
  Chern-Weil map for principal bundles over
groupoids,  {\em Math. Z}, {\bf 255}, (2007), 451-491.    




\bibitem{McKenzie} K.C.H. Mackenzie,
 General theory of Lie groupoids and Lie algebroids. {\em London
 Mathematical Society Lecture Note Series}, {\bf 213}, (2005).

\bibitem{OR} P. Ortega; T. Ratiu, Momentum map and Hamiltonian reduction.
Progress in Mathematics, {\bf 222} (2004). 

\bibitem{Sabourin} H. Sabourin, Sur la structure transverse \`a une orbite de  Lie nilpotente adjointe.
{\em Canad. J. Math.}, {\bf 57}, (2005), 468-506.


\bibitem{Sha} I.R. Shavarevich, Basic Algebraic Geometry, Springer-Verlag, (1974). 

\bibitem{TY} P. Tauvel, R. Yu, Lie algebras and algebraic
 groups. Springer Monographs in Mathematics.
 Springer-Verlag, Berlin, (2005).

\bibitem{Weinstein1} A. Weinstein, The integration problem
for complex Lie algebroids. math.DG/0601752.

\bibitem{X_BV} P. Xu,  Gerstenhaber algebras and BV-algebras in
  Poisson geometry.  
{\em Comm. Math. Phys.}  {\bf 200}  (1999),   545--560.

\bibitem{Xu:2004} P. Xu,  Momentum maps and Morita equivalence. {\em
  J. Differential Geom.} {\bf 67 }(2004), no. 2, 289--333. 

\bibitem{Xu:1995} P. Xu, On Poisson groupoids,
 {\em Internat. J.~Math.} {\bf 6} (1995) 101--124.

\end{thebibliography}
\end{document}